\documentclass[11pt, a4paper, oneside]{amsart}
\usepackage[left=25mm,right=25mm,top=25mm,bottom=25mm]{geometry}

\usepackage{amssymb,amsmath,amsfonts}
\usepackage{pgf, pgfplots,pgfplotstable}

\pgfplotsset{compat=newest,
colormap={xray}{
    rgb255(0cm)=(255,255,255);
    rgb255(1cm)=(0,0,0);},
colormap={netgen}{
        rgb255(-1cm)=(0,0,255);
        rgb255(-0.75cm)=(0,143,255);
        rgb255(-0.5cm)=(0,255,216);
        rgb255(-0.25cm)=(0,255,71);
        rgb255(0cm)=(71,255,0);
        rgb255(0.25cm)=(216,255,0);
        rgb255(0.5cm)=(255,143,0);
        rgb255(0.75cm)=(255,0,0);}}

\usepgfplotslibrary{groupplots, polar}
\usetikzlibrary{shapes,arrows,decorations,calendar,matrix,backgrounds,folding,calc,positioning,spy, arrows.meta, external, bending}

\tikzset{external/only named=true}

\definecolor{matblue}{rgb}{0 0.4470 0.7410}
\definecolor{matorange}{rgb}{0.8500 0.3250 0.0980}
\definecolor{matyellow}{rgb}{0.9290 0.6940 0.1250}
\definecolor{matpurple}{rgb}{0.4940 0.1840 0.5560}
\definecolor{matgreen}{rgb}{0.4660 0.6740 0.1880}
\definecolor{matazure}{rgb}{0.3010 0.7450 0.9330}
\definecolor{matred}{rgb}{0.6350 0.0780 0.1840}
\definecolor{amber}{rgb}{1.0, 0.75, 0.0}

\usepackage[numbers]{natbib}

\newtheorem{thm}{Theorem}[section]
\theoremstyle{remark}
\newtheorem{remark}[thm]{Remark}

\usepackage{tabularx, booktabs, hyperref}
\newcolumntype{C}{>{\centering\arraybackslash}X}
\newcolumntype{L}{>{\raggedright\arraybackslash}X}
\newcolumntype{s}{>{\hsize=.5\hsize}C}
\newcolumntype{T}{>{\hsize=.05\hsize\raggedright\arraybackslash}X}
\newcolumntype{M}{>{\hsize=.25\hsize}C}

\usepackage{lmodern, float, bm, caption, url}

\renewcommand{\hat}{\widehat}
\renewcommand{\tilde}{\widetilde}

\newcommand{\pdt}[1]{{\frac{\partial #1}{\partial t}}}

\renewcommand{\div}{{\operatorname{div}}}

\newcommand{\diag}{{\operatorname{diag}}}
\newcommand{\tr}{{\operatorname{tr}}}

\newcommand{\grad}{{\nabla}}

\renewcommand{\Re}{{\operatorname{Re}}}
\renewcommand{\Pr}{{\operatorname{Pr}}}
\newcommand{\Ma}{{\operatorname{Ma}}}

\newcommand{\St}{{\operatorname{St}}}

\newcommand{\jump}[1]{{{[\![#1]\!]}}}

\newcommand{\facets}{\mathcal{F}_h}
\newcommand{\mesh}{\mathcal{T}_h}

\newcommand{\mat}[1]{{\bm{#1}}}
\renewcommand{\vec}[1]{{\bm{#1}}}

\newcommand{\I}{{\mat{I}}}
\newcommand{\T}{{\mathrm{T}}}

\newcommand{\VEL}{{\vec{u}}}
\newcommand{\HEAT}{{\vec{q}}}
\newcommand{\TAU}{{\mat{\tau}}}
\newcommand{\EPS}{{\mat{\varepsilon}}}
\newcommand{\EIG}{{\mat{\Lambda}}}
\newcommand{\AMP}{{\vec{\mathcal{L}}}}
\newcommand{\TRA}{{\vec{\mathcal{T}}}}
\newcommand{\VIS}{{\vec{\mathcal{V}}}}

\title[CBC for HDG]{Characteristic boundary conditions for Hybridizable Discontinuous Galerkin methods}

\author{Jan Ellmenreich}
\address{Institute of Analysis and Scientific Computing, TU Wien,
Wiedner Hauptstrasse 8-10, 1040 Wien}
\email{jan.ellmenreich@tuwien.ac.at}

\author{Philip L. Lederer}
\address{Fachbereich Mathematik, Universität Hamburg,
Bundesstraße 55, 20146 Hamburg, Germany}
\email{philip.lederer@uni-hamburg.de}

\author{Matteo Giacomini}
\address[UPC]{Laboratori de C\`alcul Num\`eric (LaC\`aN), ETS de Ingenier\'ia de Caminos, Canales y Puertos, Universitat Polit\`ecnica de Catalunya, Barcelona, Spain.}
\address[CIMNE]{Centre Internacional de M\`etodes Num\`erics en Enginyeria (CIMNE), Barcelona, Spain.}
\email{matteo.giacomini@upc.edu}

\author{Antonio Huerta}
\address[UPC]{Laboratori de C\`alcul Num\`eric (LaC\`aN), ETS de Ingenier\'ia de Caminos, Canales y Puertos, Universitat Polit\`ecnica de Catalunya, Barcelona, Spain.}
\address[CIMNE]{Centre Internacional de M\`etodes Num\`erics en Enginyeria (CIMNE), Barcelona, Spain.}
\email{antonio.huerta@upc.edu}

\begin{document}

\begin{abstract}
    In this work we introduce the concept of characteristic boundary conditions (CBCs) within the framework of Hybridizable Discontinuous Galerkin (HDG) methods,
    including both the Navier-Stokes characteristic boundary conditions (NSCBCs) and a novel approach to generalized characteristic relaxation boundary conditions (GRCBCs).
    CBCs are based on the characteristic decomposition of the compressible Euler equations
    and are designed to prevent the reflection of waves at the domain boundaries.
    We show the effectiveness of the proposed method for weakly compressible flows through a series of numerical experiments
    by comparing the results with common boundary conditions in the HDG setting and reference solutions available in the literature. In particular,
    HDG with CBCs show superior performance minimizing the reflection of vortices at artificial boundaries, for both inviscid and viscous flows.
\end{abstract}

\maketitle
\section{Introduction}

The direct numerical simulation (DNS) of unsteady compressible flows, involving the simultaneous resolution of both aerodynamic
and acoustic scales, commonly referred to as direct noise computation (DNC), is a challenging task in computational aeroacoustics (CAA). 
Accurately resolving the relevant physical scales necessitates the use of low-dissipative
and low-dispersive high-order numerical methods, particularly because acoustic waves can propagate over long
distances, with minimal attenuation due to viscous dissipation~\cite{coloniusComputationalAeroacousticsProgress2004}.
Furthermore, the high spatial resolution requirements call for compact and computationally efficient numerical methods,
such as finite-difference (FD), high-order finite-volume (FV), and the finite-element (FE) methods. Among these, the finite-element method is
particularly advantageous for handling arbitrary complex geometries using non-uniform, possibly unstructured, meshes. 

In simulations, limited computational resources necessitate a trade-off between spatial resolution and the
size of the computational domain, which inevitably results in domain truncation, especially in external flows. To approximate an otherwise infinite
domain, artificial boundary conditions, also known as far-field, transparent or non-reflection boundary conditions, are
employed to mimic the behavior of physical boundaries, and prevent spurious wave reflections from propagating back into the domain.
The development of accurate and efficient boundary conditions is, therefore, essential for the success of numerical simulations, 
particularly in the context of compressible flows. A comprehensive review of boundary conditions for compressible flows can be 
found in \cite{coloniusModelingArtificialBoundary2004, coloniusBoundaryConditionsTurbulence2023}.

A common strategy to prevent reflections at artificial boundaries is to extend the domain by adding a buffer layer.
Sponge layers \cite{bodonyAnalysisSpongeZones2006, maniAnalysisOptimizationNumerical2012} damp 
the flow exponentially over time towards a target state. However, these layers are inherently reflective \cite{coloniusModelingArtificialBoundary2004}, and
require to progressively increase the sponge strength to reduce reflections. More advanced techniques, such as perfectly matched layers (PML)
\cite{huAbsorbingBoundaryConditions2008, huDevelopmentPMLAbsorbing2008}, are designed to be non-reflecting, although their strict application is limited to 
linearized problems. Despite their effectiveness, all these methods suffer the common drawback of requiring additional computational resources.

A widely used method for imposing boundary conditions is the family of characteristic boundary conditions (CBCs), based on the characteristic 
decomposition of the non-linear Euler equations along the local normal to the boundary. This approach leads to the well-known locally 
one-dimensional and inviscid (LODI) relations \cite{hedstromNonreflectingBoundaryConditions1979a, thompsonTimeDependentBoundary1987,thompsonTimedependentBoundaryConditions1990}, 
which modify incoming characteristic amplitudes to prevent boundary reflections. The core concept is that, in a 
one-dimensional framework, wave amplitudes remain constant along characteristic lines, even in non-linear scenarios. To address the issue of the 
steady-state solutions depending on the initial conditions, \cite{rudyNonreflectingOutflowBoundary1980} 
introduced a relaxation term with user-defined parameters to accelerate convergence toward a target state. 
This relaxation method was further refined in a recent study \cite{pirozzoliGeneralizedCharacteristicRelaxation2013}, 
with the development of generalized characteristic relaxation boundary conditions (GRCBCs),
which eliminated the need for user-defined parameters by incorporating a suitable Courant-Friedrichs-Lewy condition (CFL) at the boundary.
The extension to viscous flows was accomplished by~\cite{poinsotBoundaryConditionsDirect}, 
who developed the Navier-Stokes characteristic boundary conditions (NSCBCs) through the introduction of appropriate viscous conditions. Since then,
NSCBCs have been used in numerous studies \cite{yooCharacteristicBoundaryConditions2005, lodatoConditionsAuxLimites2008, lodatoThreedimensionalBoundaryConditions2008, liuNonreflectingBoundaryConditions2010, odierCharacteristicInletBoundary2019}.

Although the aforementioned characteristic boundary conditions are inherently based on a one-dimensional framework, incorporating 
multidimensional effects is essential for accurately capturing complex flow phenomena. In \cite{yooCharacteristicBoundaryConditions2005}, the authors
extend NSCBCs by introducing damped convective transverse gradients through a transverse relaxation factor~$\beta$. 
The value of $\beta$ was determined via a low Mach number asymptotic analysis of the Euler equations. However, identifying the optimal value for $\beta$ remains
an active area of research \cite{lodatoConditionsAuxLimites2008,lodatoThreedimensionalBoundaryConditions2008,lodatoOptimalInclusionTransverse2012,liuNonreflectingBoundaryConditions2010}.

In the literature, the implementation of CBCs is often confined to FD methods, since the modified characteristic 
system is solely solved for the boundary nodes. However, in the finite element context, specifically for the Discontinuous Galerkin (DG) and 
Hybridizable Discontinuous Galerkin (HDG) methods, boundary conditions are enforced in a weak sense through numerical fluxes $\hat{\vec{F}}$ \cite{cockburnDiscontinuousGalerkinMethods2017, shuDiscontinuousGalerkinMethod2014, dipietroMathematicalAspectsDiscontinuous2012}. 
This approach requires an explicit external boundary state $\vec{U}^-$. To our knowledge, only two studies have addressed this issue within the DG framework.
In \cite{toulopoulosArtificialBoundaryConditions2011}, the authors
introduced mirror elements adjacent to the boundary elements, where the modified characteristic system is solved to determine $\vec{U}^-$.
In \cite{shehadiPolynomialcorrectionNavierStokesCharacteristic2024a,shehadiNonReflectingBoundaryConditions2024}, the authors imposed
the characteristic wave amplitudes and determined the external boundary state $\vec{U}^-$ via a least-squares approach by exploiting the discrete 
differential operator of the nodal shape functions.

HDG methods have gained significant popularity for simulating compressible flows in 
recent years \cite{peraireHybridizableDiscontinuousGalerkin2010, nguyenHybridizableDiscontinuousGalerkin2012, cockburnHDGMethodsHyperbolic2016a}. 
These methods effectively combine the benefits of DG approaches — such as local conservation, 
high-order spatial accuracy, and the capability to manage complex geometries — 
while offering potential reductions in computational costs, particularly in scenarios, involving implicit time stepping. This cost reduction is achieved by introducing appropriate traces of the unknowns on the mesh skeleton, 
significantly reducing the global coupling of degrees of freedom. This decoupling is especially beneficial as it facilitates the use of techniques 
like static condensation to eliminate the elemental degrees of freedom \cite{cockburnStaticCondensationHybridization2016}.

This work extends the concept of characteristic boundary conditions (CBCs) within the HDG framework,
encompassing both the Navier-Stokes characteristic boundary conditions (NSCBCs) and a novel approach to generalized characteristic relaxation boundary conditions (GRCBCs).

In Section~\ref{sec::compressible_equations} we give an overview of the compressible Navier-Stokes equations with particular emphasis on the derivation 
of the characteristic quasi-linear form of the Euler equations. We then introduce the fundamental concept of CBCs, including NSCBCs and GRCBCs, 
with a particular focus on their application to subsonic inflows and outflows.

In the subsequent Section~\ref{sec::hdg_for_compressible}, we introduce the finite element notation, the discrete variational formulation
of the mixed HDG method for compressible viscous flows, and some common boundary conditions. Additionally, we propose a slight modification of the far-field boundary condition, which showed numerical stability improvements in our experiments.

Then, in Section~\ref{sec::cbc_for_hdg}, we give a detailed explanation of CBCs implementation within the HDG framework, where 
we first consider the one-dimensional inviscid setting (LODI) and then extend it to the multidimensional and viscous setting. Additionally, we introduce a novel approach to GRCBCs
and show that common HDG boundary conditions can be recovered as a special case of the proposed method.
In contrast to the DG implementations, our implementation does not rely on additional mirror 
elements or nodal shape functions, and can be naturally embedded into the numerical method.

Section~\ref{sec::experiments} shows the effectiveness of the proposed method through a series of numerical experiments including 
a linear acoustic pulse, a non-linear oblique pressure pulse, a zero-circulation vortex, and a laminar viscous flow around a cylinder. 
We evaluate the performance of CBCs by comparing the results with the traditional boundary conditions and reference or exact solutions.
Our focus lies on weakly compressible flows, i.e. Mach number $\Ma_\infty \leq 0.3$, and the assessment of subsonic inflow and outflow boundary conditions.  

Ultimately, in Section~\ref{sec::conclusion} we summarize our findings.

\section{Compressible flow equations and characteristic boundary conditions}
\label{sec::compressible_equations}

The fundamental equations describing the motion of an unsteady, viscous compressible flow in a space-time cylinder
$\Omega \times (0, t_{end}] \in \mathbb{R}^{d+1}$ with non-empty, bounded, $d$-dimensional spatial domain $\Omega$, with
boundary $\partial \Omega$, and final time $t_{end}$, are specified by the compressible Navier-Stokes equations. In terms of
conservative variables, $\vec{U} = \begin{pmatrix} \rho, & \rho \vec{u}, & \rho E \end{pmatrix}^\T$, with density~$\rho$, velocity~$\vec{u}$, and total specific energy~$E$,
this system can be expressed in dimensionless form as
\begin{align}
    \label{eq::Navier-Stokes}
    \pdt{\vec{U}} + \div(\vec{F}(\vec{U}) - \vec{G}(\vec{U}, \EPS, \vec{\phi}) ) = \vec{0} \quad \text{in} \quad \Omega \times (0, t_{end}] ,
\end{align}
by applying either an aerodynamic or an acoustic scaling with the reference values
\begin{subequations}
    \begin{alignat}{4}
        \label{eq::aerodynamic-scaling}
         \text{(aerodynamic)} & \qquad & \rho_{\mathrm{ref}} & := \rho^*_\infty, \qquad & u_{\mathrm{ref}} & := |\vec{u}^*_\infty|, \qquad & \theta_{\mathrm{ref}} & := \theta^*_\infty (\gamma - 1) \Ma^2_\infty,  \\
        \label{eq::acoustic-scaling}
        \text{(acoustic)} & \qquad & \rho_{\mathrm{ref}} & := \rho^*_\infty, \qquad & u_{\mathrm{ref}} & := c^*_\infty, \qquad         & \theta_{\mathrm{ref}} & := \theta^*_\infty (\gamma - 1),
    \end{alignat}
\end{subequations}
depending on the flow regime under consideration, where $\theta$ and $c$ denote the temperature and speed of sound, respectively. Note that dimensional values are marked by an asterisk.
With the Mach number $\Ma_{\infty}$, the Reynolds number $\Re_{\infty}$ and the Prandtl number $\Pr_{\infty}$ defined as
\begin{align*}
    \Ma_{\infty} & := \frac{|\vec{u}^*_\infty|}{c^*_\infty}, & \Re_{\infty} & := \frac{u_{\mathrm{ref}} \rho^*_\infty L^*}{\mu^*_\infty}, & \Pr_{\infty} & := \frac{\mu^*_\infty c_p}{\kappa^*_\infty},
\end{align*}
the convective and the viscous fluxes are given by
\begin{align}
    \label{eq::fluxes}
    \vec{F}(\vec{U}) & = \begin{pmatrix}
                             \rho \vec{u}^\T                     \\[0.8ex]
                             \rho \vec{u} \otimes \vec{u} + p \I \\[0.8ex]
                             \rho H \vec{u}^\T
                         \end{pmatrix}, &
    \vec{G}(\vec{U}, \EPS, \vec{\phi}) = \frac{1}{\Re_{\infty}}\begin{pmatrix}
                                                                   \vec{0}^\T \\[0.8ex]
                                                                   2 \EPS     \\[0.8ex]
                                                                   (2\EPS \vec{u} + \vec{\phi}/\Pr_{\infty})^\T
                                                               \end{pmatrix},
\end{align}
respectively, where $p$ denotes the pressure, $H = E + p/\rho$ the specific enthalpy,
$\EPS$ the rate-of-strain tensor, $\vec{\phi} = \nabla \theta$ the temperature gradient, $\mu$ the dynamic viscosity,
$c_p$ the specific heat capacity, $\kappa$ the thermal conductivity and $L$ a characteristic length.

To close the system of equations \eqref{eq::Navier-Stokes} we assumed a Newtonian fluid, defined by a linear relationship between the
rate-of-strain tensor $\EPS$, given by
\begin{align}
    \label{eq::rate-of-strain}
    \EPS & := \frac{\grad \vec{u} + (\grad \vec{u})^\T}{2} - \frac{1}{3} \div(\vec{u}) \I,
\end{align}
and the deviatoric stress tensor $\TAU$, namely,
\begin{align}
    \TAU = \frac{2\EPS}{\Re_\infty}.
\end{align}
Additionally, we assumed that the heat flux, $\HEAT$, obeys Fourier's law of heat conduction, i.e.
\begin{align}
    \HEAT = - \frac{ \vec{\phi}}{\Re_\infty \Pr_\infty}.
\end{align}
Lastly, we assumed ideal gas conditions given by the equation of state
\begin{align}
    \gamma p & = \rho (\gamma - 1) \theta,
\end{align}
with heat capacity ratio $\gamma=1.4$ for air.
In this setting the isentropic speed of sound $c$
is defined as
\begin{align}
    c = \sqrt{\gamma \frac{p}{\rho}} = \sqrt{(\gamma - 1) \theta}.
\end{align}
The velocity $\vec{u}$, the pressure $p$ and the temperature $\theta$
can be expressed in terms of the conservative variables~$\vec{U}$ as follows:
\begin{align}
    \label{eq::primitive_variables_as_conservative}
    \vec{u}(\vec{U}) & := \frac{\rho \vec{u}}{\rho}, & p(\vec{U}) & := (\gamma - 1) \left( \rho E - \frac{|\rho \vec{u}|^2}{2\rho} \right), &  \theta(\vec{U}) & := \gamma \left( \frac{\rho E}{\rho} - \frac{|\rho \vec{u}|^2}{2\rho^2} \right).
\end{align}

Throughout this work we assume a two-dimensional setting\footnote{Every derivation can be extended to the three-dimensional case; however, we focus on the two-dimensional case for the sake of clarity and readability.} with a local orthogonal coordinate system $(\xi, \eta)$,
relatively oriented with respect to the global Cartesian coordinate system $(x, y)$ as shown in
Figure~\ref{fig::coordinate_system}. We introduce the orthonormal basis vectors $(\vec{n}, \vec{t})$ and define the corresponding
projection of the velocity $\vec{u}$ onto these basis vectors as follows
\begin{align}
    u_n  = \vec{u} \cdot \vec{n} = u_x n_x + u_y n_y, \qquad u_t  = \vec{u} \cdot \vec{t} = u_x t_x + u_y t_y.
\end{align}
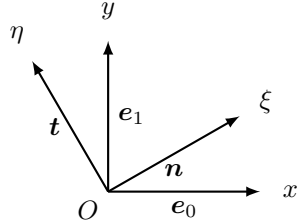
\begin{figure}[h]
    \centering
    \tikzsetnextfilename{coordinate}
    \begin{tikzpicture}[scale=2, >=latex]

        \draw[thick,->] (0,0) -- (1,0) node[midway, below] {$\vec{e}_0$};
        \draw[thick,->] (0,0) -- (0,1) node[midway, right] {$\vec{e}_1$};

        \node at (0,0) [below left] {$O$};

        \node at (1.2,0) {$x$};
        \node at (0,1.2) {$y$};

        \begin{scope}[rotate=30]
            \draw[thick,->] (0,0) -- (1,0) node[midway, below] {$\vec{n}$};
            \draw[thick,->] (0,0) -- (0,1) node[midway, left] {$\vec{t}$};

            \node at (1.2,0) {$\xi$};
            \node at (0,1.2) {$\eta$};
        \end{scope}

    \end{tikzpicture}
    \caption[Coordinate System]{Local orthogonal coordinate system $(\xi, \eta)$ with the orthonormal basis vectors $\vec{n}$ and $\vec{t}$.}
    \label{fig::coordinate_system}
\end{figure}

\subsection{Characteristic form of the governing equations}
The hyperbolic nature of the Navier-Stokes equations \eqref{eq::Navier-Stokes} with respect to time lies in the Euler
equations~\cite{hirschNumericalComputationInternal2002},
\begin{align}
    \label{eq::Euler-conservative}
    \pdt{\vec{U}} + \div(\vec{F}(\vec{U})) = \vec{0},
\end{align}
which are derived by neglecting the viscous contributions. From a characteristic point of view,
it is essential to express these equations in quasi-linear form
\begin{align}
    \label{eq::quasi-Euler-conservative}
    \pdt{\vec{U}} + \mat{A} \frac{\partial \vec{U}}{\partial \xi} + \mat{B} \frac{\partial \vec{U}}{\partial \eta} = \vec{0},
\end{align}
where the convective Jacobians given  are defined as \cite{rohdeEigenvaluesEigenvectorsEuler2001a}
\begin{subequations}
    \begin{align}
        \label{eq::convective-Jacobians_A}
        \mat{A} & := \frac{\partial (\vec{F} \vec{n})}{\partial \vec{U}} =
        \begin{pmatrix}
            0                                             & n_x                        & n_y                        & 0                \\[1.2ex]
            \frac{\gamma - 1}{2}| \vec{u} |^2 n_x - u_x u_n & u_n - (\gamma - 2) u_x n_x   & u_x n_y - (\gamma- 1) u_y n_x  & (\gamma - 1) n_x \\[1.2ex]
            \frac{\gamma - 1}{2}| \vec{u} |^2 n_y - u_y u_n & u_y n_x - (\gamma- 1) u_x n_y  & u_n - (\gamma - 2) u_y n_y   & (\gamma - 1) n_y \\[1.2ex]
            (\frac{\gamma - 1}{2}| \vec{u} |^2 - H) u_n   & H n_x - (\gamma - 1) u_x u_n & H n_y - (\gamma - 1) u_y u_n & \gamma u_n
        \end{pmatrix}, \\[1.2ex]
        \mat{B} & := \frac{\partial (\vec{F}\vec{t})}{\partial \vec{U}}  =
        \begin{pmatrix}
            0                                             & t_x                        & t_y                        & 0                \\[1.2ex]
            \frac{\gamma - 1}{2}| \vec{u} |^2 t_x - u_x u_t & u_t - (\gamma - 2) u_x t_x   & u_x t_y - (\gamma- 1) u_y t_x  & (\gamma - 1) t_x \\[1.2ex]
            \frac{\gamma - 1}{2}| \vec{u} |^2 t_y - u_y u_t & u_y t_x - (\gamma- 1) u_x t_y  & u_t - (\gamma - 2) u_y t_y   & (\gamma - 1) t_y \\[1.2ex]
            (\frac{\gamma - 1}{2}| \vec{u} |^2 - H) u_t   & H t_x - (\gamma - 1) u_x u_t & H t_y - (\gamma - 1) u_y u_t & \gamma u_t
        \end{pmatrix}.
    \end{align}
\end{subequations}
Unfortunately, the convective Jacobians $\mat{A}$ and $\mat{B}$ cannot be simultaneously diagonalized.
This arises from the fact that, in a multidimensional framework, waves can propagate in infinitely many directions,
and characteristic variables may not remain constant along characteristic lines \cite{laneyComputationalGasdynamics1998,liuNonreflectingBoundaryConditions2010}.

Nevertheless, a characteristic decomposition of $\mat{A}$ can be performed to obtain the eigenvalues $\EIG$ and
eigenvectors $\mat{P}$, defined as follows
\begin{subequations}
    \begin{align}
        \EIG         & := \diag{(\lambda_0, \lambda_1, \lambda_2, \lambda_3)} = \diag{(u_n - c, u_n, u_n, u_n + c)},                                                                                                                                                                                                                                                                              \\[1.2ex]
        \mat{P}      & :=
        \begin{pmatrix}
            \frac{1}{2 c^{2}}                                                         & \frac{1}{c^{2}}            & 0                           & \frac{1}{2 c^{2}}                                                         \\[1.2ex]
            - \frac{n_{x}}{2 c} + \frac{u_x}{2 c^{2}}                                   & \frac{u_x}{c^{2}}            & \rho n_{y}                  & \frac{n_{x}}{2 c} + \frac{u_x}{2 c^{2}}                                     \\[1.2ex]
            - \frac{n_{y}}{2 c} + \frac{u_y}{2 c^{2}}                                   & \frac{u_y}{c^{2}}            & - \rho n_{x}                & \frac{n_{y}}{2 c} + \frac{u_y}{2 c^{2}}                                     \\[1.2ex]
            \frac{1}{2(\gamma - 1)} - \frac{u_{n}}{2 c} + \frac{|\vec{u}|^2}{4 c^{2}} & \frac{|\vec{u}|^2}{2c^{2}} & \rho u_x n_{y} - \rho u_y n_{x} & \frac{1}{2(\gamma - 1)} + \frac{u_{n}}{2 c} + \frac{|\vec{u}|^2}{4 c^{2}}
        \end{pmatrix}, \\[1.2ex]
        \mat{P}^{-1} & :=
        \begin{pmatrix}
            c u_{n} + \frac{\gamma - 1}{2}|\vec{u}|^2     & - c n_{x} + u_x (1 - \gamma) & - c n_{y} + u_y (1 - \gamma) & \gamma - 1 \\[1.2ex]
            c^{2} - \frac{\gamma - 1}{2}|\vec{u}|^2       & - u_x (1 - \gamma)           & - u_y (1 - \gamma)           & 1 - \gamma \\[1.2ex]
            - \frac{u_x n_{y}}{\rho} + \frac{u_y n_{x}}{\rho} & \frac{n_{y}}{\rho}         & - \frac{n_{x}}{\rho}       & 0          \\[1.2ex]
            - c u_{n} + \frac{\gamma - 1}{2}|\vec{u}|^2   & c n_{x} + u_x (1 - \gamma)   & c n_{y} + u_y (1 - \gamma)   & \gamma - 1
        \end{pmatrix},
    \end{align}
\end{subequations}
By definition of eigendecomposition, it holds that $\mat{A} = \mat{P} \EIG \mat{P}^{-1}$, although the choice of eigenvectors $\mat{P}$ is not unique.
Furthermore, the eigenvalues $\EIG$, and consequently the convective Jacobian $\mat{A}$,
can be decomposed based on the sign of the eigenvalues
\begin{align}
    \label{eq::EIG_decomposition}
    \EIG  = \EIG^{-} + \EIG^{+}, \qquad \EIG^{\mp}  := \frac{\EIG \mp |\EIG|}{2}, \qquad \mat{A}^{\mp}  := \mat{P} \EIG^{\mp} \mat{P}^{-1}.
\end{align}
Additionally, we introduce the closely related identities and their projections onto the eigenvectors
\begin{align}
    \label{eq::Q_definition}
    \mat{I}  = \mat{I}^{-} + \mat{I}^{+}, \qquad \mat{I}^{\mp}  := \frac{\mat{I} \mp \mat{I}_n}{2}, \qquad \mat{Q}^{\mp}  := \mat{P} \mat{I}^{\mp} \mat{P}^{-1},
\end{align}
where $\mat{I}_n$ is defined as
\begin{alignat}{2}
    \label{eq::farfield_identity}
    \mat{I}_n & := \diag{(I_0, I_1, I_2, I_3)}, \qquad & I_i & = \begin{cases} 1 &  \lambda_i > 0, \\ 0 &  \lambda_i = 0, \\ -1 &  \lambda_i< 0.  \end{cases}
\end{alignat}

For the remainder of this work, the basis vector $\vec{n}$ is assumed to represent the outward unit normal vector at
any arbitrary boundary. Consequently, the basis vector $\vec{t}$ is the unit tangential vector, obtained by rotating $\vec{n}$ counterclockwise.
With these definitions, we can classify the eigenvalues~$\EIG^{\mp}$ as the speeds of incoming and outgoing waves, respectively. Additionally, the
Jacobians $\mat{A}^{-}$ and $\mat{A}^{+}$ are referred to as the incoming and outgoing convective Jacobians,
while the identities $\mat{I}^{-}$ and $\mat{I}^{+}$ are referred to as the incoming and outgoing
identities.

Notably, in the special case of $\mat{I}_n = \mat{0}$, the incoming and outgoing identities reduce to
$\mat{I}^{\mp} = \tfrac{1}{2} \mat{I}$, ensuring that both contributions are retained.

The characteristic form of the Euler equations is obtained by premultiplying the system of
equations~\eqref{eq::quasi-Euler-conservative} by the left eigenvectors $\mat{P}^{-1}$, while
exploiting the identity~$\mat{P} \mat{P}^{-1} = \I$. This procedure leads to
\begin{align*}
    \mat{P}^{-1}\pdt{\vec{U}} + \mat{P}^{-1}\mat{A} \frac{\partial \vec{U}}{\partial \xi} + \mat{P}^{-1}\mat{B} \frac{\partial \vec{U}}{\partial \eta} & =
    \mat{P}^{-1}\pdt{\vec{U}} + \mat{P}^{-1}\mat{A} \mat{P} \mat{P}^{-1}\frac{\partial \vec{U}}{\partial \xi} + \mat{P}^{-1}\mat{B}\mat{P} \mat{P}^{-1} \frac{\partial \vec{U}}{\partial \eta}                                                                                                                                                                             \\
                                                                                                                                                       & = \mat{P}^{-1}\pdt{\vec{U}} + \EIG \mat{P}^{-1}\frac{\partial \vec{U}}{\partial \xi} + \mat{\Theta} \mat{P}^{-1} \frac{\partial \vec{U}}{\partial \eta}                                                = \vec{0}, 
\end{align*}
with
\begin{align}
    \mat{\Theta} := \mat{P}^{-1} \mat{B} \mat{P} & = \begin{pmatrix}
                                                         u_t              & 0   & -\rho c^2 & 0                \\[0.8ex]
                                                         0                & u_t & 0         & 0                \\[0.8ex]
                                                         -\frac{1}{2\rho} & 0   & u_t       & -\frac{1}{2\rho} \\[0.8ex]
                                                         0                & 0   & -\rho c^2 & u_t
                                                     \end{pmatrix}.
\end{align}
We introduce the characteristic variables $\delta \vec{W}$, which correspond to the eigenvalues~$\EIG$ and are expressed in terms of variations as
\begin{align}
    \label{eq::variation-characteristic-variables}
    \delta \vec{W} & = \mat{P}^{-1} \delta \vec{U}  = \begin{pmatrix}
                                                          \delta p - \rho c \, \delta \vec{u} \cdot \vec{n} \\[0.8ex]
                                                          c^2 \delta\rho - \delta p                         \\[0.8ex]
                                                          -\delta \vec{u} \cdot \vec{t}                     \\[0.8ex]
                                                          \delta p + \rho c \,  \delta \vec{u} \cdot \vec{n}
                                                      \end{pmatrix}.
\end{align}
Here the first and last characteristic variables, $\delta \vec{W}_0$ and $\delta \vec{W}_{3}$, correspond to acoustic waves traveling upstream and downstream with
speeds $u_n - c$ and $u_n + c$, respectively. The remaining characteristic variables $\delta \vec{W}_1$ and  $\delta \vec{W}_2$ are
associated with entropy and vorticity waves, which are convected by the flow with velocity $u_n$.

With these definitions, the characteristic form of the Euler equations can be written in compact notation as
\begin{align}
    \label{eq::quasi-Euler-characteristic-2d}
    \pdt{\vec{W}} + \EIG\frac{\partial \vec{W}}{\partial \xi} + \mat{\Theta} \frac{\partial \vec{W}}{\partial \eta}=  \vec{0}. 
\end{align}
By neglecting the tangential contributions $\mat{\Theta} \partial \vec{W}/\partial \eta$, the
system decouples into $(d+2)$ independent equations,
\begin{align}
    \label{eq::quasi-Euler-characteristic-1d}
    \pdt{\vec{W}} + \EIG\frac{\partial \vec{W}}{\partial \xi} =  \vec{0}, 
\end{align}
which, rewritten in terms of the incoming and outgoing characteristic variables $\delta \vec{W}^{\mp} := \I^\mp \delta \vec{W}$, reads as
\begin{align}
    \label{eq::quasi-Euler-characteristic-1d_pm}
    \pdt{\vec{W}^\mp} + \EIG^\mp \frac{\partial \vec{W}^\mp}{\partial \xi} =  \vec{0}. 
\end{align}
\begin{remark}
    Premultiplying equation \eqref{eq::quasi-Euler-characteristic-1d} by $\mat{I}^{\mp}$ effectively acts as a filter, isolating the corresponding wave component.
\end{remark}
From the theory of hyperbolic systems, it is well-established that in the simplified one-dimensional case, the characteristic variables $\delta \vec{W}_i$ remain constant along the characteristic lines $d \xi/dt = \lambda_i$.
However, as noted earlier, in a multidimensional context, the eigenvectors $\mat{P}$ fail to diagonalize the tangential convective Jacobian $\mat{B}$,
resulting in a non-diagonal matrix $\mat{\Theta}$. The off-diagonal terms in $\mat{\Theta}$ act as source terms and introduce coupling between the
characteristic variables $\delta \vec{W}$. Consequently, in a multidimensional setting, the characteristic quantities do not remain constant along characteristic lines derived from the simplified system of equations \eqref{eq::quasi-Euler-characteristic-1d}.

By repeating the above procedure for the full set of equations \eqref{eq::Navier-Stokes}
we obtain the characteristic form of the Navier-Stokes equations
\begin{align}
    \label{eq::quasi-NS-characteristic-2d}
    \pdt{\vec{W}} + \EIG\frac{\partial \vec{W}}{\partial \xi} + \mat{\Theta} \frac{\partial \vec{W}}{\partial \eta}=  \mat{P}^{-1} \VIS,  
\end{align}
where the viscous contributions are defined as
\begin{align}
    \VIS := \div(\vec{G}(\vec{U}, \EPS, \vec{\phi})) = \frac{\partial (\vec{G} \vec{n})}{\partial \xi} + \frac{\partial (\vec{G} \vec{t})}{\partial \eta}.
\end{align}
Note that, rigorously speaking, the term \emph{characteristic} is not appropriate due to the elliptic nature of the viscous terms $\VIS$ \cite{poinsotBoundaryConditionsDirect}.

\subsection{Navier-Stokes characteristic boundary conditions}
\label{sec::NSCBC}
For the development of time-dependent boundary conditions, particularly non-reflective boundary conditions,
the characteristic Navier-Stokes system \eqref{eq::quasi-NS-characteristic-2d} has been widely used in literature \cite{poinsotBoundaryConditionsDirect, yooCharacteristicBoundaryConditions2005,thompsonTimeDependentBoundary1987,thompsonTimedependentBoundaryConditions1990,
    yooCharacteristicBoundaryConditions2007, lodatoThreedimensionalBoundaryConditions2008}. This approach is commonly referred to as characteristic Navier-Stokes boundary conditions (NSCBCs).
The system of equations~\eqref{eq::quasi-NS-characteristic-2d} is often reformulated in terms of characteristic amplitudes~$\AMP$ and transverse convective contributions~$\TRA$ as follows:
\begin{align}
    \pdt{\vec{W}} + \AMP + \TRA & =  \mat{P}^{-1} \VIS, 
\end{align}
where the terms $\AMP$ and $\TRA$ are defined as
\begin{subequations}
    \begin{align}
        \label{eq::characteristic_amplitudes_definition}
        \AMP & := \EIG \frac{\partial \vec{W}}{\partial \xi} =   \begin{pmatrix}
                                                                     (u_n - c) \left[ \frac{\partial p}{\partial \xi} - \rho c \frac{\partial \vec{u}}{\partial \xi} \cdot \vec{n} \right] \\[0.8ex]
                                                                     u_n  \left[ c^2 \frac{\partial \rho}{\partial \xi} - \frac{\partial p}{\partial \xi} \right]                          \\[0.8ex]
                                                                     u_n  \left[- \frac{\partial \vec{u}}{\partial \xi} \cdot \vec{t} \right]                                              \\[0.8ex]
                                                                     (u_n + c) \left[ \frac{\partial p}{\partial \xi} + \rho c \frac{\partial \vec{u}}{\partial \xi} \cdot \vec{n} \right]
                                                                 \end{pmatrix} = \begin{pmatrix}
                                                                                     \mathcal{L}_{acou} \\[0.8ex] \mathcal{L}_{ent} \\[0.8ex] \mathcal{L}_{vort} \\[0.8ex] \mathcal{L}_{acou}
                                                                                 \end{pmatrix}, \\[0.8ex]
        \TRA & := \mat{\Theta} \frac{\partial \vec{W}}{\partial \eta} =
        \begin{pmatrix}
            u_t  \left[ \frac{\partial p}{\partial \eta} - \rho c \frac{\partial \vec{u}}{\partial \eta} \cdot \vec{n} \right] + \rho c^2 \frac{\partial \vec{u}}{\partial \eta} \cdot \vec{t} \\[1.2ex]
            u_t \left[ c^2 \frac{\partial \rho}{\partial \eta} - \frac{\partial p}{\partial \eta} \right]                                                                                      \\[0.8ex]
            u_t  \left[- \frac{\partial \vec{u}}{\partial \eta} \cdot \vec{t} \right] - \frac{1}{\rho} \frac{\partial p}{\partial \eta}                                                        \\[0.8ex]
            u_t  \left[ \frac{\partial p}{\partial \eta} + \rho c \frac{\partial \vec{u}}{\partial \eta} \cdot \vec{n} \right] + \rho c^2 \frac{\partial \vec{u}}{\partial \eta} \cdot \vec{t}
        \end{pmatrix}.
    \end{align}
\end{subequations}

By further decoupling the characteristic amplitudes~$\AMP$ and transverse convective contributions~$\TRA$ into incoming and outgoing contributions, we have
\begin{subequations}
    \begin{alignat}{2}
        \AMP & =   \AMP^{-} +  \AMP^{+}, \qquad & \AMP^{\mp} & := \EIG^{\mp} \frac{\partial \vec{W}}{\partial \xi},     \label{eq::characteristic_amplitudes}                      \\
        \TRA & =  \TRA^{-} +  \TRA^{+},  \qquad & \TRA^\mp   & := \mat{I}^\mp \mat{\Theta} \frac{\partial \vec{W}}{\partial \eta}. \label{eq::transverse-convective-contributions}
    \end{alignat}
\end{subequations}
Thus, the system of equations can be expressed as
\begin{align}
    \label{eq::quasi-linear-incoming-outgoing}
    \pdt{\vec{W}} + \AMP^{-} +  \AMP^{+} + \TRA^{-} +  \TRA^{+} & =  \mat{P}^{-1} \VIS.
\end{align}

To ensure a well-posed problem, it is essential from a physical standpoint to specify the incoming wave information on artificial boundaries.
The NSCBC approach models the incoming amplitudes $\AMP^{-}$ and appropriately relaxes the
transverse terms $\TRA^-$ in equation \eqref{eq::quasi-linear-incoming-outgoing}, which contribute to the time evolution of the incoming characteristic wave
$\delta \vec{W}^-$.
Therefore, we rewrite the system of equations \eqref{eq::quasi-linear-incoming-outgoing} as
\begin{align}
    \label{eq::NSCBC-2d}
    \pdt{\vec{W}} + \tilde{\AMP}^{-} + \AMP^{+} + \beta \TRA^{-} + \TRA^{+} & =  \mat{P}^{-1} \VIS, 
\end{align}
where $\tilde{\AMP}^-$ denotes the modeled incoming characteristic amplitudes and $\beta > 0$ is a relaxation coefficient for the transverse contribution.
It is important to note that for outgoing waves $\delta \vec{W}^+$ no modeling is necessary,
as they are determined by the conditions prevailing within the interior.

The incoming characteristic amplitudes $\AMP^{-}$ are typically approximated by assuming a relaxation
toward a given target state $\vec{U}^{-}$. By applying a first-order, one-sided upwind finite-difference scheme along a characteristic length $L$, an approximation
to the derivative is obtained:
\begin{align}
    \label{eq::FD-approximation}
    \frac{\partial \vec{W}}{\partial \xi} \stackrel{\text{FD approx.}}{\approx} -\mat{P}^{-1} \frac{\vec{U} - \vec{U}^-}{L}.
\end{align}
Consequently,
the modeled incoming characteristic amplitudes $\widetilde{\AMP}^{-}$ are obtained by inserting the above approximation into equation \eqref{eq::characteristic_amplitudes}
\begin{align}
    \label{eq::modeled_amplitudes}
    \AMP^{-} \approx \widetilde{\AMP}^{-} & :=  \mat{D}^- \mat{P}^{-1} \left( \vec{U} - \vec{U}^- \right),
\end{align}
where a positive diagonal eigenvalue-like relaxation matrix is introduced:
\begin{align}
    \label{eq::relaxation_matrix}
    \mat{D}^- \sim -\frac{1}{L}\EIG^{-}.
\end{align}
The precise definition of the relaxation matrix depends on the specific type of boundary condition and will be addressed in the following sections.

Once the modeled characteristic amplitudes $\widetilde{\AMP}^{-}$ and the transverse relaxation coefficient $\beta$ are defined,
the time-dependent problem described by equation \eqref{eq::NSCBC-2d} is solved at the boundary.

In a viscous setting, the authors of \cite{poinsotBoundaryConditionsDirect} propose to modify the viscous contributions
$\VIS$ in the conservative equations \eqref{eq::NSCBC-2d} by appropriate viscous conditions. Usually, these
modifications include a zero-normal gradient condition for a stress component $\TAU_{ij}$ and/or a zero-normal gradient
condition for a heat flux component $\HEAT_i$ as these conditions gradually transition to the inviscid state when
viscosity and conductivity approach zero.

\begin{remark}
    To properly impose the modeled characteristic amplitudes $\widetilde{\AMP}^{-}$, it is necessary to specify both the target
    state~$\vec{U}^{-}$ and the relaxation matrix~$\mat{D}^-$, depending on the specific Mach regime and the type of boundary
    condition being applied. While selecting the target state is generally straightforward based on the boundary
    condition type, the relaxation matrix typically contains parameters that may require tuning to suit the specific circumstances.
    Additionally, the value of the tangential relaxation parameter $\beta \sim \mathcal{O}(1)$ must be adjusted accordingly.
\end{remark}

It is worth noticing that, by neglecting tangential and viscous contributions in \eqref{eq::NSCBC-2d}, we recover the locally one-dimensional
inviscid (LODI) approach, given by:
\begin{align}
    \label{eq::lodi}
    \pdt{\vec{W}} + \tilde{\AMP}^{-} + \AMP^{+} & = \vec{0}. 
\end{align}
While the LODI framework \eqref{eq::lodi} offers computationally efficiency due to the decoupling of the characteristic
variables, the one-dimensional inviscid assumption may be too restrictive for modeling multidimensional flows. This is because it
oversimplifies the complex interactions that occur in flows with variations across multiple dimensions.

In the following subsections, we present examples of boundary conditions for subsonic inflow and outflow conditions, which are
based on the NSCBC framework.

\subsubsection{Subsonic outflow}
\label{sec::subsonic_outflow}
\begin{figure}[h]
    \centering
    \tikzsetnextfilename{subsonic_outflow_amplitudes}
    \begin{tikzpicture}[>=latex]
        \fill[amber!30] (0,0) -- (3,0) to[out=40, in=-90] (4.6,2) to[out=90, in=-135] (4,4) -- (0, 4) -- cycle;
        \draw[thin, black, dashed] (4,4) -- (0,4) -- (0,0) --  (3,0);
        \draw[thick, black] (3,0) to[out=40, in=-90] (4.6,2) to[out=90, in=-135] (4,4);

        \draw[black, solid,line width=1pt, ->](4.6,2) -- (5.6,2);
        \node at (5.1, 2.2) {$\vec{n}$};
        \node at (0.5, 3.5) {$\Omega$};
        \node at (4.2, 3.5) {$\Gamma$};

        \draw[domain=2.5:3.5, smooth, variable=\x, matblue] plot (\x, {0.1*sin(250 * \x - 180) + 1});
        \draw[domain=2.5:3.5, smooth, variable=\x, matblue] plot (\x, {0.1*sin(500 * \x - 180) + 2});
        \draw[domain=2.5:3.5, smooth, variable=\x, matblue] plot (\x, {0.1*sin(500 * \x - 180) + 3});
        \draw[black, solid,line width=1pt, ->](2.5,1) -- (3.7,1) node [midway, above] {$u_n + c$};
        \draw[black, solid,line width=1pt, ->](2.5,2) -- (3.7,2) node [midway, above] {$u_n$};
        \draw[black, solid,line width=1pt, ->](2.5,3) -- (3.7,3) node [midway, above] {$u_n$};
        \node at (2, 2) {$\AMP^+$};

        \draw[domain=6.1:7.1, smooth, variable=\x, matred] plot (\x, {0.1*sin(750 * \x) + 2});
        \draw[black, solid,line width=1pt, ->](7.1,2) -- (5.9,2) node [midway, above] {$u_n - c$};

        \node at (8, 2) {$\AMP^{-}_{outflow}$};
    \end{tikzpicture}


    \caption{Schematic representation of the incoming, $\AMP^-_{outflow}$, and outgoing, $\AMP^+$, characteristic amplitudes at a subsonic outflow boundary $\Gamma$.}
    \label{fig::characteristic_waves_subsonic_outflow}

\end{figure}
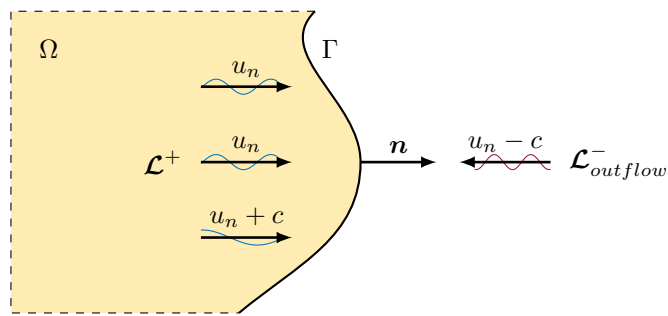

For a subsonic outflow, ensuring well-posedness requires specifying a
single physical boundary condition, typically the target pressure $p_\infty$. However, imposing this pressure strongly
can lead to unwanted reflections at the boundary. A more consistent approach from the physical standpoint
consists of modeling the single characteristic wave of purely acoustic nature entering the domain from downstream
\begin{align}
    \AMP^{-}_{outflow} & := \begin{pmatrix}
                                \mathcal{L}_{acou} \\[0.8ex] 0 \\[0.8ex] 0 \\[0.8ex] 0
                            \end{pmatrix} =
    \begin{pmatrix}
        (u_n - c) \left[ \frac{\partial p}{\partial \xi} - \rho c \frac{\partial \vec{u}}{\partial \xi} \cdot \vec{n} \right] \\[0.8ex] 0 \\[0.8ex] 0 \\[0.8ex] 0
    \end{pmatrix}.
\end{align}
A schematic representation of the characteristic amplitudes $\AMP$ for a domain $\Omega$ is shown in Figure~\ref{fig::characteristic_waves_subsonic_outflow}.

A widely used boundary condition~\cite{rudyNonreflectingOutflowBoundary1980,
    poinsotBoundaryConditionsDirect,yooCharacteristicBoundaryConditions2005} models the incoming acoustic
amplitude~$\AMP^{-}_{outflow}$~as:
\begin{align}
    \label{eq::L_acou_strikverda}
    \AMP^{-}_{outflow} \approx
    \widetilde{\AMP}^{-}_{outflow}  =
    \begin{pmatrix}
        D_p  (p - p_\infty) \\[0.8ex]
        0                   \\[0.8ex]
        0                   \\[0.8ex]
        0
    \end{pmatrix},
\end{align}
where the relaxation matrix $\mat{D}^-$ is given by:
\begin{align}
    \mat{D}^-                       = \begin{pmatrix}
                                          D_p & 0 & 0 & 0 \\
                                          0   & 0 & 0 & 0 \\
                                          0   & 0 & 0 & 0 \\
                                          0   & 0 & 0 & 0
                                      \end{pmatrix}.
\end{align}
The pressure difference $(p - p_\infty)$ ensures relaxation toward the target pressure $p_\infty$.
A specific form of the relaxation factor $D_p$ was introduced in~\cite{rudyNonreflectingOutflowBoundary1980}:
\begin{align}
    \label{eq::L_acou_strikverda_relaxation}
    D_p & := \frac{c \left(1 - \Ma_\infty^2\right)  \sigma}{L},
\end{align}
where $L$ is a characteristic length and $\sigma > 0$ is a non-negative constant. The optimal value of $\sigma$ was determined to be $\sigma_{opt} \approx 0.28$, \cite{rudyNonreflectingOutflowBoundary1980, shehadiNonReflectingBoundaryConditions2024},
though different values have been reported in the literature depending on the test case.
A vanishing constant, $\sigma = 0$, results in a perfectly non-reflecting boundary in a one-dimensional setting but entails a drift from the target pressure $p_\infty$,
as no information propagates upstream. Conversely, as the constant $\sigma$ increases significantly, $\sigma \gg 0$, the boundary becomes stiffer and more reflective.


If incoming tangential effects are non-negligible, such as in the wake of a bluff body, the inclusion of appropriately
damped tangential contributions can significantly improve the non-reflecting behavior of an outflow boundary. These contributions
can be modeled as follows:
\begin{align}
    \TRA^- & = \begin{pmatrix}
                   \mathcal{T}^-_{acou} \\[0.8ex] 0 \\[0.8ex] 0 \\[0.8ex] 0
               \end{pmatrix} =  \begin{pmatrix}
                                    u_t \left[ \frac{\partial p}{\partial \eta} - \rho c \frac{\partial \vec{u}}{\partial \eta} \cdot \vec{n}\right] + \rho c^2 \frac{\partial \vec{u}}{\partial \eta} \cdot \vec{t} \\[0.8ex] 0 \\[0.8ex] 0 \\[0.8ex] 0
                                \end{pmatrix},
\end{align}
where the incoming acoustic tangential component is given by:

\begin{align}
    \mathcal{T}^-_{acou} & \approx \tilde{\mathcal{T}}^-_{acou} := \beta_l u_t \left[ \frac{\partial p}{\partial \eta} - \rho c \frac{\partial \vec{u}}{\partial \eta} \cdot \vec{n} \right] + \beta_t \rho c^2 \frac{\partial \vec{u}}{\partial \eta} \cdot \vec{t}.
\end{align}
While the optimal choice of relaxation parameters $\beta_l, \beta_t$ remains an active area of research
\cite{liuNonreflectingBoundaryConditions2010,lodatoThreedimensionalBoundaryConditions2008, lodatoOptimalInclusionTransverse2012}, in this work, the approach proposed
in \cite{yooCharacteristicBoundaryConditions2007} for low Mach number flows is preferred. In that study, the authors perform a low-Mach asymptotic expansion of the incoming acoustic wave equation and
demonstrated that suitable transverse relaxation coefficients for the incoming acoustic wave are given by:
\begin{align}
    \beta = \beta_l & = \beta_t = \Ma_{\infty}.
\end{align}

In the viscous setting, reference \cite{poinsotBoundaryConditionsDirect} proposes to impose boundary conditions that enforce vanishing
normal derivatives of the normal-tangential component of the deviatoric stress tensor $\TAU$ and vanishing normal derivatives of the normal component of the
heat flux $\HEAT$, namely,
\begin{align}
    \frac{\partial(\vec{t}^\T \TAU \vec{n})}{\partial \xi}  = 0, \qquad \frac{\partial(\HEAT \cdot \vec{n})}{\partial \xi}  = 0.
\end{align}

\subsubsection{Subsonic inflow}
\label{sec::subsonic_inflow}
A subsonic inflow consists of $d+1$ incoming characteristics: acoustic, entropy and vorticity waves.
Following the rationale applied to outflow boundary conditions, the objective is to approximate
\begin{align}
    \label{eq::L_inflow_yoo}
    \vec{\mathcal{L}}^{-}_{\text{inflow}} := \begin{pmatrix}
                                                 \mathcal{L}^{-}_{acou}      \\[0.8ex]
                                                 \mathcal{L}^{-}_{entropy}   \\[0.8ex]
                                                 \mathcal{L}^{-}_{vorticity} \\[0.8ex]
                                                 0
                                             \end{pmatrix}
    = \begin{pmatrix}
          (u_n - c) \left[ \frac{\partial p}{\partial \xi} - \rho c \frac{\partial \vec{u}}{\partial \xi} \cdot \vec{n} \right] \\[0.8ex]
          u_n  \left[ c^2 \frac{\partial \rho}{\partial \xi} - \frac{\partial p}{\partial \xi} \right]                          \\[0.8ex]
          u_n  \left[- \frac{\partial \vec{u}}{\partial \xi} \cdot \vec{t} \right]                                              \\[0.8ex]
          0
      \end{pmatrix}.
\end{align}

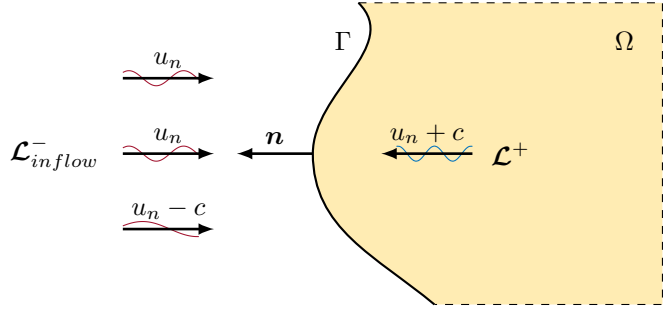
\begin{figure}[h]
    \centering
    \tikzsetnextfilename{subsonic_inflow_amplitudes}
    \begin{tikzpicture}[>=latex]
        \begin{scope}[yscale=1, xscale=-1]
            \fill[amber!30] (0,0) -- (3,0) to[out=40, in=-90] (4.6,2) to[out=90, in=-135] (4,4) -- (0, 4) -- cycle;
            \draw[thin, black, dashed] (4,4) -- (0,4) -- (0,0) --  (3,0);
            \draw[thick, black] (3,0) to[out=40, in=-90] (4.6,2) to[out=90, in=-135] (4,4);

            \draw[black, solid,line width=1pt, ->](4.6,2) -- (5.6,2);
            \node at (5.1, 2.2) {$\vec{n}$};
            \node at (0.5, 3.5) {$\Omega$};
            \node at (4.2, 3.5) {$\Gamma$};

            \draw[domain=2.5:3.5, smooth, variable=\x, matblue] plot (\x, {0.1*sin(750 * \x) + 2});
            \draw[black, solid,line width=1pt, ->](2.5,2) -- (3.7,2) node [midway, above] {$u_n + c$};

            \node at (2, 2) {$\AMP^+$};

            \draw[domain=6.1:7.1, smooth, variable=\x, matred] plot (\x, {0.1*sin(250 * \x + 180) + 1});
            \draw[domain=6.1:7.1, smooth, variable=\x, matred] plot (\x, {0.1*sin(500 * \x + 180) + 2});
            \draw[domain=6.1:7.1, smooth, variable=\x, matred] plot (\x, {0.1*sin(500 * \x + 180) + 3});
            \draw[black, solid,line width=1pt, ->](7.1,1) -- (5.9,1) node [midway, above] {$u_n - c$};
            \draw[black, solid,line width=1pt, ->](7.1,2) -- (5.9,2) node [midway, above] {$u_n$};
            \draw[black, solid,line width=1pt, ->](7.1,3) -- (5.9,3) node [midway, above] {$u_n$};

            \node at (8, 2) {$\AMP^{-}_{inflow}$};
        \end{scope}
    \end{tikzpicture}

    \caption{Schematic representation of the incoming, $\AMP^-_{inflow}$, and outgoing, $\AMP^+$, characteristic amplitudes at a subsonic inflow boundary $\Gamma$.}
    \label{fig::characteristic_waves_subsonic_inflow}
\end{figure}

A common approach in the finite difference framework is to model the incoming characteristic amplitudes
$\AMP^{-}_{inflow}$ by imposing the target velocity $\VEL_\infty$ and the static temperature $\theta_\infty$, leading to
\begin{align}
    \label{eq::L_inflow_u_T}
    \AMP^{-}_{\text{inflow}} \approx
    \tilde{\AMP}^{-}_{\text{inflow}}                   =
    \begin{pmatrix}
        -  D_{u_n} \rho c \left(\vec{u}- \VEL_{\infty}\right) \cdot \vec{n} \\[0.8ex]
        -  D_{\theta} \rho \left( \theta - \theta_\infty \right)            \\[0.8ex]
        -  D_{u_t} \left(\vec{u}  - \VEL_{\infty}\right) \cdot \vec{t}      \\[0.8ex]
        0
    \end{pmatrix}, \qquad \mat{D}^-  = \begin{pmatrix}
                                           D_{u_n} & 0          & 0       & 0 \\
                                           0       & D_{\theta} & 0       & 0 \\
                                           0       & 0          & D_{u_t} & 0 \\
                                           0       & 0          & 0       & 0
                                       \end{pmatrix},
\end{align}
where the precise definition of the relaxation factors $D_{u_n}, D_{\theta}, D_{u_t}$ can be found in \cite*[A.4]{yooCharacteristicBoundaryConditions2007}.
In this particular case, no modification of the viscous contributions $\VIS$ is necessary, since density $\rho$ and pressure $p$
are determined by the continuity equation, which excludes any viscous terms \cite{poinsotBoundaryConditionsDirect}.

A different approach, discussed in \cite{shehadiPolynomialcorrectionNavierStokesCharacteristic2024a,
    shehadiNonReflectingBoundaryConditions2024}, imposes the mass flow rate by specifying the target
density~$\rho_\infty$ and the target velocity $\VEL_\infty$, leading to
\begin{align}
    \label{eq::L_inflow_u_rho}
    \AMP^{-}_{\text{inflow}} \approx
    \tilde{\AMP}^{-}_{\text{inflow}}                     =
    \begin{pmatrix}
        -  D_{u_n} \rho c \left( \vec{u} - \VEL_{\infty}\right) \cdot \vec{n} \\[0.8ex]
        D_{\rho} c^2 \left( \rho - \rho_\infty \right)                        \\[0.8ex]
        -  D_{u_t} \left( \vec{u} - \VEL_{\infty} \right) \cdot \vec{t}       \\[0.8ex]
        0
    \end{pmatrix}, \qquad \mat{D}^-  = \begin{pmatrix}
                                           D_{u_n} & 0        & 0       & 0 \\
                                           0       & D_{\rho} & 0       & 0 \\
                                           0       & 0        & D_{u_t} & 0 \\
                                           0       & 0        & 0       & 0
                                       \end{pmatrix},
\end{align}
where the relaxation factors $D_{u_n}, D_{\rho}, D_{u_t}$ are defined as
\begin{align}
    \label{eq::L_inflow_u_rho_relaxation}
    D_{u_n}  =  \sigma_u \frac{c \left(1 - \Ma_\infty^2\right)}{L}, \qquad D_{\rho} = \sigma_\rho \frac{c}{L}, \qquad D_{u_t}  = \sigma_t \frac{c}{L},
\end{align}
with $\sigma_u$, $\sigma_p$ and $\sigma_t$ being tunable parameters that influence the rate of relaxation towards the target values.
In contrast to the previous ansatz, reference \cite{poinsotBoundaryConditionsDirect} proposes an alternative strategy that enforces
a vanishing normal derivative of the normal-normal component of the deviatoric stress tensor $\TAU$ at the boundary, namely,
\begin{align}
    \frac{\partial(\vec{n}^\T \TAU \vec{n})}{\partial \xi} & = 0.
\end{align}
This condition ensures that normal stress components do not introduce spurious numerical artifacts at the inflow boundary,
particularly in cases where viscosity plays a significant role.

For imposition of conditions on an inflow, including total temperature $\theta_t$, total pressure $p_t$ and
flow direction $\vec{u}/\|\vec{u}\|$ refer to
\cite{shehadiNonReflectingBoundaryConditions2024,odierCharacteristicInletBoundary2019}.

The treatment of tangential contributions $\TRA^-$ for incoming waves in subsonic inflows remains barely discussed in literature. Typically, these contributions are
neglected by setting the tangential relaxation coefficient $\beta = 0$.

\subsection{Generalized characteristic relaxation boundary conditions}
\label{sec::GRCBC}
In their seminal work \cite{pirozzoliGeneralizedCharacteristicRelaxation2013} from 2013, the authors introduced the concept of Generalized
Characteristic Relaxation Boundary Conditions (GRCBCs) within a finite difference (FD) framework.
The primary distinction between GRCBCs and the NSCBCs lies in the definition of the relaxation matrix $\mat{D}^-$ in equation~\eqref{eq::relaxation_matrix}.

GRCBCs originate from a first-order, one-sided upwind finite difference scheme, where the characteristic length $L = |\vec{x}_b - \vec{x}_g|$
represents the distance between a boundary point~$\vec{x}_b$ and a ghost point~$\vec{x}_g$.
Inside the domain, the CFL condition dictates that, for the maximum propagation speed $\lambda_{max} = \max{\lbrace |u_n + c|, |u_n - c|\rbrace}$, the time step
$\Delta t$ must adhere to the stability constraint
\begin{align*}
    \Delta t  \lesssim \frac{L_i}{\lambda_{max}}, \qquad L_i := |\vec{x}_b - \vec{x}_i|,
\end{align*}
where $\vec{x}_i$ is the first interior mesh point.
From this consideration, it follows that the length $L$ must satisfy $L \geq L_i$ to comply with the CFL condition.
The natural choice, therefore, is to set $L = L_i$, leading to a relaxation matrix that depends solely on the mesh size, for instance,
\begin{align}
    \label{eq::relaxation_matrix_pirozzoli}
    \mat{D}^- & := -\frac{1}{L_i} \EIG^-.
\end{align}

Since the only difference between GRCBCs and NSCBCs lies in the definition of the relaxation matrix~$\mat{D}^-$, other aspects (including the imposition of the
transverse relaxation parameter $\beta$, the target state $\vec{U}^-$, and the viscous contributions $\VIS$) remain unchanged
from the NSCBC framework.

It is important to note that, in the finite element (FE) context, the concept of a first interior mesh point $\vec{x}_i$ is inherently ambiguous due to
differences in the underlying spatial discretization.  Addressing this issue within the FE framework
constitutes a key novelty of this work and will be extensively discussed in Section~\ref{sec::grcbc_hdg}.

\section{Hybridizable Discontinuous Galerkin method for compressible flows}
\label{sec::hdg_for_compressible}

\subsection{Mesh partition and functional spaces}
Consider a partition $\mesh$ of the domain $\Omega$ into a shape-regular, uniform and non-overlapping $d$-simplices with
\begin{align*}
    \max_{T \in \mesh} \frac{h_T^d}{|T|} & \leq C    \qquad C > 0 \quad \forall T \in \mesh, & \qquad
    h_T                                  & \leq C h  \qquad C > 0 \quad \forall T \in \mesh,
\end{align*}
where $h_T$ denotes the element diameter, and the global mesh size $h$ is defined as
\begin{align*}
    h_T  := \text{diam}(T), \qquad h  := \max\limits_{T \in \mesh}{h_T}.
\end{align*}
Further, let $\facets$ denote the set containing the mesh skeleton,
which consists of the sets of interior facets $\facets^{\text{int}}$ and exterior facets $\facets^{\text{ext}}$, given by
\begin{align*}
    \facets              := \bigcup_{T \in \mesh} \partial T,   \qquad
    \facets^{\text{int}} := \facets \backslash \partial \Omega, \qquad
    \facets^{\text{ext}} := \facets \cap \partial \Omega.
\end{align*}
On a shared interface $F = \partial T^+ \cap \partial T^- \in \facets^{\text{int}}$ between two adjacent elements $T^+$ and $T^-$ with outward
unit normal vectors $\vec{n}^+ = -\vec{n}^-$, let $\vec{V}^\pm$ denote the traces of the restriction of a matrix-valued function $\vec{v}$ onto the respective elements.
The jump operator for the traces, see \cite{dipietroMathematicalAspectsDiscontinuous2012,cockburnStaticCondensationHybridization2016}, is then defined as
\begin{align*}
    \jump{\vec{v} \vec{n}} := \vec{v}^+\vec{n}^+ + \vec{v}^-\vec{n}^-.
\end{align*}

Let $L^2(\Omega)$ denote the space of real-valued square-integrable functions over $\Omega$, equipped with the inner product and norm
\begin{align*}
    \left( u, v \right) := \sum_{T \in \mesh} \int_{T} u v \, d\bm{x}, \qquad
    \|u\|               := \sqrt{( u, u )},
\end{align*}
and define $L^2(\facets)$ as the space of real-valued square-integrable functions on $\facets$, with inner product and norm
\begin{align*}
    {\left\langle u, v \right\rangle}_{\facets} := \sum_{F \in \facets} \int_{F} u v \, d\bm{s}, \qquad
    \|u\|_{\facets}                             := \sqrt{\left\langle u, v \right\rangle}_{\facets}.
\end{align*}
For an arbitrary space $V$, define $V(\Omega, \mathbb{R}^n)$ as the vector-valued space, whose $n$-components are in~$V$, and $V(\Omega, \mathbb{R}^{n \times n})$
as the matrix-valued space, whose $n^2$-components are in $V$. Then $L^2(\Omega, \mathbb{R}^{n \times n}_{\mathrm{sym}})$ denotes the space of 
real-valued square-integrable symmetric tensors over $\Omega$, defined as
\begin{align*}
    L^2(\Omega, \mathbb{R}^{n \times n}_{\mathrm{sym}}) := \left\{ \mat{v} \in L^2(\Omega, \mathbb{R}^{n \times n}) : \mat{v} = \mat{v}^\T \right\}.
\end{align*}

The spaces of polynomials of degree at most $k$, restricted to an element $T$ or a facet $F$, are denoted by $\mathbb{P}^k(T)$ and $\mathbb{P}^k(F)$, respectively.
The broken polynomial spaces on $\mesh$ and $\facets$ are then given by:
\begin{align*}
    \mathbb{P}^k(\mesh)                 & := \lbrace v \in L^2(\Omega) : v|_T \in \mathbb{P}^k(T) \, \forall T \in \mesh \rbrace,                                            \\
    \mathbb{P}^k(\facets)               & := \lbrace v \in L^2(\facets) : v|_F \in \mathbb{P}^k(F) \, \forall F \in \facets \rbrace,                                         \\
    \mathbb{P}^k(\mesh, \mathbb{R}^n)   & := \lbrace \vec{v} \in L^2(\Omega, \mathbb{R}^n) : \vec{v}|_T \in \mathbb{P}^k(T, \mathbb{R}^n) \, \forall T \in \mesh \rbrace,    \\
    \mathbb{P}^k(\facets, \mathbb{R}^n) & := \lbrace \vec{v} \in L^2(\facets, \mathbb{R}^n) : \vec{v}|_F \in \mathbb{P}^k(F, \mathbb{R}^n) \, \forall F \in \facets \rbrace, \\
    \mathbb{P}^k(\mesh, \mathbb{R}^{n \times n}_{\mathrm{sym}})   & := \lbrace \vec{v} \in L^2(\Omega, \mathbb{R}^{n \times n}_{\mathrm{sym}}) : \vec{v}|_T \in \mathbb{P}^k(T, \mathbb{R}^{n \times n}_{\mathrm{sym}}) \, \forall T \in \mesh \rbrace.
\end{align*}

\subsection{Hybridizable Discontinuous Galerkin method}
We consider a recently developed mixed Hybridizable Discontinuous Galerkin method presented in
\cite{vila-perezHybridisableDiscontinuousGalerkin2021}, which, in addition to the conservative variables~$\vec{U}$,
introduces as mixed variables the rate-of-strain~tensor~$\EPS$, see \eqref{eq::rate-of-strain}, and the temperature gradient $\vec{\phi} = \grad{\theta}$.
The discrete variational formulation of the compressible Navier-Stokes equations
\begin{align}
    \label{eq::Navier-Stokes-hdg}
    \pdt{\vec{U}} + \div(\vec{F}(\vec{U}) - \vec{G}(\vec{U}, \EPS, \vec{\phi})) = \vec{0}   \quad \text{in} \quad \Omega \times (0, t_{end}] ,
\end{align}
is derived by the standard procedure of multiplying by a suitable test function, integrating by parts, and enforcing numerical flux continuity. The resulting formulation reads:

Find $\left(\vec{U}_h,\hat{\vec{U}}_h, \EPS_h, \vec{\phi}_h \right) \in U_h \times \hat{U}_h \times \Xi_h \times \Theta_h$ such
that
\begin{subequations}
    \label{eq::discrete_variational_formulation}
    \begin{align}
        \label{eq::discrete_variational_formulation_interior}
        \sum_{T \in \mesh} \int_{T} \pdt{\vec{U}_h} \cdot \vec{V}_h \, d\bm{x} -
        \int_{T} \left(\vec{F}(\vec{U}_h) - \vec{G}(\vec{U}_h, \EPS_h, \vec{\phi}_h)\right)  : \grad{\vec{V}_h} \, d\bm{x}+
        \int_{\partial T} (\hat{\vec{F}}_h - \hat{\vec{G}}_h) \vec{n} \cdot \vec{V}_h   \, d\bm{s}                                                                                                                                             & = 0, \\
        \label{eq::discrete_variational_formulation_gamma}
        - \sum_{F \in \facets^{\text{int}}} \int_{F} \jump{(\hat{\vec{F}}_h - \hat{\vec{G}}_h) \vec{n}} \cdot \hat{\vec{V}}_h \, d\bm{s}
        + \sum_{F \in \facets^{\text{ext}}} \int_{F} \hat{\vec{\Gamma}}_h \cdot \hat{\vec{V}}_h  \, d\bm{s}                                                                                                                     & = 0, \\
        \sum_{T  \in \mesh} \int_{T} \EPS_h : \mat{\zeta}_h \, d\bm{x} + \int_{T} \VEL_h \cdot \div(\mat{\zeta}_h - \frac{1}{3}\tr(\mat{\zeta}_h)\I) \, d\bm{x}  - \int_{\partial T} \hat{\VEL}_h \cdot \left[\mat{\zeta}_h - \frac{1}{3}\tr(\mat{\zeta}_h)\I \right] \vec{n} \, d\bm{s} & = 0, \\
        \sum_{T  \in \mesh} \int_{T} \vec{\phi}_h \cdot \vec{\varphi}_h \, d\bm{x} + \int_{T} \theta_h \div(\vec{\varphi}_h) \, d\bm{x} - \int_{\partial T} \hat{\theta}_h \vec{\varphi}_h \cdot \vec{n}    \, d\bm{s}                                    & = 0,
    \end{align}
\end{subequations}
for all $\left(\vec{V}_h,\hat{\vec{V}}_h, \mat{\zeta}_h, \vec{\varphi}_h \right) \in U_h \times \hat{U}_h \times \Xi_h \times \Theta_h$, with the discrete spaces chosen as
\begin{alignat*}{2}
    U_h       & := L^2\left( (0, t_{end}] ; \mathbb{P}^k(\mesh, \mathbb{R}^{4})          \right),      \qquad &
    \hat{U}_h & := L^2\left( (0, t_{end}] ; \mathbb{P}^k(\facets, \mathbb{R}^{4})         \right),              \\
    \Xi_h     & := L^2\left( (0, t_{end}] ; \mathbb{P}^k(\mesh, \mathbb{R}^{2 \times 2}_{\mathrm{sym}}) \right),  \qquad &
    \Theta_h  & := L^2\left( (0, t_{end}] ; \mathbb{P}^k(\mesh, \mathbb{R}^{2})                 \right).
\end{alignat*}
To define the discrete velocities $\vec{u}_h := \vec{u}(\vec{U}_h)$, $\hat{\vec{u}}_h := \vec{u}(\hat{\vec{U}}_h)$, and the discrete temperatures $\theta_h := \theta(\vec{U}_h)$, $\hat{\theta}_h := \theta(\hat{\vec{U}}_h)$, we used equation \eqref{eq::primitive_variables_as_conservative}.
In the formulation, $\hat{\vec{\Gamma}}_h$ represents the boundary operator, detailed below, imposing the boundary conditions, whereas
the numerical fluxes are defined as
\begin{subequations}
    \begin{align}
        \hat{\vec{F}}_h \vec{n} & := \vec{F}(\hat{\vec{U}}_h) \vec{n} + \mat{S}_c (\vec{U}_h - \hat{\vec{U}}_h),                       \\
        \hat{\vec{G}}_h \vec{n} & := \vec{G}(\hat{\vec{U}}_h, \EPS_h, \vec{\phi}_h) \vec{n} - \mat{S}_d (\vec{U}_h - \hat{\vec{U}}_h),
    \end{align}
\end{subequations}
where $\mat{S}_c, \mat{S}_d$ are convective and viscous stabilization matrices, respectively, chosen to ensure
the correct physical behavior of the flow. A detailed discussion on suitable stabilization
matrices and their connection to Riemann solvers can be found in \cite{vila-perezHybridisableDiscontinuousGalerkin2021,
    peraireEmbeddedDiscontinuousGalerkin2011, peraireHybridizableDiscontinuousGalerkin2010}, whereas the specific expressions employed 
in this study are detailed in Section~\ref{sec::experiments}.

\begin{remark}
    Following the discussion in \cite{vila-perezHybridisableDiscontinuousGalerkin2021, giacominiTutorialHybridizableDiscontinuous2020},
    the symmetry of the discrete rate-of-strain tensor~$\EPS_h = \EPS_h^T$ is enforced strongly by properly rearranging the unique off-diagonal components using Voigt notation.
    For instance in two dimensions, this approach decreases the number of components in the matrix-valued space $\Xi_h$ from four to three, thereby significantly reducing the computational cost.
\end{remark}

\subsection{Boundary conditions and stability considerations}
In HDG methods boundary conditions are imposed weakly via the
boundary operator~$\hat{\vec{\Gamma}}_h$ in \eqref{eq::discrete_variational_formulation_gamma}. Common boundary
conditions include: 
\begin{subequations}
    \begin{alignat}{4} \label{eq::adiabatic_wall}
        \text{(Adiabatic wall)} & \qquad &  \hat{\vec{\Gamma}}_{h} = \hat{\vec{\Gamma}}_{\textit{ad}} & := \begin{pmatrix}
                                                                           \hat{\rho}_h - \rho_h \\[0.8ex] \hat{\rho \VEL}_h \\[0.8ex] 1/(\Re_\infty \Pr_\infty)\vec{\phi}_h \cdot \vec{n}
                                                                       \end{pmatrix} - \mat{S}_d \begin{pmatrix}
                                                                        0 \\ \bm{0} \\ \rho E_h - \hat{\rho E}_h
                                                                       \end{pmatrix},  \\[0.8ex]
        \text{(Isothermal wall)}& \qquad &  \hat{\vec{\Gamma}}_{h} =\hat{\vec{\Gamma}}_{\textit{iso}} & := \hat{\vec{U}}_h -  \begin{pmatrix}
                                                                                            \rho_h \\[0.8ex] \vec{0} \\[0.8ex] \rho_h \theta_{\text{wall}}/\gamma
                                                                                        \end{pmatrix}, \\[0.8ex]
        \text{(Inviscid wall)} & \qquad & \hat{\vec{\Gamma}}_{h} =\hat{\vec{\Gamma}}_{\textit{inv}} & := \hat{\vec{U}}_h -\begin{pmatrix}
                                                                                          \rho_h \\[0.8ex] (\I - \vec{n} \otimes \vec{n}) \rho \VEL_h \\[0.8ex] \rho E_h
                                                                                      \end{pmatrix}.
    \end{alignat}
\end{subequations}
In addition, in aerodynamics, it is common to impose the static pressure on a subsonic outflow boundary by setting
\begin{align}
    \label{eq::subsonic_outflow}
    \hat{\vec{\Gamma}}_{h} = \hat{\vec{\Gamma}}_{p_\infty} := \hat{\vec{U}}_h - \vec{U}_{p_\infty},
\end{align}
where the pressure target state $\vec{U}_{p_\infty}$ is defined as
\begin{align}
    \vec{U}_{p_\infty}                 := \begin{pmatrix}
                                         \rho_h \\[0.8ex] \rho \VEL_h \\[0.8ex] p_\infty/(\gamma - 1) + \rho_h |\VEL_h|^2/2
                                     \end{pmatrix}.
\end{align}
This boundary condition is particularly effective in the wake of a body, as it allows entropy waves and vorticity waves to leave the domain. However, it has a significant drawback:
it causes the full outgoing acoustic wave to be reflected back into the domain, as the imposed target pressure $p_\infty$ acts as a hard boundary
for pressure waves.

A more suitable boundary condition for compressible flow, commonly referred to as the far-field boundary condition \cite{vila-perezHybridisableDiscontinuousGalerkin2021, nguyenHybridizableDiscontinuousGalerkin2012,
    peraireHybridizableDiscontinuousGalerkin2010}, uses a characteristic
approach 
\begin{align}
    \label{eq::farfield}
    \hat{\vec{\Gamma}}_{h} = \hat{\vec{\Gamma}}_{\infty} := \hat{\mat{A}}^+(\hat{\vec{U}}_h - \vec{U}_h) -  \hat{\mat{A}}^-(\hat{\vec{U}}_h - \vec{U}_\infty),
\end{align}
where the far-field target state $\vec{U}_\infty$ is given by
\begin{align}
    \vec{U}_{\infty}                 := \begin{pmatrix}
                                       \rho_\infty \\[0.8ex] \rho_\infty \VEL_\infty \\[0.8ex] \rho_\infty E_\infty
                                   \end{pmatrix}.
\end{align}
Here, $\hat{\mat{A}}$ is the convective Jacobian, as defined in \eqref{eq::convective-Jacobians_A}, evaluated at the intermediate state $\hat{\vec{U}}_h$. The incoming and outgoing convective Jacobians, denoted by $\hat{\mat{A}}^{\mp}$, 
are obtained via the decomposition of $\hat{\mat{A}}$ evaluated at the intermediate state $\hat{\vec{U}}_h$ into its characteristic wave contributions, as given in \eqref{eq::EIG_decomposition}.
This decomposition into incoming and outgoing contributions, based on the Jacobians $\hat{\mat{A}}^{\mp}$, ensures an upwind discretization, as it correctly determines the wave propagation 
direction that is used to determine the intermediate state $\hat{\vec{U}}_h$. Specifically, outgoing waves, $\hat{\mat{A}}^+(\hat{\vec{U}}_h - \vec{U}_h)$, 
are determined based on the interior state $\vec{U}_h$, while incoming waves, $\hat{\mat{A}}^-(\hat{\vec{U}}_h - \vec{U}_\infty)$, are determined based on the exterior state $\vec{U}_\infty$.

In supersonic flows, where information cannot propagate upstream, the specified boundary condition ensures the correct physical behavior.
Specifically, in an outflow scenario, $\hat{\mat{A}}^{-} = \mat{0}$, meaning only outgoing characteristics are considered. In an inflow scenario, $\hat{\mat{A}}^{+} = \mat{0}$, ensuring
that only the upstream state influences the intermediate state $\hat{\vec{U}}_h$.

For subsonic flows, the far-field boundary operator \eqref{eq::farfield} effectively absorbs planar acoustic waves,
but significant reflection may occur when oblique waves and vortices cross the boundary, as observed in the numerical experiments in Section~\ref{sec::experiments}.

A limitation of the formulation in \eqref{eq::farfield} arises when the boundary is aligned in flow direction, i.e. $u_n = 0$, as the
convective Jacobian $\hat{\mat{A}}$ then has $d$ zero eigenvalues. This might lead to difficulties when solving the
system, e.g. when using an iterative solver. A practical workaround has been described in \cite{vila-perezHybridisableDiscontinuousGalerkin2021} by introducing a
threshold $\varepsilon \ll 1$, such that
\begin{align*}
    \tilde{\lambda_{i}}^+ = \max(\lambda_{i}^+, \varepsilon), \qquad \tilde{\lambda_{i}}^- = \min(\lambda_{i}^-, -\varepsilon).
\end{align*}

Ensuring that the convective Jacobian $\hat{\mat{A}}$ has non-zero eigenvalues, we propose a numerically more stable discretization by premultiplying \eqref{eq::farfield} by the inverse of
the absolute convective Jacobian $|\hat{\mat{A}}|^{-1} = \hat{\bm{P}} |\hat{\EIG}|^{-1} \hat{\bm{P}}^{-1}$, where the 
hat on all matrices is to recall that they are computed at the intermediate state $\hat{\vec{U}}_h$. This results in the modified boundary condition:
\begin{align*}
    |\hat{\mat{A}}|^{-1} \left[ \hat{\mat{A}}^+(\hat{\vec{U}}_h - \vec{U}_h) -  \hat{\mat{A}}^-(\hat{\vec{U}}_h - \vec{U}_\infty) \right]                                                   &     
    = |\hat{\mat{A}}|^{-1} (\hat{\mat{A}}^+ -  \hat{\mat{A}}^-) \hat{\vec{U}}_h - |\hat{\mat{A}}|^{-1} \hat{\mat{A}}^+ \vec{U}_h + |\hat{\mat{A}}|^{-1} \hat{\mat{A}}^- \vec{U}_\infty.
\end{align*}
Since $\hat{\mat{A}}^+ -  \hat{\mat{A}}^- = |\hat{\mat{A}}|$, this simplifies to 
\begin{align*}
    \hat{\vec{U}}_h - |\hat{\mat{A}}|^{-1} \hat{\mat{A}}^+ \vec{U}_h + |\hat{\mat{A}}|^{-1} \hat{\mat{A}}^- \vec{U}_\infty                                                           = \vec{0}.
\end{align*}
From the definitions in \eqref{eq::Q_definition}, we identify $\hat{\mat{Q}}^+ = |\hat{\mat{A}}|^{-1} \hat{\mat{A}}^+$ and $\hat{\mat{Q}}^- = -|\hat{\mat{A}}|^{-1} \hat{\mat{A}}^-$. Then, the final 
stabilized boundary operator is expressed as
\begin{align}
    \label{eq::farfield_mod}
    \hat{\vec{\Gamma}}_{h} = \hat{\vec{\Gamma}}_{\infty}^{\mat{Q}} & := \hat{\vec{U}}_h - \hat{\mat{Q}}^+\vec{U}_h -  \hat{\mat{Q}}^- \vec{U}_\infty
    = \hat{\vec{U}}_h - \frac{\vec{U}_h + \vec{U}_\infty}{2} - \hat{\mat{P}}\hat{\mat{I}}_n \hat{\mat{P}}^{-1} \frac{\vec{U}_h - \vec{U}_\infty}{2}.
\end{align}

\begin{remark}
    Assuming $\hat{\mat{Q}}^+$ and $\hat{\mat{Q}}^-$ are constant, then \eqref{eq::farfield_mod} resembles the solution of
    a linear hyperbolic problem with constant coefficients.
\end{remark}

\begin{remark}
    In the case of a singular vanishing eigenvalue $\lambda_i = 0$, the matrix $\hat{\mat{I}}_n$ has a zero diagonal entry, and
    \eqref{eq::farfield_mod} simplifies to averaging between the two states $\vec{U}_h$ and $\vec{U}_\infty$
    for the corresponding wave, analogously to the Lax-Friedrichs flux \cite{vila-perezHybridisableDiscontinuousGalerkin2021}. The same 
    holds true for multiple vanishing eigenvalues.
\end{remark}
\section{Characteristic boundary conditions for Hybridizable Discontinuous Galerkin method}
\label{sec::cbc_for_hdg}
In this work we devise for the first time the Navier-Stokes characteristic boundary conditions (NSCBCs) of Subsection~\ref{sec::NSCBC}
and the generalized characteristic relaxation boundary conditions~(GRCBCs) of Subsection~\ref{sec::GRCBC} in a Hybridizable Discontinuous Galerkin (HDG) setting.
To simplify the explanation of the implementation, we initially consider the one-dimensional, inviscid case \eqref{eq::lodi} in conservative form, thus
\begin{align}
    \label{eq::lodi_conservative}
    \mat{P} \left(\pdt{\vec{W}} + \AMP^- + \AMP^+ \right) = \pdt{\vec{U}} + \mat{P}\AMP^- + \mat{P}\AMP^+ & = \vec{0} \quad \text{on} \quad \partial \Omega \times (0, t_{end}],
\end{align}
while the extension to the multidimensional viscous setting is addressed subsequently.

For ease of notation, every scalar, vector or matrix evaluated at the intermediate state $\hat{\vec{U}}_h$ is denoted by
a hat symbol, e.g. the eigenvectors $\hat{\mat{P}}$.

We begin by discretizing equation \eqref{eq::lodi_conservative} at the boundary $\partial \Omega$ using an implicit Euler method
considering a time step $\Delta t$ on a time interval $0 = t_0 < t_n < \ldots \leq t_{end}$ with $t_n := n \Delta t$. Then the discrete system of equations reads as
\begin{align}
    \frac{\hat{\vec{U}}^{n+1}_h -  \hat{\vec{U}}^{n}_h}{\Delta t} + \hat{\mat{P}}^{n+1} \hat{\AMP}^{n+1, -} + \hat{\mat{P}}^{n+1} \hat{\AMP}^{n+1, +}   = \vec{0},
\end{align}
where the superscript $n \in \lbrace 0, 1, \ldots \rbrace$ refers to an evaluation at time $t_{n}$.
Unless otherwise stated, we omit the superscript $n+1$, referring to the time evaluation $t_{n+1}$, for the remainder of this section.

\begin{remark}
    \label{rem::cbc_time_discretization}
    The choice of time marching scheme is not restricted to an implicit Euler, but has been chosen for sake of simplicity in the
    derivation. We recommend using the same time marching scheme as for the elementwise conservative variable $\vec{U}_h$ in
    \eqref{eq::discrete_variational_formulation_interior}.
\end{remark}


We rearrange terms by multiplying with the denominator of the time marching scheme i.e. $\Delta t$, leading to
\begin{align}
    \label{eq::discrete_lodi}
    \hat{\vec{U}}_h -  \hat{\vec{U}}^{n}_h + \Delta t \hat{\mat{P}}\hat{\AMP}^- +  \Delta t \hat{\mat{P}}\hat{\AMP}^+ = \vec{0}.
\end{align}

By mimicking the filter from~\eqref{eq::farfield}, outgoing waves can be simply determined by $\hat{\mat{A}}^+(\hat{\vec{U}}_h - \vec{U}_h)$, whereas
we aim to model incoming waves to achieve better accuracy. Therefore, we premultiply~\eqref{eq::discrete_lodi} by the incoming Jacobian $\hat{\mat{A}}^-$, resulting in
\begin{align}
    \label{eq::discrete_lodi_orthogonality}
    \hat{\mat{A}}^- \left(\hat{\vec{U}}_h -  \hat{\vec{U}}^{n}_h \right) + \Delta t \hat{\mat{A}}^-\hat{\mat{P}}\hat{\AMP}^- +  \Delta t \hat{\mat{A}}^-\hat{\mat{P}}\hat{\AMP}^+ = \vec{0}.
\end{align}
On the one hand, we observe that outgoing characteristic amplitudes $\hat{\AMP}^+$, see \eqref{eq::characteristic_amplitudes}, vanish, since
\begin{align}
    \label{eq::filtering_plus}
    \hat{\mat{A}}^-\hat{\mat{P}}\hat{\AMP}^+ = \hat{\mat{P}}\hat{\EIG}^- \hat{\EIG}^+ \hat{\mat{P}}^{-1} \frac{\partial \vec{U}_h}{\partial \xi} = \mat{0},
\end{align}
due to the property $\hat{\EIG}^- \hat{\EIG}^+ = \mat{0}$ of the incoming and outgoing eigenvalues. On the other hand, the incoming characteristic
amplitudes $\hat{\AMP}^-$ can be expressed by the full set of characteristic amplitudes $\hat{\AMP}$, since the following identity holds
\begin{align}
    \label{eq::filtering_minus}
    \hat{\mat{A}}^-\hat{\mat{P}}\hat{\AMP}^- = \hat{\mat{P}} \EIG^- \hat{\AMP}^- = \hat{\mat{P}} \EIG^- \hat{\AMP} =  \hat{\mat{A}}^-\hat{\mat{P}}\hat{\AMP}.
\end{align}
Therefore, the system of equations \eqref{eq::discrete_lodi_orthogonality} filters outgoing waves, while incoming waves are retained.

Using \eqref{eq::filtering_plus} and \eqref{eq::filtering_minus} we rewrite the system as
\begin{align}
    \label{eq::discrete_lodi_orthogonality_full}
    \hat{\mat{A}}^- \left(\hat{\vec{U}}_h -  \hat{\vec{U}}^{n}_h + \Delta t \hat{\mat{P}}\hat{\AMP} \right) = \vec{0},
\end{align}
and in the spirit of \eqref{eq::modeled_amplitudes}, we model the full set of characteristic amplitudes
\begin{align}
    \label{eq::modeled_amplitudes_full}
    \hat{\AMP} \approx \tilde{\hat{\AMP}} := \hat{\mat{D}} \hat{\mat{P}}^{-1} \left( \hat{\vec{U}}_h - \vec{U}^-  \right), \qquad  \hat{\mat{D}} & \sim -\frac{\hat{\EIG}}{L},
\end{align}
since the correct incoming characteristic amplitudes are selected based on the eigenvalues of the incoming Jacobian $\hat{\mat{A}}^-$. This behavior can be particularly advantageous
at boundaries where flow inversion occurs, as the number of incoming characteristic amplitudes either increases or decreases,
depending on whether the boundary transitions to a subsonic inflow or outflow condition, respectively.

Inserting definition \eqref{eq::modeled_amplitudes_full} of the modeled
characteristic amplitudes $\tilde{\hat{\AMP}}$ into \eqref{eq::discrete_lodi_orthogonality_full}, yields
\begin{align}
    \label{eq::discrete_lodi_incoming}
    \hat{\mat{A}}^- \left(\hat{\vec{U}}_h -  \hat{\vec{U}}^{n}_h + \Delta t \hat{\mat{P}}\hat{\mat{D}} \hat{\mat{P}}^{-1} \left( \hat{\vec{U}}_h - \vec{U}^-  \right) \right) = \vec{0}.
\end{align}

Ultimately, we couple the above expression with the outgoing waves $\hat{\mat{A}}^+(\hat{\vec{U}}_h - \vec{U}_h)$,
and define the characteristic boundary operator $\hat{\vec{\Gamma}}_{\mathcal{C}}$ by combining equation \eqref{eq::discrete_lodi_incoming}, resulting in
\begin{align}
    \label{eq::cbc_boundary_operator}
    \hat{\vec{\Gamma}}_h = \hat{\vec{\Gamma}}_{\mathcal{C}} & := \underbrace{\hat{\mat{A}}^+(\hat{\vec{U}}_h - \vec{U}_h)}_{\text{outgoing}} -  \underbrace{\hat{\mat{A}}^- \left(\hat{\vec{U}}_h -  \hat{\vec{U}}^{n}_h + \Delta t \hat{\mat{P}}\hat{\mat{D}} \hat{\mat{P}}^{-1} \left( \hat{\vec{U}}_h - \vec{U}^-  \right) \right)}_{\text{incoming}}.
\end{align}
Note that the resulting characteristic boundary operator~\eqref{eq::cbc_boundary_operator} resembles the far-field operator~\eqref{eq::farfield}, with
a modified exterior state modeling the non-reflecting behavior.


Alternatively, following the same procedure as for the modified far-field operator \eqref{eq::farfield_mod}, we can express the
characteristic boundary operator \eqref{eq::cbc_boundary_operator} in a more stable variant
\begin{align}
    \label{eq::cbc_boundary_operator_mod}
    \hat{\vec{\Gamma}}_h = \hat{\vec{\Gamma}}^{\mat{Q}}_{\mathcal{C}} & := \hat{\vec{U}}_h - \hat{\mat{Q}}^+\vec{U}_h -  \hat{\mat{Q}}^- \left[\hat{\vec{U}}^{n}_h - \Delta t \hat{\mat{P}}\hat{\mat{D}} \hat{\mat{P}}^{-1} \left( \hat{\vec{U}}_h - \vec{U}^-  \right)\right].
\end{align}

\subsection{Implementation of the Navier-Stokes characteristic boundary conditions}
\label{sec::nscbc_hdg}
With the definition of characteristic boundary operator $\hat{\vec{\Gamma}}_{\mathcal{C}}$ in~\eqref{eq::cbc_boundary_operator},
the Navier-Stokes characteristic boundary conditions from Subsection~\ref{sec::NSCBC} are immediately recovered by appropriately
choosing the relaxation matrix $\hat{\mat{D}}$ and the target state $\vec{U}^-$.

\subsubsection{Subsonic outflow}
The incoming acoustic wave on a subsonic outflow is modeled by imposing the target pressure $p_\infty$ as in~\eqref{eq::L_acou_strikverda}, namely,
\begin{align}
    \label{eq::D_outflow_standard}
    \hat{\mat{D}}  = \begin{pmatrix}
                         \hat{D}_p & 0          & 0          & 0                      \\
                         0         & -\hat{u}_n & 0          & 0                      \\
                         0         & 0          & -\hat{u}_n & 0                      \\
                         0         & 0          & 0          & -(\hat{u}_n + \hat{c})
                     \end{pmatrix},                                                                    \qquad
    \vec{U}^-         = \vec{U}_{p_\infty} := \begin{pmatrix}
                                                  \rho_h \\[0.8ex] \rho \VEL_h \\[0.8ex] p_\infty/(\gamma - 1) + \rho_h |\VEL_h|^2/2
                                              \end{pmatrix}.
\end{align}

\subsubsection{Subsonic inflow}
For a subsonic inflow the incoming characteristic amplitudes can be modeled
by either imposing the target density $\rho_\infty$ and the target velocity $\VEL_\infty$ as in \eqref{eq::L_inflow_u_rho}, yielding
\begin{align}
    \hat{\mat{D}} = \begin{pmatrix}
                        \hat{D}_{u_n} & 0              & 0             & 0                      \\
                        0             & \hat{D}_{\rho} & 0             & 0                      \\
                        0             & 0              & \hat{D}_{u_t} & 0                      \\
                        0             & 0              & 0             & -(\hat{u}_n + \hat{c})
                    \end{pmatrix},                                                                                                       \qquad
    \vec{U}^-         = \vec{U}_{\rho \VEL_\infty} := \begin{pmatrix} \rho_\infty \\[0.8ex] \rho \VEL_\infty \\[0.8ex]  p_h/(\gamma - 1) + \rho_\infty |\VEL_\infty|^2/2
                                                      \end{pmatrix},
\end{align}
or by imposing the target density $\rho_\infty$, the target velocity
$\VEL_\infty$ and the target temperature $\theta_\infty$ as in \eqref{eq::L_inflow_u_T}, that is
\begin{align}
    \hat{\mat{D}} = \begin{pmatrix}
                        \hat{D}_{u_n} & 0                & 0             & 0                      \\
                        0             & \hat{D}_{\theta} & 0             & 0                      \\
                        0             & 0                & \hat{D}_{u_t} & 0                      \\
                        0             & 0                & 0             & -(\hat{u}_n + \hat{c})
                    \end{pmatrix}, \qquad
    \vec{U}^-         = \vec{U}_{\theta_\infty} := \rho(\theta_h) \begin{pmatrix} 1 \\ \VEL_\infty \\[0.8ex] \theta_h/\gamma + |\VEL_\infty|^2/2
                                                                  \end{pmatrix},
\end{align}
where the density,  $\rho(\theta_h)$, has been derived under the assumption of isentropic flow conditions
\begin{align}
    \label{eq::isentropic_density}
    \rho(\theta_h) & = \rho_\infty \left( \frac{\theta_h}{\theta_{\infty}} \right)^{\frac{1}{\gamma - 1}}.
\end{align}

\subsection{Implementation of the generalized characteristic relaxation boundary conditions}
\label{sec::grcbc_hdg}
Extending the concept of generalized relaxation boundary conditions (GRCBCs) from Subsection~\ref{sec::GRCBC} in a finite element framework
is more challenging due to the ambiguous definition of the first inner mesh point $\vec{x}_i$ to the boundary, required to determine the characteristic length $L_i$ in
the solely mesh-dependent relaxation matrix $\hat{\mat{D}}$, see \eqref{eq::relaxation_matrix_pirozzoli}. Numerical experiments conducted with element local diameter $L_i = h_T$ and $L_i = h_T/(k+1)$ consistently resulted in a hard boundary.

To address this issue, we propose a novel approach. Following the rationale of \cite{pirozzoliGeneralizedCharacteristicRelaxation2013},
any characteristic wave should fulfill a CFL condition to prevent stiffness in the numerical approximation, i.e., in matrix notation,
\begin{align}
    \label{eq::CFL_condition}
    \Delta t \mat{L}^{-1} |\hat{\EIG}| \leq \mat{C},  \qquad \mat{0} <  \mat{C} \leq \I,
\end{align}
where $\mat{L}$ and $\mat{C}$ are diagonal matrices and the symbols $<, \leq$ should be understood to act component-wise on the diagonal.

Let $\mat{L}$ be a matrix of different characteristic lengths. Assuming a given CFL-matrix $\mat{C}$ and
comparing \eqref{eq::CFL_condition} with \eqref{eq::modeled_amplitudes_full}, we can define the relaxation matrix~$\hat{\mat{D}}$ as
\begin{align}
    \label{eq::D_generalized}
    \hat{\mat{D}} & := \frac{1}{\Delta t} \mat{C} = \frac{1}{\Delta t} \begin{pmatrix}
                                                                           C_0 & 0   & 0   & 0   \\
                                                                           0   & C_1 & 0   & 0   \\
                                                                           0   & 0   & C_2 & 0   \\
                                                                           0   & 0   & 0   & C_3
                                                                       \end{pmatrix}.
\end{align}

The emerging relaxation matrix $\hat{\mat{D}}$ can be interpreted as if the first-order, one-sided upwind finite difference scheme
from \eqref{eq::FD-approximation} is performed at the characteristic level for different characteristic lengths $\mat{L}$, instead of a fixed length $L$, leading to
\begin{align}
    \label{eq::FD-approximation-grcbc}
    \frac{\partial \vec{W}}{\partial \xi} \stackrel{\text{FD approx.}}{\approx} -\mat{L}^{-1} \hat{\mat{P}}^{-1} (\hat{\vec{U}}_h - \vec{U}^-).
\end{align}
By assuming the characteristic lengths $\mat{L}$ to be equal to
\begin{align*}
    \mat{L}^{-1} = \frac{1}{\Delta t} \mat{C} |\hat{\EIG}|^{-1},
\end{align*}
and inserting \eqref{eq::FD-approximation-grcbc} into the definition of the characteristic amplitudes \eqref{eq::characteristic_amplitudes_definition}, we recover
the modeled characteristic amplitudes $\tilde{\hat{\AMP}}$ from \eqref{eq::modeled_amplitudes_full}
\begin{align}
    \label{eq::derivation_modeled_amplitudes_grcbc}
    \hat{\AMP} \approx \tilde{\hat{\AMP}} =  \frac{1}{\Delta t}\mat{C} \hat{\mat{P}}^{-1} (\hat{\vec{U}}_h - \vec{U}^-) = \hat{\mat{D}} \hat{\mat{P}}^{-1} (\hat{\vec{U}}_h - \vec{U}^-) .
\end{align}
In the derivation of \eqref{eq::derivation_modeled_amplitudes_grcbc} we exploited the fact, that due to the filtering property of the incoming Jacobian $\hat{\mat{A}}^-$ in \eqref{eq::filtering_minus}, it holds
\begin{align}
    - \hat{\mat{A}}^{-} \hat{\mat{P}} \hat{\EIG} |\hat{\EIG}|^{-1} \frac{1}{\Delta t}\mat{C}  =
    - \hat{\mat{P}} \hat{\EIG}^- \hat{\EIG} |\hat{\EIG}|^{-1} \frac{1}{\Delta t}\mat{C} =
    \hat{\mat{P}} \hat{\EIG}^- \frac{1}{\Delta t}\mat{C} =
    \hat{\mat{A}}^{-} \hat{\mat{P}} \frac{1}{\Delta t}\mat{C},
\end{align}
hence we can omit the term $-\hat{\EIG} |\hat{\EIG}|^{-1}$.

In comparison to the solely mesh-dependent relaxation matrix $\hat{\mat{D}}$ described in Subsection~\ref{sec::GRCBC},
the proposed method requires the user to specify the CFL-like relaxation factors $0 < C_i \leq 1$, for each
incoming characteristic wave. To facilitate the choice of relaxation factors $C_i$, we provide a prediction on the lower and
upper bounds:

\begin{itemize}
    \item \textbf{No relaxation} $\mat{C} \longrightarrow \mat{0}$:

        On the one hand, the lack of a relaxation
        matrix $\hat{\mat{D}}$ leads to a perfectly non-reflecting boundary for planar waves, as will be seen in the numerical experiments. On the other hand,
        as outlined in Subsection~\ref{sec::subsonic_outflow}, it may lead to a drift from the target state~$\vec{U}^-$.
    
    \item \textbf{Full relaxation} $\mat{C} \longrightarrow \I$
    
        In this scenario, the contributions of the incoming characteristic amplitudes are fully retained and have the same order of magnitude
        as the outgoing characteristic amplitudes, allowing the target state~$\vec{U}^-$ to be reached more quickly.
        However, depending on the target state $\vec{U}^-$, this can result in a highly reflecting boundary condition.
\end{itemize}


This approach recovers the common HDG boundary conditions, such as the subsonic outflow \eqref{eq::subsonic_outflow} and the far-field condition \eqref{eq::farfield},
as special cases. Assume $\hat{\mat{D}} = \frac{1}{\Delta t}\I$, then the boundary operator \eqref{eq::cbc_boundary_operator} reduces to
\begin{align*}
    \hat{\vec{\Gamma}}_h & = \hat{\mat{A}}^+(\hat{\vec{U}}_h - \vec{U}_h) -  \hat{\mat{A}}^- \left( 2 \hat{\vec{U}}_h -  \hat{\vec{U}}^{n}_h  - \vec{U}^- \right),
\end{align*}
which for $\hat{\vec{U}}^{n+1}_h \approx  \hat{\vec{U}}^{n}_h$ and target state $\vec{U}^- = \vec{U}_\infty$ resembles the far-field operator in \eqref{eq::farfield}.

The choice of the CFL factors $C_i$ can be seen as a trade-off between the two extremes of no relaxation and full relaxation, where the
improvement lies in the fact that the user can now adjust the relaxation factors $C_i$ to the specific problem at hand, allowing for a more flexible boundary condition.
Further, the inclusion of additional tangential and viscous contributions is possible and has a significant impact on the numerical solution, as demonstrated in the forthcoming numerical experiments.

\subsection{Multidimensional and viscous setting}
The extension to the complete Navier-Stokes equations follows the same line of thought. On the boundary, we discretize
the conservative system of equations~\eqref{eq::NSCBC-2d} by
\begin{align}
    \frac{\hat{\vec{U}}_h -  \hat{\vec{U}}^{n}_h}{\Delta t} + \hat{\mat{P}}\hat{\AMP}^- + \hat{\mat{P}}\hat{\AMP}^+ + \beta \hat{\mat{P}}\hat{\TRA}^- + \hat{\mat{P}}\hat{\TRA}^+    - \VIS = \vec{0}.
\end{align}
We rearrange terms and filter the equations by the incoming Jacobian $\hat{\mat{A}}^-$, yielding
\begin{align}
    \label{eq::discrete_nscbc_orthogonality}
    \hat{\mat{A}}^- \left(\hat{\vec{U}}_h -  \hat{\vec{U}}^{n}_h \right) + \Delta t \left[\hat{\mat{A}}^-\hat{\mat{P}}\hat{\AMP}^- + \hat{\mat{A}}^-\hat{\mat{P}}\hat{\AMP}^+ + \beta \hat{\mat{A}}^-\hat{\mat{P}}\hat{\TRA}^-  + \hat{\mat{A}}^-\hat{\mat{P}}\hat{\TRA}^+ - \hat{\mat{A}}^-\VIS \right]  = \vec{0}.
\end{align}
Similarly to the one-dimensional inviscid setting, we observe that outgoing tangential contributions $\hat{\TRA}^+$
vanish, while incoming contributions $\hat{\TRA}^-$ are retained, i.e.
\begin{align*}
    \hat{\mat{A}}^- \hat{\mat{P}}\hat{\TRA}^+  = \hat{\mat{A}}^-  \hat{\mat{Q}}^+  \hat{\mat{B}}\frac{\partial \vec{U}^-}{\partial \eta} = \vec{0},                                                                 \qquad
    \hat{\mat{A}}^- \hat{\mat{P}}\hat{\TRA}^-  = \hat{\mat{A}}^-  \hat{\mat{Q}}^- \hat{\mat{B}}\frac{\partial \vec{U}^-}{\partial \eta} = \hat{\mat{A}}^- \hat{\mat{B}}\frac{\partial \vec{U}^-}{\partial \eta} = \hat{\mat{A}}^-  \hat{\mat{P}}\hat{\TRA},
\end{align*}
hence we rewrite the system of equations \eqref{eq::discrete_nscbc_orthogonality} as
\begin{align}
    \hat{\mat{A}}^- \left(\hat{\vec{U}}_h -  \hat{\vec{U}}^{n}_h + \Delta t \left[\hat{\mat{P}}\hat{\AMP}  + \beta \hat{\mat{P}}\hat{\TRA} - \VIS \right]  \right) = \vec{0}.
\end{align}
We insert the modeled characteristic amplitudes \eqref{eq::modeled_amplitudes_full} and the conservative definition of the transverse contributions $\hat{\mat{P}}\hat{\TRA}$, to get
\begin{align}
    \hat{\mat{A}}^- \left(\hat{\vec{U}}_h -  \hat{\vec{U}}^{n}_h + \Delta t \left[ \hat{\mat{P}}\hat{\mat{D}} \hat{\mat{P}}^{-1} \left( \hat{\vec{U}}_h - \vec{U}^-  \right) +\beta \hat{\mat{B}}\frac{\partial \vec{U}^-}{\partial \eta}  - \VIS \right]  \right) = \vec{0}.
\end{align}

Ultimately, the definition of the boundary operator \eqref{eq::cbc_boundary_operator} remains unchanged, aside from the inclusion of additional tangential and viscous contributions
\begin{align}
    \label{eq::cbc_boundary_operator_nscbc}
    \hat{\vec{\Gamma}}_h = \hat{\vec{\Gamma}}_{\mathcal{C}}^{\TRA \VIS} & := \underbrace{\hat{\mat{A}}^+(\hat{\vec{U}}_h - \vec{U}_h)}_{\text{outgoing}} -  \underbrace{\hat{\mat{A}}^- \left(\hat{\vec{U}}_h -  \hat{\vec{U}}^{n}_h + \Delta t \left[ \hat{\mat{P}}\hat{\mat{D}} \hat{\mat{P}}^{-1} \left( \hat{\vec{U}}_h - \vec{U}^-  \right) +\beta \hat{\mat{B}}\frac{\partial \vec{U}^-}{\partial \eta}  - \VIS \right]  \right)}_{\text{incoming}} .
\end{align}


\subsubsection[Treatment of the tangential convective contributions]{Treatment of the tangential convective contributions $\TRA$}
Unlike the normal derivative, the tangential derivative along the boundary of the facet variable $\hat{\vec{U}}_h$ is defined, and therefore we discretize the tangential contributions as
\begin{align}
    \label{eq::tangential_convection}
    \hat{\mat{B}}\frac{\partial \vec{U}^-}{\partial \eta} & := \hat{\mat{B}}\frac{\partial \hat{\vec{U}}_h}{\partial \eta}.
\end{align}

\subsubsection[Treatment of the viscous contributions]{Treatment of the viscous contributions $\VIS$}
For the implementation of the viscous contributions $\VIS$, we exploit the fact that in our discretization the viscous flux $\vec{G}(\vec{U}, \EPS, \vec{\phi})$ is
linear in the conservative variables~$\vec{U}$, the rate-of-strain tensor~$\EPS$ and the temperature gradient~$\vec{\phi}$, see equation~\eqref{eq::fluxes}.
By regrouping the components of the rate-of-strain tensor $\EPS$ and the components of the temperature gradient $\vec{\phi}$ in a Voigt-like form
\begin{align*}
    \vec{\chi} :=  (\varepsilon_{xx}, \varepsilon_{xy}, \varepsilon_{yy}, \phi_{x}, \phi_{y} )^\T,
\end{align*}
we can express $\VIS$ in quasi-linear form as
\begin{align}
    \label{eq::quasi-linear-form-viscous}
    \VIS =  \frac{\partial (\vec{G} \vec{n})}{\partial \xi} + \frac{\partial (\vec{G} \vec{t})}{\partial \eta} =
    \frac{1}{\Re_{\infty}} \left[ \mat{A}_n^{\vec{U}} \frac{\partial \vec{U}}{\partial \xi} + \mat{A}_n^{\vec{\chi}} \frac{\partial \vec{\chi}}{\partial \xi} +
    \mat{A}_t^{\vec{U}} \frac{\partial \vec{U}}{\partial \eta} + \mat{A}_t^{\vec{\chi}} \frac{\partial \vec{\chi}}{\partial \eta} \right].
\end{align}
The viscous Jacobians in \eqref{eq::quasi-linear-form-viscous} are given as
\begin{subequations}
    \begin{alignat}{2}
        \mat{A}_n^{\vec{U}} & := \frac{\partial (\mat{G}\vec{n})}{\partial \vec{U}},  \qquad & \mat{A}_n^{\vec{\chi}} & := \frac{\partial (\mat{G}\vec{n})}{\partial \vec{\chi}}, \\
        \mat{A}_t^{\vec{U}} & := \frac{\partial (\mat{G}\vec{t})}{\partial \vec{U}},  \qquad & \mat{A}_t^{\vec{\chi}} & := \frac{\partial (\mat{G}\vec{t})}{\partial \vec{\chi}},
    \end{alignat}
\end{subequations}
where for an arbitrary unit vector $\vec{r} = \begin{bmatrix}r_x,& r_y \end{bmatrix}^\T$ it holds
\begin{align*}
    \mat{A}_r^{\vec{U}}    & := \frac{\partial (\mat{G}\vec{r})}{\partial \vec{U}}  = \frac{2}{\rho} \begin{pmatrix}
                                                                                                         0                            & 0                                           & 0                                           & 0 \\
                                                                                                         0                            & 0                                           & 0                                           & 0 \\
                                                                                                         0                            & 0                                           & 0                                           & 0 \\
                                                                                                         -\vec{r}^\T\EPS\vec{u} \cdot & \varepsilon_{xx} r_x + \varepsilon_{xy} r_y & \varepsilon_{xy} r_x + \varepsilon_{yy} r_y & 0
                                                                                                     \end{pmatrix}, \\
    \mat{A}_r^{\vec{\chi}} & := \frac{\partial (\mat{G}\vec{r})}{\partial \vec{\chi}}  = \begin{pmatrix}
                                                                                             0        & 0                   & 0        & 0                        & 0                        \\
                                                                                             2r_x     & 2r_y                & 0        & 0                        & 0                        \\
                                                                                             0        & 2r_x                & 2r_y     & 0                        & 0                        \\
                                                                                             2r_x u_x & 2r_x u_y + 2r_y u_x & 2r_y u_y & \frac{r_x}{\Pr_{\infty}} & \frac{r_y}{\Pr_{\infty}}
                                                                                         \end{pmatrix}.
\end{align*}
In order to modify the viscous contributions $\VIS$ to accommodate zero-normal gradient viscous conditions as described in Subsection~\ref{sec::NSCBC},
we reformulate $\partial_\xi \vec{\chi}$ in terms of normal-normal, normal-tangential and tangential-tangential components.
Together with the constant rotation matrix $\mat{R}$, emerging from the assumption of a local orthogonal basis, we express these as
\begin{align}
    \mat{R} := \begin{pmatrix}
                   n_x & t_x \\
                   n_y & t_y
               \end{pmatrix},    \qquad
    \begin{pmatrix}
        \varepsilon_{nn} & \varepsilon_{nt} \\
        \varepsilon_{nt} & \varepsilon_{tt}
    \end{pmatrix} =
    \mat{R}^\T \EPS \mat{R}, & \qquad
    \begin{pmatrix}
        \phi_{n} \\
        \phi_{t}
    \end{pmatrix} =
    \mat{R}^\T \vec{\phi}.
\end{align}
Therefore, by rearranging terms we can express the Cartesian derivatives in terms of normal-normal, normal-tangential and tangential-tangential derivatives, namely,
\begin{align}
    \begin{pmatrix}
        \frac{\partial \varepsilon_{xx}}{\partial \xi} & \frac{\partial \varepsilon_{xy}}{\partial \xi} \\
        \frac{\partial \varepsilon_{xy}}{\partial \xi} & \frac{\partial \varepsilon_{yy}}{\partial \xi}
    \end{pmatrix} =
    \mat{R}
    \begin{pmatrix}
        \frac{\partial \varepsilon_{nn}}{\partial \xi} & \frac{\partial \varepsilon_{nt}}{\partial \xi} \\
        \frac{\partial \varepsilon_{nt}}{\partial \xi} & \frac{\partial \varepsilon_{tt}}{\partial \xi}
    \end{pmatrix}
    \mat{R}^\T, \qquad
    \begin{pmatrix}
        \frac{\partial \phi_x}{\partial \xi} \\
        \frac{\partial \phi_y}{\partial \xi}
    \end{pmatrix} =
    \mat{R}
    \begin{pmatrix}
        \frac{\partial \phi_n}{\partial \xi} \\
        \frac{\partial \phi_t}{\partial \xi}
    \end{pmatrix},
\end{align}
or in explicit form
\begin{align}
    \label{eq::normal_chi}
    \frac{\partial \vec{\chi}}{\partial{\xi}} =
    \begin{pmatrix}
        \frac{\partial \varepsilon_{xx}}{\partial \xi} \\[0.8ex]
        \frac{\partial \varepsilon_{xy}}{\partial \xi} \\[0.8ex]
        \frac{\partial \varepsilon_{yy}}{\partial \xi} \\[0.8ex]
        \frac{\partial \phi_{x}}{\partial \xi}         \\[0.8ex]
        \frac{\partial \phi_{y}}{\partial \xi}
    \end{pmatrix} =
    \begin{pmatrix}
        \frac{\partial \varepsilon_{nn}}{\partial \xi} n_x^2 + 2\frac{\partial \varepsilon_{nt}}{\partial \xi} n_x t_x + \frac{\partial \varepsilon_{tt}}{\partial \xi} t_x^2                \\[0.8ex]
        \frac{\partial \varepsilon_{nn}}{\partial \xi} n_x n_y + \frac{\partial \varepsilon_{nt}}{\partial \xi} (t_y n_x + t_x n_y) + \frac{\partial \varepsilon_{tt}}{\partial \xi} t_x t_y \\[0.8ex]
        \frac{\partial \varepsilon_{nn}}{\partial \xi} n_y^2 + \frac{\partial \varepsilon_{nt}}{\partial \xi} n_y t_y + \frac{\partial \varepsilon_{tt}}{\partial \xi} t_y^2                 \\[0.8ex]
        \frac{\partial \phi_{n}}{\partial \xi} n_x + \frac{\partial \phi_{t}}{\partial \xi} t_x                                                                                              \\[0.8ex]
        \frac{\partial \phi_{n}}{\partial \xi} n_y + \frac{\partial \phi_{t}}{\partial \xi} t_y
    \end{pmatrix}.
\end{align}

To impose the viscous conditions, we modify \eqref{eq::normal_chi} by eliminating the
appropriate  normal-normal, normal-tangential and tangential-tangential derivative of the rate-of-strain tensor
and/or the normal derivative of the temperature gradient, which for a constant viscosity $\mu$ and constant heat conductivity $\kappa$
is equivalent to imposing the deviatoric stress tensor gradient $\TAU$ and heat flux gradient $\HEAT$, respectively. For example, imposing
\begin{align*}
    \frac{\partial(\vec{t}^\T \TAU \vec{n})}{\partial \xi}  = 0, \qquad \frac{\partial(\HEAT \cdot \vec{n})}{\partial \xi} = 0,
\end{align*}
on a subsonic outflow (see Subsection~\ref{sec::subsonic_outflow}), yields the following modified normal derivative
\begin{align}
    \label{eq::normal_chi_modified}
    \frac{\partial \tilde{\vec{\chi}}}{\partial{\xi}} =
    \begin{pmatrix}
        \frac{\partial \varepsilon_{nn}}{\partial \xi} n_x^2 +  \frac{\partial \varepsilon_{tt}}{\partial \xi} t_x^2    \\[0.8ex]
        \frac{\partial \varepsilon_{nn}}{\partial \xi} n_x n_y + \frac{\partial \varepsilon_{tt}}{\partial \xi} t_x t_y \\[0.8ex]
        \frac{\partial \varepsilon_{nn}}{\partial \xi} n_y^2  + \frac{\partial \varepsilon_{tt}}{\partial \xi} t_y^2    \\[0.8ex]
        \frac{\partial \phi_{t}}{\partial \xi} t_x                                                                      \\[0.8ex]
        \frac{\partial \phi_{t}}{\partial \xi} t_y
    \end{pmatrix}.
\end{align}

Ultimately, we replace $\VIS$ in \eqref{eq::cbc_boundary_operator_nscbc} by the modified viscous contributions
\begin{align}
    \tilde{\VIS} = \frac{1}{\Re_{\infty}} \left[ \mat{A}_n^{\vec{U}_h} \frac{\partial \vec{U}_h}{\partial \xi} +
    \mat{A}_n^{\vec{\chi}_h} \frac{\partial \tilde{\vec{\chi}}_h}{\partial{\xi}} + \mat{A}_t^{\vec{U}_h} \frac{\partial \vec{U}_h}{\partial \eta} +
    \mat{A}_t^{\vec{\chi}_h} \frac{\partial \vec{\chi}_h}{\partial \eta} \right],
\end{align}
where all derivatives are evaluated at the interior state $\vec{U}_h$ and $\vec{\chi}_h$.

To conclude this section, we provide a summary of all introduced target states and boundary conditions in Table~\ref{tab::states_summary} and Table~\ref{tab::boundary_conditions_summary},
respectively. Recall that the difference between the Navier-Stokes characteristic boundary conditions (NSCBCs) of Subsection~\ref{sec::nscbc_hdg}
and the generalized relaxation boundary conditions (GRCBCs) of Subsection~\ref{sec::grcbc_hdg} solely lies in the definition of
the relaxation matrix $\hat{\mat{D}}$. Furthermore, the definition of the characteristic boundary operator (CBC) in Table~\ref{tab::boundary_conditions_summary} has been derived assuming an
implicit Euler time marching scheme, see Remark~\ref{rem::cbc_time_discretization}. Adapting the characteristic boundary operator (CBC) to a different time discretization scheme is straightforward and can be accomplished by following the same procedure outlined in this section.

\begin{table}[h]
    \centering
    \begin{tabularx}{\textwidth}{TML}
                                        & Target state    & $\vec{U}^-$                                                                                                                                                      \\
        \midrule
        $\vec{U}_{p_\infty}$         & Pressure state  & $\begin{pmatrix}\rho_h,  & \rho \VEL_h,  & p_\infty/(\gamma - 1) + \rho_h |\VEL_h|^2/2\end{pmatrix}^\T$                                                          \\[1.5ex]
        $\vec{U}_{\infty}$           & Far-field state & $\begin{pmatrix} \rho_\infty, &\rho_\infty \VEL_\infty, & \rho_\infty E_\infty\end{pmatrix}^\T$                                                                  \\[1.5ex]
        $\vec{U}_{\rho \VEL_\infty}$ & Mass state      & $\begin{pmatrix} \rho_\infty, &\rho \VEL_\infty,  & p_h/(\gamma - 1) + \rho_\infty |\VEL_\infty|^2/2\end{pmatrix}^\T$                                            \\[1.5ex]
        $\vec{U}_{\theta_\infty}$    & Entropy State   & $\rho_\infty \left(  \theta_h/\theta_{\infty} \right)^{1/(\gamma - 1)} \begin{pmatrix}1,  & \VEL_\infty, & \theta_h/\gamma+ |\VEL_\infty|^2/2 \end{pmatrix}^\T $
    \end{tabularx}
    \caption{Summary of target states}
    \label{tab::states_summary}
\end{table}

\begin{table}[ht]
    \centering
    \begin{tabularx}{\textwidth}{TML}
                                                       & Boundary type    & $\hat{\Gamma}_h$                                                                                                                                                                                                                                                                                                             \\
        \midrule
        $\hat{\vec{\Gamma}}_{\textit{ad}}$             & Adiabatic wall   & $\begin{pmatrix} \hat{\rho}_h - \rho_h, & \hat{\rho \VEL}_h,  & 1/(\Re_\infty \Pr_\infty)\vec{\phi}_h \cdot \vec{n} \end{pmatrix}^\T - \mat{S}_d \begin{pmatrix} 0, & \bm{0}, & \rho E_h - \hat{\rho E}_h\end{pmatrix}^\T$                                                                                                   \\[1.5ex]
        $\hat{\vec{\Gamma}}_{\textit{iso}}$            & Isothermal wall  & $\hat{\vec{U}}_h -  \begin{pmatrix} \rho_h, & \vec{0}, & \rho_h \theta_{\text{wall}}/\gamma \end{pmatrix}^\T$                                                                                                                                                                                                    \\[1.5ex]
        $\hat{\vec{\Gamma}}_{\textit{inv}}$            & Inviscid wall    & $\hat{\vec{U}}_h -\begin{pmatrix} \rho_h, & (\I - \vec{n} \otimes \vec{n}) \rho \VEL_h, & \rho E_h \end{pmatrix}^\T $                                                                                                                                                                                                        \\[1.5ex]
        $\hat{\vec{\Gamma}}_{p_\infty}$                & Subsonic outflow & $\hat{\vec{U}}_h - \vec{U}_{p_\infty}$                                                                                                                                                                                                                                                                                       \\[1.5ex]
        $\hat{\vec{\Gamma}}_{\infty}$                  & Far-field        & $\hat{\mat{A}}^+(\hat{\vec{U}}_h - \vec{U}_h) -  \hat{\mat{A}}^-(\hat{\vec{U}}_h - \vec{U}_\infty)$                                                                                                                                                                                                                          \\[1.5ex] 
        $\hat{\vec{\Gamma}}_{\mathcal{C}}$             &                  & $\hat{\mat{A}}^+(\hat{\vec{U}}_h - \vec{U}_h) -  \hat{\mat{A}}^- \left(\hat{\vec{U}}_h -  \hat{\vec{U}}^{n}_h + \Delta t \left[ \hat{\mat{P}}\hat{\mat{D}} \hat{\mat{P}}^{-1} \left( \hat{\vec{U}}_h - \vec{U}^-  \right) \right]  \right)$                                                                                  \\[1.5ex]
        $\hat{\vec{\Gamma}}_{\mathcal{C}}^{\TRA}$      & CBC              & $\hat{\mat{A}}^+(\hat{\vec{U}}_h - \vec{U}_h) -  \hat{\mat{A}}^- \left(\hat{\vec{U}}_h -  \hat{\vec{U}}^{n}_h + \Delta t \left[ \hat{\mat{P}}\hat{\mat{D}} \hat{\mat{P}}^{-1} \left( \hat{\vec{U}}_h - \vec{U}^-  \right) +\beta \hat{\mat{B}}\frac{\partial \hat{\vec{U}}}{\partial \eta} \right]  \right)$                 \\[1.5ex]
        $\hat{\vec{\Gamma}}_{\mathcal{C}}^{\TRA \VIS}$ &                  & $\hat{\mat{A}}^+(\hat{\vec{U}}_h - \vec{U}_h) -  \hat{\mat{A}}^- \left(\hat{\vec{U}}_h -  \hat{\vec{U}}^{n}_h + \Delta t \left[ \hat{\mat{P}}\hat{\mat{D}} \hat{\mat{P}}^{-1} \left( \hat{\vec{U}}_h - \vec{U}^-  \right) +\beta \hat{\mat{B}}\frac{\partial \hat{\vec{U}}}{\partial \eta}  - \tilde{\VIS} \right]  \right)$
    \end{tabularx}
    \caption{Summary of boundary conditions}
    \label{tab::boundary_conditions_summary}
\end{table}

\section{Numerical experiments}
\label{sec::experiments}
The implementation of the characteristic boundary conditions described in Section~\ref{sec::cbc_for_hdg} has been
evaluated by four benchmarks in a two-dimensional setting, inspired by~\cite{pirozzoliGeneralizedCharacteristicRelaxation2013,shehadiPolynomialcorrectionNavierStokesCharacteristic2024a}. Every numerical experiment includes a comparison of the
numerical solution against the exact solution, if available; otherwise, it is compared against a reference solution. The latter is computed in a larger reference domain $\Omega_{r}$
equipped with appropriate non-reflecting closures at the boundaries and, if deemed necessary, with sponge layers to dampen incoming disturbances.
One of the metrics we use are the scaled $L_2(\Omega)$-errors of the discrete pressure $p_h$, with respect to the reference/exact solution $p_r$, that is
\begin{align*}
    e_p & := \frac{\| p_h - p_r \|}{\| p_r \|}.
\end{align*}
Additionally, the results are compared to common boundary conditions in the HDG framework, namely the subsonic outflow (SO) boundary condition \eqref{eq::subsonic_outflow}
and the far-field (FF) boundary condition \eqref{eq::farfield}. To distinguish between the Navier-Stokes characteristic boundary conditions (NSCBCs) of Subsection~\ref{sec::nscbc_hdg} and the 
generalized characteristic relaxation boundary conditions (GRCBCs) of Subsection~\ref{sec::grcbc_hdg}, the different relaxation matrices $\hat{\mat{D}}$, the different target states $\vec{U}^-$, the inclusion of tangential
convection contributions $\TRA^-$ and the viscous terms $\VIS$ we introduce the following notation:
\begin{itemize}
    \item The type of characteristic boundary condition is denoted by the abbreviation NSCBC or GRCBC.
    \item The value of the user-required inputs of the relaxation matrices $\hat{\mat{D}}$, is displayed as subscript, e.g., GRCBC$_{0.01}$.
          For the sake of simplicity, we assumed a single user-value for all components of the relaxation matrix i.e. for the GRCBCs we solely imposed $C_i = C$, see \eqref{eq::D_generalized}, while for the NSCBCs we set $\sigma_i = \sigma$ and $L=1$, see \eqref{eq::L_acou_strikverda_relaxation} and \eqref{eq::L_inflow_u_rho_relaxation}.
    \item  For subsonic outflows we introduced the target states $\vec{U}_{p_{\infty}}$ and $\vec{U}_{\infty}$, while for
          subsonic inflows we introduced $\vec{U}_{\infty}$, $\vec{U}_{\rho \VEL_{\infty}}$ and $\vec{U}_{\theta_{\infty}}$. The type of target state, $\vec{U}^-$,
          is denoted by the corresponding subscript, e.g., GRCBC$_{0.01 p_{\infty}}$ for the target state $\vec{U}_{p_{\infty}}$.
    \item If tangential convective contributions $\TRA$ are included, the symbol $\TRA$ is denoted in the superscript, e.g., GRCBC$_{0.01 p_{\infty}}^{\TRA}$. The tangential relaxation parameter $\beta$ is set to $\Ma_\infty$ in all cases.
    \item If viscous contributions $\VIS$ are included, the symbol $\VIS$ is denoted in the superscript, e.g., GRCBC$_{0.01 p_{\infty}}^{\TRA \VIS}$.
\end{itemize}



A calorically perfect gas with specific heat ratio $\gamma = 1.4$, Prandtl number $\Pr_\infty = 0.72$ and with constant viscosity $\mu=1$ is considered.
Unless specified otherwise, the spatial discretization employs a polynomial order $k=4$, time-stepping is performed using an implicit BDF2 scheme with a suitable time
step $\Delta t$, and the convective and viscous stabilization matrices are defined as $\mat{S}_c = \hat{\mat{A}}^+$ and
$\mat{S}_d = \diag(0, \vec{1}, 1/\Pr_\infty)/\Re_\infty$, respectively.

The numerical experiments have been conducted using the open-source python extension \texttt{DreAm}, which is
built on the open-source finite element library \texttt{NGSolve}~\cite{schoberlC++11ImplementationFinite2014}. Simulation data is publicly 
available in the \texttt{Zenodo} repository~\cite{ellmenreichSimulationDataPaper2025}.

\subsection{Planar acoustic pulse}
\label{sec::experiments::planar_acoustic_pulse}
\graphicspath{{./figures/planar_acoustic_pulse}}
\pgfplotsset{table/search path={./figures/planar_acoustic_pulse/data}}

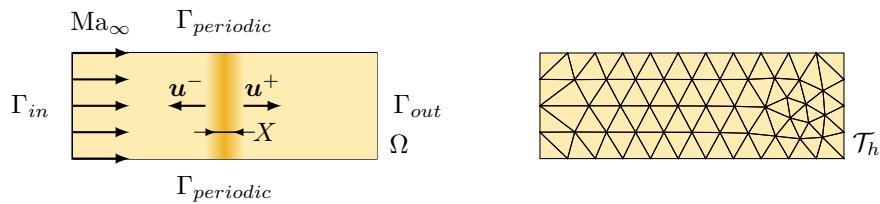
\begin{figure}[htb]
    \centering
    \tikzsetnextfilename{planar_acoustic_pulse_domain}
\begin{tikzpicture}[scale=5, >=latex]
    \draw[thick] (-0.4, -0.14) rectangle (0.4, 0.14);
    \fill[amber!30] (-0.4, -0.14) rectangle (0.4, 0.14);

    \node at (-0.51, 0) {$\Gamma_{in}$};
    \node at (0.51, 0) {$\Gamma_{out}$};
    \node at (0, 0.22) {$\Gamma_{periodic}$};
    \node at (0, -0.22) {$\Gamma_{periodic}$};

    \draw[thick, ->] (-0.4, -0.14) -- (-0.25, -0.14);
    \draw[thick, ->] (-0.4, -0.07) -- (-0.25, -0.07) ;
    \draw[thick, ->] (-0.4, 0) --     (-0.25, 0);
    \draw[thick, ->] (-0.4, 0.07) --  (-0.25, 0.07);
    \draw[thick, ->] (-0.4, 0.14) --  (-0.25, 0.14);
    \node at (-0.325, 0.22) {$\Ma_\infty$};

    \shade[left color=amber!30, right color=amber!30, middle color=matyellow] (-0.05, -0.14) rectangle (0.05, 0.14);
    \draw[thick, ->] (0.05, 0) --  (0.15, 0) node[midway, above] {$\VEL^+$};
    \draw[thick, ->] (-0.05, 0) --  (-0.15, 0) node[midway, above] {$\VEL^-$};

    \draw[>-<] (-0.05, -0.07) -- (0.05,-0.07) node[right] {$X$};
    \draw[] (-0.08, -0.07) -- (0.08,-0.07);

    \node at (0.46, -0.1) {$\Omega$};

    \begin{scope}[xshift=35]

        \fill[amber!30] (-0.4, -0.14) rectangle (0.4, 0.14);
        \node at (0.46, -0.1) {$\mesh$};

        \pgfplotstableread[col sep=comma]{mesh.csv}\mesh
        \pgfplotstablegetrowsof{\mesh}
        \pgfmathsetmacro{\rows}{\pgfplotsretval-1}

        \foreach \row in {0,...,\rows} {
                \pgfplotstablegetelem{\row}{x0}\of\mesh
                \let\xi=\pgfplotsretval
                \pgfplotstablegetelem{\row}{y0}\of\mesh
                \let\yi=\pgfplotsretval
                \pgfplotstablegetelem{\row}{x1}\of\mesh
                \let\xj=\pgfplotsretval
                \pgfplotstablegetelem{\row}{y1}\of\mesh
                \let\yj=\pgfplotsretval
                \pgfplotstablegetelem{\row}{x2}\of\mesh
                \let\xk=\pgfplotsretval
                \pgfplotstablegetelem{\row}{y2}\of\mesh
                \let\yk=\pgfplotsretval
                \draw (\xi,\yi) -- (\xj,\yj) -- (\xk,\yk) -- cycle;
            }








    \end{scope}

\end{tikzpicture}
    \caption[]{Representation of the domain $\Omega$ with corresponding boundaries $\Gamma$ and the triangulation $\mesh$ with mesh size $h = 0.035$.}
    \label{fig::pulse_1d_domain}
\end{figure}

The first benchmark consists of an inviscid one-dimensional downstream $\VEL^+$/upstream $\VEL^-$ propagating linear acoustic pulse, similar to \cite{shehadiPolynomialcorrectionNavierStokesCharacteristic2024a}, 
in a domain $\Omega = (-4X, 4X) \times (-0.14, 0.14)$ depicted in Figure~\ref{fig::pulse_1d_domain}. In order to mimic a one-dimensional domain, the upper and lower boundaries $\Gamma_{periodic}$ are coupled periodically.

The acoustic pulse is initialized at $t=0$ with the initial conditions given as
\begin{align*}
    \rho_0 = \rho_\infty, \qquad \VEL^\pm_0 = \VEL_\infty \left( 1 \pm \frac{p_\infty}{ \rho_\infty c_\infty} \frac{\alpha}{| \VEL_\infty |} \exp\left(-x^2/X^2\right)\right), \qquad p_0  = p_\infty \left(1 + \alpha \exp\left(-x^2/X^2\right) \right),
\end{align*}
where $X=0.1$ denotes the extent of the acoustic pulse and $\alpha$ the pulse strength. The constant uniform flow 
$\VEL_\infty$ ensures the correct physical behavior at the subsonic inflow boundary $\Gamma_{in}$ and subsonic outflow boundary $\Gamma_{out}$.
Note that we use acoustic scaling, see \eqref{eq::acoustic-scaling}, with dimensionless quantities $\rho_\infty=1$, $|\VEL_\infty| = \Ma_\infty = 0.03$, $c_\infty = 1$ and $p_\infty = 1/\gamma$.

Every simulation is performed up to a dimensionless time $t_{end}=5$ with time step $\Delta t = 2\cdot 10^{-3}$. As a reference solution we use a
computed solution in a larger domain $\Omega_{r} =  (-14X, 14X) \times (-0.14, 0.14)$ in which the acoustic pressure pulse propagates
undisturbed through the domain during the time interval $(0, t_{end}]$ without impinging on the boundaries. In addition to the scaled $L_2(\Omega)$-errors of the pressure $e_p$,
we rely on the mean incoming acoustic $\overline{w}^-$ and outgoing $\overline{w}^+$ acoustic characteristics at the boundaries $\Gamma_{in}$ and $\Gamma_{out}$ to assess the quality of the numerical solution.
These are defined as
\begin{align*}
    w^\mp := (p_h - p_\infty) \mp \rho_h c_h (\VEL_h - \VEL_\infty) \cdot \vec{n}, \qquad \overline{w}^\mp  := \frac{1}{|\Gamma|} \int_{\Gamma} w^\mp \, d\bm{s},
\end{align*}
respectively.  Further, for the upstream propagating pulse, we consider the mean entropy wave 
\begin{align*}
    w^s := c_h^2 (\rho_h - \rho_\infty) - (p_h - p_\infty) , \qquad \overline{w}^s  := \frac{1}{|\Gamma_{in}|} \int_{\Gamma_{in}} w^s\, d\bm{s},
\end{align*}
at the boundary $\Gamma_{in}$, since this characteristic wave is determined from upstream in the case of subsonic inflow.

In order to investigate the effect of the relaxation matrix $\hat{\mat{D}}$ on the reflection of the acoustic pulse, we consider three different
relaxation coefficients $C = \left[ 0.1, 0.01, 0.0 \right]$ for the GRCBCs and $\sigma = \left[ 1.0, \sigma_{opt}, 0.0 \right]$ for the NSCBCs, 
where $\sigma_{opt} = 0.28$ is approximately the optimal relaxation parameter for a subsonic outflow, as described in Subsection~\ref{sec::subsonic_outflow}.

\begin{figure}[htb]
    \centering
    \tikzsetnextfilename{planar_acoustic_pulse_w_alpha0.001_right}
\begin{tikzpicture}
    \pgfplotsset{every axis/.append style={
                width=4cm,
                height=4cm,
                mark repeat = 20,
                mark size = 1pt,
                mark options={solid},
                xmin=0,
                xmax=1.0,
                max space between ticks=20,
                semithick,
                xtick pos=bottom,
                xmajorgrids=true, ymajorgrids=true,
                axis line style={thin},
                tick label style = {font=\tiny},
                ytick pos=left,
                xtick = {0,0.25,0.5,0.75,1.0}},
        legend style={font=\tiny},
        legend cell align={left},
    }

    \begin{groupplot}[group style={group size=4 by 2, x descriptions at=edge bottom, vertical sep=0.5cm, horizontal sep=0.5cm}]

        \nextgroupplot[ymin=-0.0002, ymax=0.0015, ylabel = $\overline{w}^+$, legend to name=bcs]
        \addplot[black,    mark phase = 0, mark=square*, mark options={solid}] table [y=w+_right, x=t]{grcbc_farfield_reference_Ma0.03_alpha0.001_wave_right_dt0.002.dat}; \addlegendentry{Reference};
        \addplot[matblue,  mark phase = 8, mark=square*]  table [y=w+_right, x=t]{farfield_inflow_and_outflow_Ma0.03_alpha0.001_wave_right_dt0.002.dat}; \addlegendentry{FF}; \label{leg::w-::ff}
        \addplot[matred, mark phase = 15, mark=square*] table [y=w+_right, x=t]{farfield_inflow_and_pressure_outflow_Ma0.03_alpha0.001_wave_right_dt0.002.dat};\addlegendentry{SO}; \label{leg::w-::so}

        \addlegendimage{matgreen, mark=square*, mark options={solid}};
        \addlegendentry{GRCBC$_{0.1\infty}$};  \label{leg::w-::grcbc_ff}
        \addlegendimage{matgreen, mark=*, dashed, mark options={solid}};
        \addlegendentry{GRCBC$_{0.01 \infty}$};
        \addlegendimage{matgreen, mark=diamond*, dash dot, mark options={solid}};
        \addlegendentry{GRCBC$_{0.0\infty}$};

        \addlegendimage{matyellow, mark=square*, mark options={solid}};
        \addlegendentry{GRCBC$_{0.1p_{\infty}}$}; \label{leg::w-::grcbc_so}
        \addlegendimage{matyellow, mark=*, dashed, mark options={solid}};
        \addlegendentry{GRCBC$_{0.01p_{\infty}}$};
        \addlegendimage{matyellow, mark=diamond*, dash dot, mark options={solid}};
        \addlegendentry{GRCBC$_{0.0p_{\infty}}$};

        \addlegendimage{matorange, mark=square*, mark options={solid}};
        \addlegendentry{NSCBC$_{1.0p_{\infty}}$}; \label{leg::w-::nscbc_so}
        \addlegendimage{matorange, mark=*, dashed, mark options={solid}};
        \addlegendentry{NSCBC$_{0.28p_{\infty}}$};
        \addlegendimage{matorange, mark=diamond*, dash dot, mark options={solid}};
        \addlegendentry{NSCBC$_{0.01p_{\infty}}$};

        \nextgroupplot[ymin=-0.0002, ymax=0.0015, yticklabel=\empty, scaled y ticks = false]
        \addplot[black,  mark phase = 0, mark=square*, mark options={solid}] table [y=w+_right, x=t]{grcbc_farfield_reference_Ma0.03_alpha0.001_wave_right_dt0.002.dat};
        \addplot[matgreen, mark phase = 6, mark=square*]           table [y=w+_right_0.1, x=t]{grcbc_farfield_inflow_and_outflow_Ma0.03_alpha0.001_wave_right_dt0.002.dat};
        \addplot[matgreen, mark phase = 11, mark=*, dashed]   table [y=w+_right_0.01, x=t]{grcbc_farfield_inflow_and_outflow_Ma0.03_alpha0.001_wave_right_dt0.002.dat};
        \addplot[matgreen, mark phase = 16, mark=diamond*, dash dot] table [y=w+_right_0.0, x=t]{grcbc_farfield_inflow_and_outflow_Ma0.03_alpha0.001_wave_right_dt0.002.dat};

        \nextgroupplot[ymin=-0.0002, ymax=0.0015, yticklabel=\empty, scaled y ticks = false]
        \addplot[black,     mark phase = 0,  mark=square*, mark options={solid}] table [y=w+_right, x=t]{grcbc_farfield_reference_Ma0.03_alpha0.001_wave_right_dt0.002.dat};
        \addplot[matyellow, mark phase = 6,  mark=square*]           table [y=w+_right_0.1, x=t]{grcbc_farfield_inflow_and_pressure_outflow_Ma0.03_alpha0.001_wave_right_dt0.002.dat};
        \addplot[matyellow, mark phase = 11, mark=*, dashed]   table [y=w+_right_0.01, x=t]{grcbc_farfield_inflow_and_pressure_outflow_Ma0.03_alpha0.001_wave_right_dt0.002.dat};
        \addplot[matyellow, mark phase = 16, mark=diamond*, dash dot] table [y=w+_right_0.0, x=t]{grcbc_farfield_inflow_and_pressure_outflow_Ma0.03_alpha0.001_wave_right_dt0.002.dat};

        \nextgroupplot[ymin=-0.0002, ymax=0.0015, yticklabel=\empty, scaled y ticks = false]
        \addplot[black,    mark phase = 0,  mark=square*, mark options={solid}] table [y=w+_right, x=t]{grcbc_farfield_reference_Ma0.03_alpha0.001_wave_right_dt0.002.dat};
        \addplot[matorange, mark phase = 6,  mark=square*]           table [y=w+_right_1.0, x=t]{nscbc_farfield_inflow_and_pressure_outflow_Ma0.03_alpha0.001_wave_right_dt0.002.dat};
        \addplot[matorange, mark phase = 11, mark=*, dashed]   table [y=w+_right_0.28, x=t]{nscbc_farfield_inflow_and_pressure_outflow_Ma0.03_alpha0.001_wave_right_dt0.002.dat};
        \addplot[matorange, mark phase = 16, mark=diamond*, dash dot] table [y=w+_right_0.01, x=t]{nscbc_farfield_inflow_and_pressure_outflow_Ma0.03_alpha0.001_wave_right_dt0.002.dat};
        \coordinate (top) at (rel axis cs:0,1);

        \nextgroupplot[ymin=-0.0015, ymax=0.0002,  ylabel = $\overline{w}^-$, xlabel = $t$]
        \addplot[black,    mark phase = 0, mark=square*, mark options={solid}] table [y=w-_right, x=t]{grcbc_farfield_reference_Ma0.03_alpha0.001_wave_right_dt0.002.dat};
        \addplot[matblue,  mark phase = 8, mark=square*] table [y=w-_right, x=t]{farfield_inflow_and_outflow_Ma0.03_alpha0.001_wave_right_dt0.002.dat};
        \addplot[matred, mark phase = 15, mark=square*] table [y=w-_right, x=t]{farfield_inflow_and_pressure_outflow_Ma0.03_alpha0.001_wave_right_dt0.002.dat};

        \nextgroupplot[ymin=-0.0015, ymax=0.0002, yticklabel=\empty, scaled y ticks = false, xlabel = $t$]
        \addplot[black,  mark phase = 0,  mark=square*, mark options={solid}] table [y=w-_right, x=t]{grcbc_farfield_reference_Ma0.03_alpha0.001_wave_right_dt0.002.dat};
        \addplot[matgreen, mark phase = 6,  mark=square*]           table [y=w-_right_0.1, x=t]{grcbc_farfield_inflow_and_outflow_Ma0.03_alpha0.001_wave_right_dt0.002.dat};
        \addplot[matgreen, mark phase = 11, mark=*, dashed]   table [y=w-_right_0.01, x=t]{grcbc_farfield_inflow_and_outflow_Ma0.03_alpha0.001_wave_right_dt0.002.dat};
        \addplot[matgreen, mark phase = 16, mark=diamond*, dash dot] table [y=w-_right_0.0, x=t]{grcbc_farfield_inflow_and_outflow_Ma0.03_alpha0.001_wave_right_dt0.002.dat};

        \nextgroupplot[ymin=-0.0015, ymax=0.0002, yticklabel=\empty, scaled y ticks = false, xlabel = $t$]
        \addplot[black,     mark phase = 0, mark=square*, mark options={solid}] table [y=w-_right, x=t]{grcbc_farfield_reference_Ma0.03_alpha0.001_wave_right_dt0.002.dat};
        \addplot[matyellow, mark phase = 6, mark=square*]           table [y=w-_right_0.1, x=t]{grcbc_farfield_inflow_and_pressure_outflow_Ma0.03_alpha0.001_wave_right_dt0.002.dat};
        \addplot[matyellow, mark phase = 11,mark=*, dashed]   table [y=w-_right_0.01, x=t]{grcbc_farfield_inflow_and_pressure_outflow_Ma0.03_alpha0.001_wave_right_dt0.002.dat};
        \addplot[matyellow, mark phase = 16,mark=diamond*, dash dot] table [y=w-_right_0.0, x=t]{grcbc_farfield_inflow_and_pressure_outflow_Ma0.03_alpha0.001_wave_right_dt0.002.dat};

        \nextgroupplot[ymin=-0.0015, ymax=0.0002, yticklabel=\empty, scaled y ticks = false, xlabel = $t$]
        \addplot[black,    mark phase = 0,  mark=square*, mark options={solid}] table [y=w-_right, x=t]{grcbc_farfield_reference_Ma0.03_alpha0.001_wave_right_dt0.002.dat};
        \addplot[matorange, mark phase = 6,  mark=square*]           table [y=w-_right_1.0, x=t]{nscbc_farfield_inflow_and_pressure_outflow_Ma0.03_alpha0.001_wave_right_dt0.002.dat};
        \addplot[matorange, mark phase = 11, mark=*, dashed]   table [y=w-_right_0.28, x=t]{nscbc_farfield_inflow_and_pressure_outflow_Ma0.03_alpha0.001_wave_right_dt0.002.dat};
        \addplot[matorange, mark phase = 16, mark=diamond*, dash dot] table [y=w-_right_0.01, x=t]{nscbc_farfield_inflow_and_pressure_outflow_Ma0.03_alpha0.001_wave_right_dt0.002.dat};
        \coordinate (bot) at (rel axis cs:1,0);

    \end{groupplot}
    \path (top)--(bot) coordinate[midway] (group center);
    \node[right=1em,inner sep=0pt] at(group center -| current bounding box.east) {\pgfplotslegendfromname{bcs}};

\end{tikzpicture}
    \caption[]{Downstream propagating linear acoustic pulse with strength $\alpha = 0.001$: incoming $\overline{w}^-$ and outgoing $\overline{w}^+$ acoustic characteristics at the boundary $\Gamma_{out}$.}
    \label{fig::w_alpha0.001_right}
\end{figure}
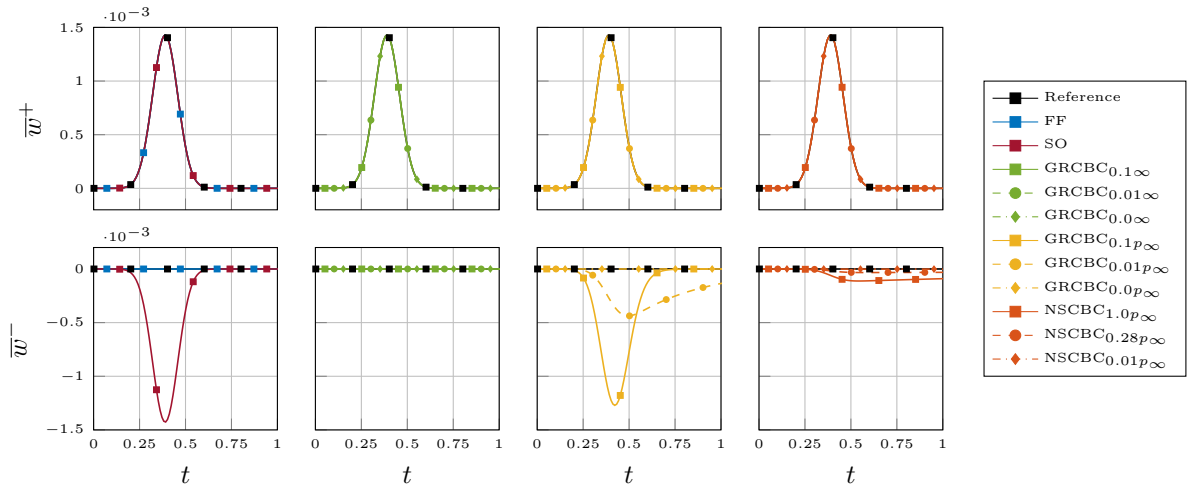

\begin{figure}[htb]
    \centering
    \input{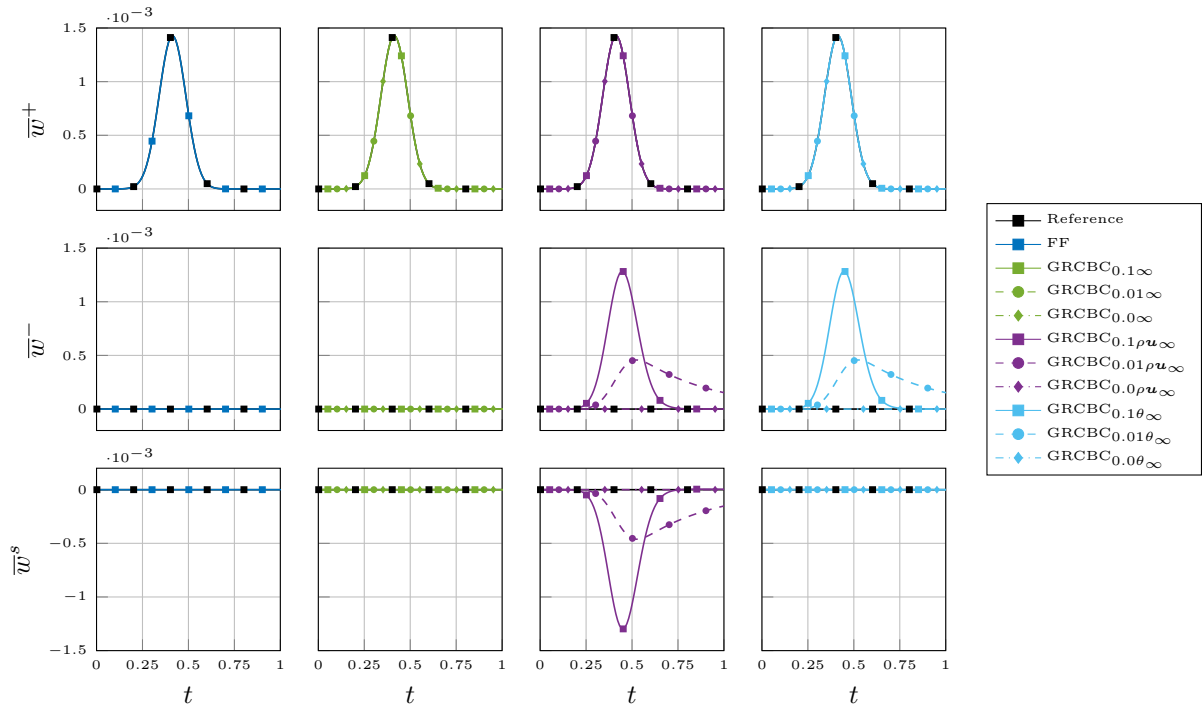}
    \caption[]{Upstream propagating linear acoustic pulse with strength $\alpha = 0.001$: incoming $\overline{w}^-$ and outgoing $\overline{w}^+$ acoustic characteristics and incoming entropy wave $\overline{w}^s$ at the boundary $\Gamma_{in}$.}
    \label{fig::w_alpha0.001_left}
\end{figure}

For a linear acoustic pulse with strength $\alpha = 0.001$, the mean incoming $\overline{w}^-$
and outgoing $\overline{w}^+$ acoustic characteristics
for the downstream propagating pulse are shown in Figure~\ref{fig::w_alpha0.001_right},
while the upstream propagation pulse are shown in Figure~\ref{fig::w_alpha0.001_left}.
The scaled errors $e_p$ for the downstream and upstream propagating pulse are shown in Figure~\ref{fig::p1d_0.001_right} and Figure~\ref{fig::p1d_0.001_left}, respectively.
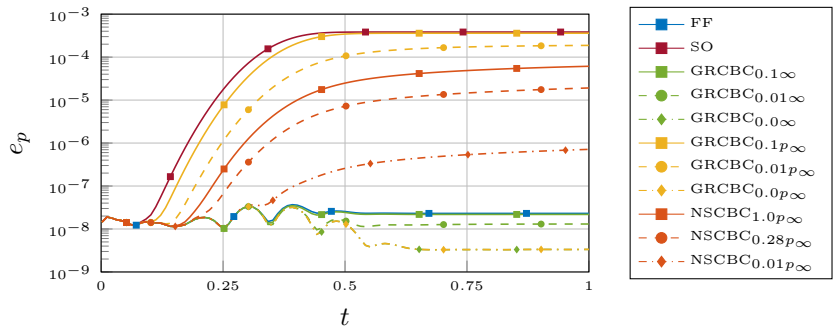
\begin{figure}[htb]
    \centering
    \tikzsetnextfilename{planar_acoustic_pulse_alpha0.001_right}
\begin{tikzpicture}
    \pgfplotsset{every axis/.append style={
                width=8cm,
                height=5cm,
                mark repeat = 20,
                mark size = 1pt,
                mark options={solid},
                xmin=0,
                xmax=1.0,
                max space between ticks=20,
                semithick,
                xtick pos=bottom,
                xmajorgrids=true, ymajorgrids=true,
                axis line style={thin},
                tick label style = {font=\tiny},
                ytick pos=left,
                xtick = {0,0.25,0.5,0.75,1.0}},
        legend style={font=\tiny},
        legend cell align={left},
        legend columns=1,
    }

    \begin{groupplot}[group style={group size=1 by 1, x descriptions at=edge bottom, vertical sep=0.5cm, horizontal sep=1cm}]

        \nextgroupplot[ymode=log, ymin=1e-9, ymax=1e-3, ylabel = $e_p$, xlabel = $t$, legend to name=pulse_right]

        \addplot[matblue,  mark phase=8, mark=square*, mark options={solid}]  table [y=e_p, x=t]{farfield_inflow_and_outflow_Ma0.03_alpha0.001_wave_right_dt0.002.dat};         \addlegendentry{FF};
        \addplot[matred, mark phase=15, mark=square*, mark options={solid}] table [y=e_p, x=t]{farfield_inflow_and_pressure_outflow_Ma0.03_alpha0.001_wave_right_dt0.002.dat};  \addlegendentry{SO};

        \addplot[matgreen, mark phase=6, mark=square*, mark options={solid}]           table [y=e_p_0.1, x=t]{grcbc_farfield_inflow_and_outflow_Ma0.03_alpha0.001_wave_right_dt0.002.dat};        \addlegendentry{GRCBC$_{0.1\infty}$};  
        \addplot[matgreen, mark phase=11, mark=*, dashed, mark options={solid}]   table [y=e_p_0.01, x=t]{grcbc_farfield_inflow_and_outflow_Ma0.03_alpha0.001_wave_right_dt0.002.dat};            \addlegendentry{GRCBC$_{0.01 \infty}$};
        \addplot[matgreen, mark phase=6, mark=diamond*, dash dot, mark options={solid}] table [y=e_p_0.0, x=t]{grcbc_farfield_inflow_and_outflow_Ma0.03_alpha0.001_wave_right_dt0.002.dat};      \addlegendentry{GRCBC$_{0.0\infty}$};

        \addplot[matyellow, mark phase=6, mark=square*, mark options={solid}]           table [y=e_p_0.1, x=t]{grcbc_farfield_inflow_and_pressure_outflow_Ma0.03_alpha0.001_wave_right_dt0.002.dat};   \addlegendentry{GRCBC$_{0.1p_{\infty}}$};
        \addplot[matyellow, mark phase=11,mark=*, dashed, mark options={solid}]   table [y=e_p_0.01, x=t]{grcbc_farfield_inflow_and_pressure_outflow_Ma0.03_alpha0.001_wave_right_dt0.002.dat};   \addlegendentry{GRCBC$_{0.01p_{\infty}}$};
        \addplot[matyellow, mark phase=11,mark=diamond*, dash dot, mark options={solid}] table [y=e_p_0.0, x=t]{grcbc_farfield_inflow_and_pressure_outflow_Ma0.03_alpha0.001_wave_right_dt0.002.dat};   \addlegendentry{GRCBC$_{0.0p_{\infty}}$};

        \addplot[matorange, mark phase=6, mark=square*, mark options={solid}]           table [y=e_p_1.0, x=t]{nscbc_farfield_inflow_and_pressure_outflow_Ma0.03_alpha0.001_wave_right_dt0.002.dat};    \addlegendentry{NSCBC$_{1.0p_{\infty}}$};
        \addplot[matorange, mark phase=11,mark=*, dashed, mark options={solid}]   table [y=e_p_0.28, x=t]{nscbc_farfield_inflow_and_pressure_outflow_Ma0.03_alpha0.001_wave_right_dt0.002.dat};     \addlegendentry{NSCBC$_{0.28p_{\infty}}$};
        \addplot[matorange, mark phase=16,mark=diamond*, dash dot, mark options={solid}] table [y=e_p_0.01, x=t]{nscbc_farfield_inflow_and_pressure_outflow_Ma0.03_alpha0.001_wave_right_dt0.002.dat};    \addlegendentry{NSCBC$_{0.01p_{\infty}}$};

        \coordinate (top) at (rel axis cs:1,1);
        \coordinate (bottom) at (rel axis cs:1,0);

    \end{groupplot}
    \path (top)--(bottom) coordinate[midway] (group center);
    \node[right=1em,inner sep=0pt] at(group center -| current bounding box.east) {\pgfplotslegendfromname{pulse_right}};

\end{tikzpicture}
    \caption[]{Downstream propagating linear acoustic pulse with strength $\alpha = 0.001$: scaled $L_2(\Omega)$-errors of the pressure $e_p$.}
    \label{fig::p1d_0.001_right}
\end{figure}
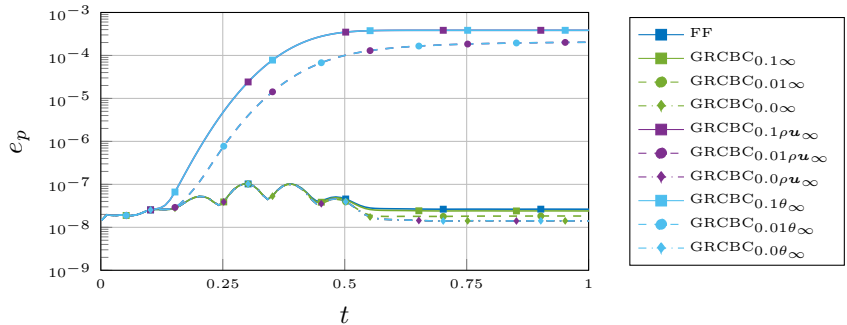
\begin{figure}[htb]
    \centering
    \tikzsetnextfilename{planar_acoustic_pulse_alpha0.001_left}
\begin{tikzpicture}
    \pgfplotsset{every axis/.append style={
                width=8cm,
                height=5cm,
                mark repeat = 20,
                mark size = 1pt,
                mark options={solid},
                xmin=0,
                xmax=1.0,
                max space between ticks=20,
                semithick,
                xtick pos=bottom,
                xmajorgrids=true, ymajorgrids=true,
                axis line style={thin},
                tick label style = {font=\tiny},
                ytick pos=left,
                xtick = {0,0.25,0.5,0.75,1.0}},
        legend style={font=\tiny},
        legend cell align={left},
        legend columns=1,
    }

    \begin{groupplot}[group style={group size=1 by 1, x descriptions at=edge bottom, vertical sep=0.5cm, horizontal sep=1cm}]
        \nextgroupplot[ymode=log, ymin=1e-9, ymax=1e-3,  ylabel = $e_p$, xlabel = $t$,  legend to name=pulse_left]
        
        \addplot[matblue, mark phase=11, mark=square*]  table [y=e_p, x=t]{farfield_inflow_and_outflow_Ma0.03_alpha0.001_wave_left_dt0.002.dat};  \addlegendentry{FF};
        
        \addplot[matgreen, mark phase=6, mark=square*]           table [y=e_p_0.1, x=t]{grcbc_farfield_inflow_and_outflow_Ma0.03_alpha0.001_wave_left_dt0.002.dat};    \addlegendentry{GRCBC$_{0.1\infty}$};  
        \addplot[matgreen, mark phase=11,mark=*, dashed]   table [y=e_p_0.01, x=t]{grcbc_farfield_inflow_and_outflow_Ma0.03_alpha0.001_wave_left_dt0.002.dat};  \addlegendentry{GRCBC$_{0.01 \infty}$};
        \addplot[matgreen, mark phase=16,mark=diamond*, dash dot] table [y=e_p_0.0, x=t]{grcbc_farfield_inflow_and_outflow_Ma0.03_alpha0.001_wave_left_dt0.002.dat};  \addlegendentry{GRCBC$_{0.0\infty}$};
        
        \addplot[matpurple, mark phase=11, mark=square*]           table [y=e_p_0.1, x=t]{grcbc_mass_inflow_Ma0.03_alpha0.001_wave_left_dt0.002.dat};    \addlegendentry{GRCBC$_{0.1\rho \VEL_\infty}$};
        \addplot[matpurple, mark phase=16,mark=*, dashed]   table [y=e_p_0.01, x=t]{grcbc_mass_inflow_Ma0.03_alpha0.001_wave_left_dt0.002.dat}; \addlegendentry{GRCBC$_{0.01\rho \VEL_\infty}$};
        \addplot[matpurple, mark phase=6,mark=diamond*, dash dot] table [y=e_p_0.0, x=t]{grcbc_mass_inflow_Ma0.03_alpha0.001_wave_left_dt0.002.dat}; \addlegendentry{GRCBC$_{0.0\rho \VEL_\infty}$};

        \addplot[matazure, mark phase=16, mark=square*]           table [y=e_p_0.1, x=t]{grcbc_temperature_inflow_Ma0.03_alpha0.001_wave_left_dt0.002.dat};  \addlegendentry{GRCBC$_{0.1 \theta_\infty}$};
        \addplot[matazure, mark phase=6,mark=*, dashed]   table [y=e_p_0.01, x=t]{grcbc_temperature_inflow_Ma0.03_alpha0.001_wave_left_dt0.002.dat};\addlegendentry{GRCBC$_{0.01 \theta_\infty}$};
        \addplot[matazure, mark phase=11,mark=diamond*, dash dot] table [y=e_p_0.0, x=t]{grcbc_temperature_inflow_Ma0.03_alpha0.001_wave_left_dt0.002.dat};\addlegendentry{GRCBC$_{0.0 \theta_\infty}$};

        \coordinate (top) at (rel axis cs:1,1);
        \coordinate (bottom) at (rel axis cs:1,0);

    \end{groupplot}
    \path (top)--(bottom) coordinate[midway] (group center);
    \node[right=1em,inner sep=0pt] at(group center -| current bounding box.east) {\pgfplotslegendfromname{pulse_left}};

\end{tikzpicture}
    \caption[]{Upstream propagating linear acoustic pulse with strength $\alpha = 0.001$: scaled $L_2(\Omega)$-errors of the pressure $e_p$.}
    \label{fig::p1d_0.001_left}
\end{figure}

From the figures it is evident that the far-field boundary condition (\ref*{leg::w-::ff} FF) and the characteristic boundary conditions with target state $\vec{U}_\infty$ (\ref*{leg::w-::grcbc_ff} GRCBC$_\infty$) closely match the reference solution for both the downstream and upstream propagating acoustic pulse.
In contrast, simulations enforcing the static pressure $\vec{U}_{p_\infty}$ at the outflow $\Gamma_{out}$ (\ref*{leg::w-::so}~SO, \ref*{leg::w-::grcbc_so}~GRCBC$_{p_\infty}$, \ref*{leg::w-::nscbc_so}~NSCBC$_{p_\infty}$) and either
the mass flux $\vec{U}_{\rho \VEL}$ or the static temperature $\vec{U}_{\theta_\infty}$ at the inflow $\Gamma_{in}$ (\ref*{leg::w-::grcbc_ma}~GRCBC$_{\rho u_\infty}$,  \ref*{leg::w-::grcbc_te}~GRCBC$_{\theta_\infty}$), result in total/partial reflection of the acoustic pulse.
This behavior is expected, as the imposition of hard boundaries inevitably leads to a reflection of the acoustic pulse.
Interestingly, the reflection at the inflow $\Gamma_{in}$ for the mass flux target state $\vec{U}_{\rho \VEL}$ leads to an incoming entropy wave $\overline{w}^s$,
which is absent when imposing the static temperature $\vec{U}_{\theta_\infty}$, due to the underlying isentropic relation \eqref{eq::isentropic_density}.

A way to overcome the reflective issue within the context of a one-dimensional perturbation
is by tuning the amount of relaxation $\hat{\mat{D}}$ towards the target state $\vec{U}^-$, as discussed in Subsection~\ref{sec::subsonic_outflow}.
For every characteristic boundary condition (GRCBC, NSCBC) under evaluation, we observe a significant reduction in the reflection for
vanishing relaxation matrix $\hat{\mat{D}} \longrightarrow \mat{0}$, to the point of perfect non-reflecting behavior.
Nevertheless, the results of the linear acoustic pulse suggest that for a purely planar acoustic perturbation, 
the far-field target state, $\vec{U}_\infty$, appears to be the most appropriate choice, without requiring precise tuning of the relaxation matrix~$\hat{\mat{D}}$.

\subsection{Oblique pressure pulse}
\label{sec::experiments::oblique_pressure_pulse}
\graphicspath{{./figures/oblique_pressure_pulse/data}}
\pgfplotsset{table/search path={./figures/oblique_pressure_pulse/data}}
This benchmark involves analyzing the propagation of a cylindrical non-linear pressure pulse within a squared domain $\Omega=(-5R, 5R)\times(-5R, 5R)$
as shown in Figure~\ref{fig::pulse_domain}, where $R=0.2$ is the radius of the initial perturbation. The objective is to examine the
non-reflective performance of the newly implemented characteristic boundary conditions when oblique waves impinge on the boundaries,
as well as the alternating inflow and outflow behavior at the boundaries triggered by the pressure pulse. The initial conditions are given as
\begin{align*}
    \rho_0  = \rho_\infty,                                             \qquad
    \VEL_0  = \vec{0},                                                 \qquad
    p_0     = p_\infty \left(1 + \alpha \exp\left(-r^2/R^2\right) \right),
\end{align*}
where $\alpha = 0.5$ is the pulse strength and $r^2=x^2 + y^2$ the polar coordinate \cite{shehadiPolynomialcorrectionNavierStokesCharacteristic2024a}. Initially, the fluid is at rest,
which translates to $\Ma_\infty = 0$, necessitating the adoption of acoustic scaling \eqref{eq::acoustic-scaling} of the Euler equations with dimensionless quantities $\rho_\infty=1$, $c_\infty = 1$ and $p_\infty = 1/\gamma$, since for an
aerodynamic scaling the speed of sounds tends to $c_\infty \rightarrow \infty$. Further, due to the absence of a baseflow, the convective stabilization matrix is 
chosen as a Lax-Friedrichs \cite{vila-perezHybridisableDiscontinuousGalerkin2021} $\mat{S}_c = \max(|\hat{u}_n| + \hat{c}, 0) \I$.
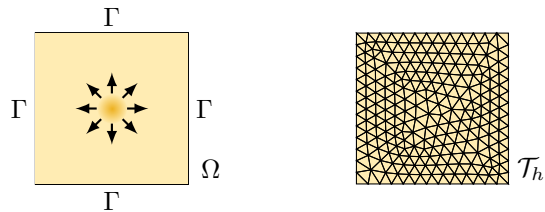
\begin{figure}[htb]
    \centering
    \tikzsetnextfilename{oblique_pressure_pulse_domain}
\begin{tikzpicture}[scale=2, >=latex]
    \draw[thick] (-0.5, -0.5) rectangle (0.5, 0.5);
    \fill[amber!30] (-0.5, -0.5) rectangle (0.5, 0.5);

    \node at (-0.61, 0) {$\Gamma$};
    \node at (0.61, 0) {$\Gamma$};
    \node at (0, 0.61) {$\Gamma$};
    \node at (0, -0.61) {$\Gamma$};

    \shade[outer color=amber!30, inner color=matyellow] (0, 0) circle (0.1);
    
    \node at (0.65, -0.4) {$\Omega$};

    \foreach \theta in {0, 45, 90, 135, 180, 225, 270, 315}
    {
        \draw[thick, ->] ({0.1*cos(\theta)}, {0.1*sin(\theta)}) -- ({0.24*cos(\theta)}, {0.24*sin(\theta)});;
    }

    \begin{scope}[scale=0.5, xshift = 120]
        \fill[amber!30] (-1, -1) rectangle (1, 1);
        \node at (1.3, -0.8) {$\mesh$};

        \pgfplotstableread[col sep=comma]{mesh.csv}\mesh
        \pgfplotstablegetrowsof{\mesh}
        \pgfmathsetmacro{\rows}{\pgfplotsretval-1}  

        \foreach \row in {0,...,\rows} {
            \pgfplotstablegetelem{\row}{x0}\of\mesh
            \let\xi=\pgfplotsretval
            \pgfplotstablegetelem{\row}{y0}\of\mesh
            \let\yi=\pgfplotsretval
            \pgfplotstablegetelem{\row}{x1}\of\mesh
            \let\xj=\pgfplotsretval
            \pgfplotstablegetelem{\row}{y1}\of\mesh
            \let\yj=\pgfplotsretval
            \pgfplotstablegetelem{\row}{x2}\of\mesh
            \let\xk=\pgfplotsretval
            \pgfplotstablegetelem{\row}{y2}\of\mesh
            \let\yk=\pgfplotsretval
            \draw (\xi,\yi) -- (\xj,\yj) -- (\xk,\yk) -- cycle;
        }

    \end{scope}

\end{tikzpicture}
    \caption[]{Representation of the domain $\Omega$ with corresponding boundaries $\Gamma$ and triangulation $\mesh$ with mesh size $h = 0.15$.}
    \label{fig::pulse_domain}
\end{figure}

\begin{figure}[htb]
    \centering
    \tikzsetnextfilename{oblique_pressure_pulse_alpha0.5}
\begin{tikzpicture}
    \pgfplotsset{every axis/.append style={
        width=8cm, 
        height=5cm, 
        mark repeat = 80,
        mark size = 1pt,
        mark options={solid},
        xmin=0, 
        xmax=4.0, 
         max space between ticks=20,
         semithick,
         xtick pos=bottom,
         xmajorgrids=true, ymajorgrids=true,
         axis line style={thin}, 
         tick label style = {font=\tiny},
         ytick pos=left,
        xtick = {0, 2, 3, 4}},        
        legend style={font=\tiny},
        legend cell align={left},
        legend columns=2,
        transpose legend,
        extra x ticks={0.752, 1.002, 1.252},
        extra x tick style={%
            grid=major,
        },
        extra x tick labels={
            $t_0$, $t_1$, $t_2$
        }
            }

        \begin{groupplot}[group style={group size=2 by 1, x descriptions at=edge bottom, vertical sep=0.5cm, horizontal sep=2cm}]

            \nextgroupplot[ymode=log, ymin=1e-6, ymax=1e1,ylabel = $e_p$,  xlabel = $t$, title=$(a)$]
            
            \addplot[matblue, mark=square*, mark options={solid}]  table [y=e_p, x=t]{standard_farfield_Ma0.0_alpha0.5.dat};      
            \addplot[matred, mark=square*, mark options={solid}] table [y=e_p, x=t]{standard_outflow_Ma0.0_alpha0.5.dat}; 
            
            \addplot[matgreen,mark phase = 20, mark=*, dashed, mark options={solid}]   table [y=e_p_0.01, x=t]{grcbc_farfield_Ma0.0_alpha0.5.dat};     

            \addplot[matyellow, mark phase = 40, mark=*, dashed, mark options={solid}]   table [y=e_p_0.01, x=t]{grcbc_outflow_Ma0.0_alpha0.5.dat};   

            \addplot[matorange, mark phase = 60, mark=*, dashed, mark options={solid}]   table [y=e_p_0.28, x=t]{nscbc_pressure_outflow_Ma0.0_alpha0.5.dat};     
            
            \addplot[matpurple, mark phase = 20, mark=*, dashed, mark options={solid}]   table [y=e_p_0.01, x=t]{grcbc_mass_Ma0.0_alpha0.5.dat};    
            
            \addplot[matazure, mark phase = 40, mark=*, dashed, mark options={solid}]   table [y=e_p_0.01, x=t]{grcbc_temperature_Ma0.0_alpha0.5.dat};    

            \coordinate (left) at (rel axis cs:0,0);
            
            
            \nextgroupplot[ymode=log, ylabel = $\|\Delta \tilde{p}\|$, xlabel = $t$, title=$(b)$, legend to name=pulse2d_alpha0.5, ymin=-0.005, ymax=0.065]S

            \addplot[black, mark=square*, mark options={solid}]  table [y=p, x=t]{grcbc_reference_Ma0.0_alpha0.5.dat};       \addlegendentry{Reference}; \label{leg::p2d::ref}
            
            \addplot[matblue, mark=square*, mark options={solid}]  table [y=p, x=t]{standard_farfield_Ma0.0_alpha0.5.dat};       \addlegendentry{FF}; \label{leg::p2d::ff}
            \addplot[matred, mark=square*, mark options={solid}] table [y=p, x=t]{standard_outflow_Ma0.0_alpha0.5.dat};  \addlegendentry{SO}; \label{leg::p2d::so}
        
            \addplot[matgreen, mark phase = 20, mark=*, dashed, mark options={solid}]   table [y=p_0.01, x=t]{grcbc_farfield_Ma0.0_alpha0.5.dat};      \addlegendentry{GRCBC$_{0.01 \infty}$}; \label{leg::p2d::grcbc_ff}

            \addplot[matyellow, mark phase = 40, mark=*, dashed, mark options={solid}]   table [y=p_0.01, x=t]{grcbc_outflow_Ma0.0_alpha0.5.dat};       \addlegendentry{GRCBC$_{0.01p_{\infty}}$}; \label{leg::p2d::grcbc_so}

            \addplot[matorange, mark phase = 60, mark=*, dashed, mark options={solid}]   table [y=p_0.28, x=t]{nscbc_pressure_outflow_Ma0.0_alpha0.5.dat};        \addlegendentry{NSCBC$_{0.28p_{\infty}}$}; \label{leg::p2d::nscbc_so}

            \addplot[matpurple, mark phase = 20, mark=*, dashed, mark options={solid}]   table [y=p_0.01, x=t]{grcbc_mass_Ma0.0_alpha0.5.dat};    \addlegendentry{GRCBC$_{0.01\rho \VEL_\infty}$}; \label{leg::p2d::grcbc_ma}

            \addplot[matazure, mark phase = 40, mark=*, dashed, mark options={solid}]   table [y=p_0.01, x=t]{grcbc_temperature_Ma0.0_alpha0.5.dat};    \addlegendentry{GRCBC$_{0.01 \theta_\infty}$}; \label{leg::p2d::grcbc_te}


            \coordinate (right) at (rel axis cs:1,0);

        \end{groupplot}
        \path (left)--(right) coordinate[midway] (group center);
        \node[inner sep=0pt,yshift=-4.5em] at(group center) {\pgfplotslegendfromname{pulse2d_alpha0.5}};

\end{tikzpicture}
    \caption[]{Pressure pulse with strength $\alpha = 0.5$: $(a)$ scaled $L_2(\Omega)$-errors of the pressure $e_p$ and $(b)$ $L_2(\Omega)$-norm of the scaled fluctuating pressure $\Delta \tilde{p}$.}
    \label{fig::p2d_0.5}
\end{figure}
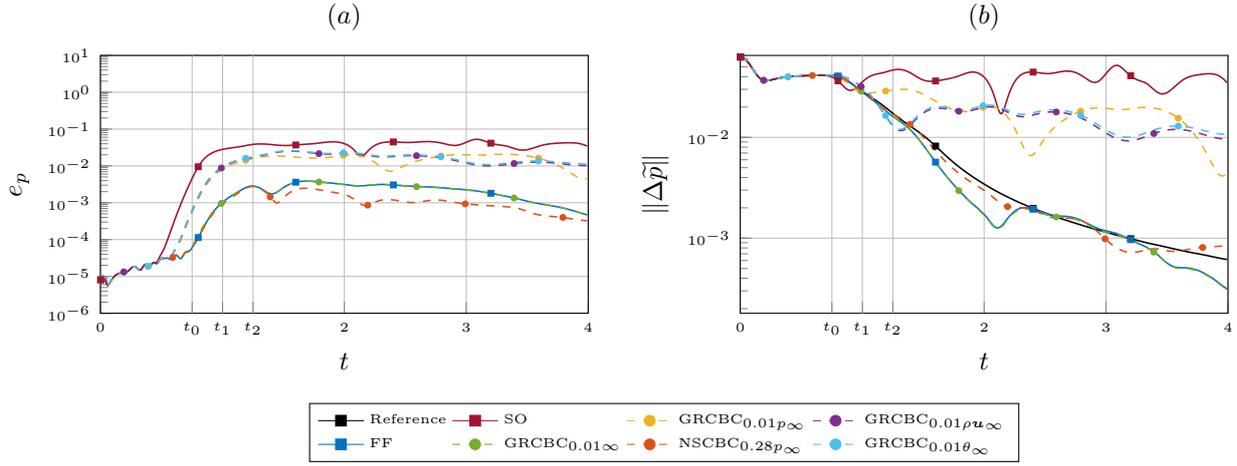

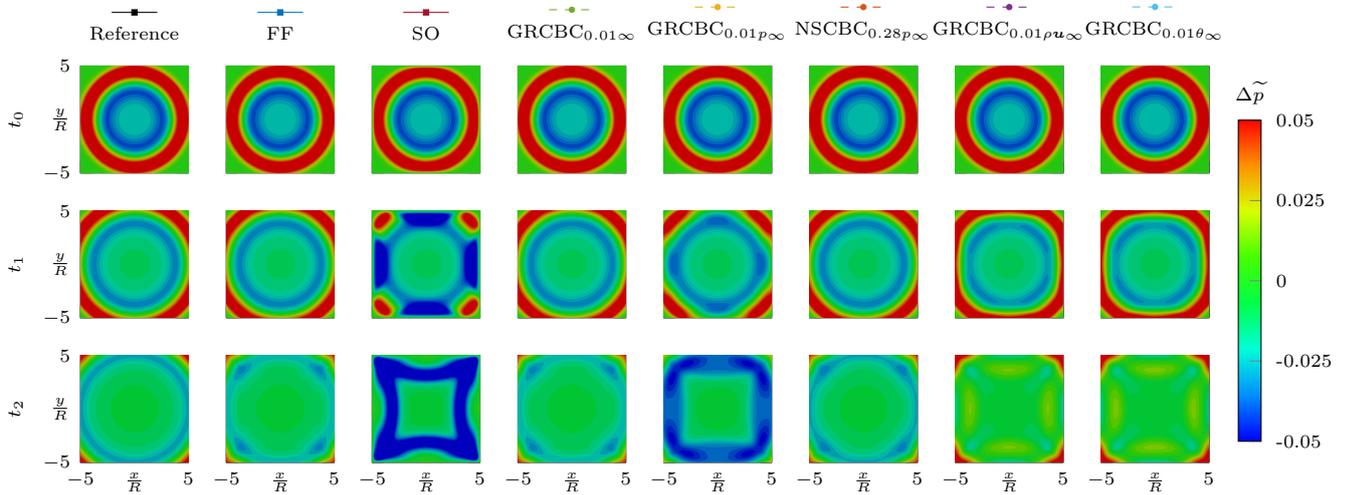
\begin{figure}[htb]
    \centering
\begin{tikzpicture}
    \pgfplotsset{every axis/.append style={
        width=3cm, 
        height=3cm, 
        mark repeat = 15,
        mark size = 1pt,
        mark options={solid},
        tick label style = {font=\tiny},
        xmin=-0.5, 
        xmax=0.5, 
        ymin=-0.5,
        ymax=0.5, 
        yticklabel=\empty, 
        xticklabel=\empty, 
        xtick pos=bottom,
        xmajorgrids=true, ymajorgrids=true,
        title style={font=\tiny, align=center},
        ylabel style={font=\tiny},
        ytick = {-0.5, 0, 0.5},
        xtick = {-0.5, 0, 0.5}, 
        yticklabels = {$-5$, $\frac{y}{R}$, $5$},
        xticklabels = {$-5$, $\frac{x}{R}$, $5$}, 
         },
        }

        \begin{groupplot}[group style={group size=8 by 3, x descriptions at=edge bottom, vertical sep=0.5cm, horizontal sep=0.5cm}, colorbar to name=pulse2d_colorbar,   
         colorbar style={title=$\Delta \widetilde{p}$, height=3*\pgfkeysvalueof{/pgfplots/parent axis height}, width=0.3cm, ytick = {-0.05, -0.025, 0 , 0.025, 0.05}, 
         yticklabels = {-0.05, -0.025, 0 , 0.025, 0.05}, ytick pos=right, scaled y ticks=false, yticklabel style={anchor=west, xshift=0.5em}}, point meta max = 0.05, point meta min = -0.05]
            
            \nextgroupplot[title=\ref*{leg::p2d::ref}\\Reference, ylabel=$t_0$, colorbar, xticklabel=\empty]
            \addplot graphics[xmin=-0.5,ymin=-0.5,xmax=0.5,ymax=0.5] {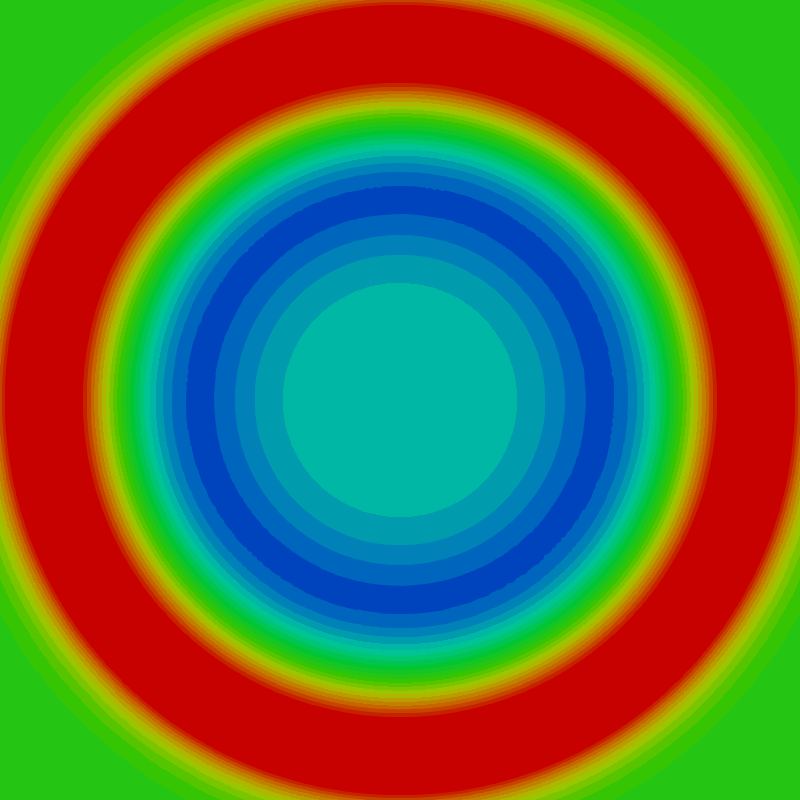};

            \nextgroupplot[title=\ref*{leg::p2d::ff}\\FF, yticklabel=\empty, xticklabel=\empty]
            \addplot graphics[xmin=-0.5,ymin=-0.5,xmax=0.5,ymax=0.5] {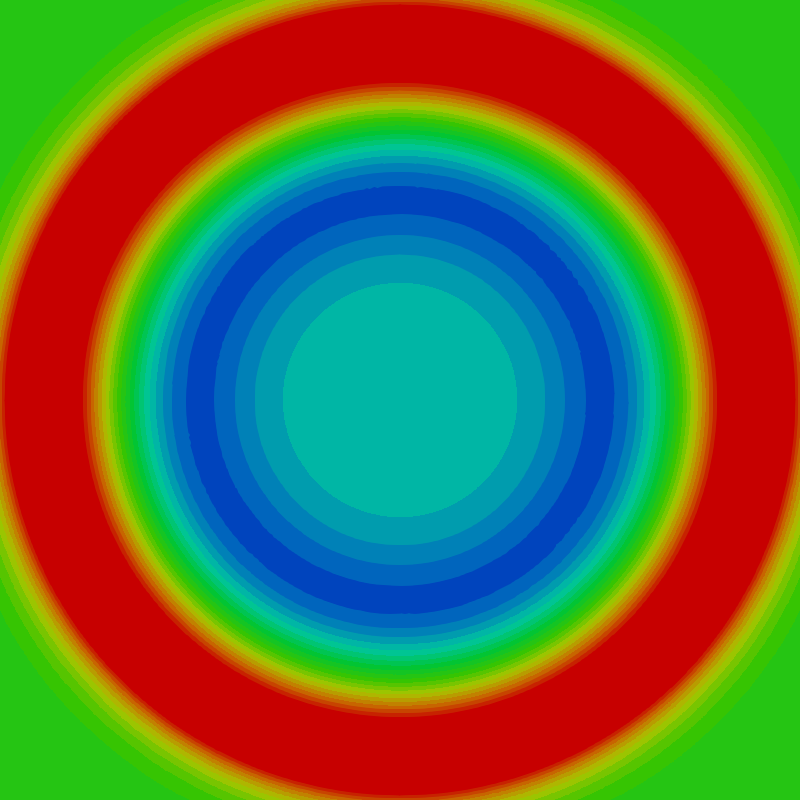};

            \nextgroupplot[title=\ref*{leg::p2d::so}\\SO, yticklabel=\empty, xticklabel=\empty]
            \addplot graphics[xmin=-0.5,ymin=-0.5,xmax=0.5,ymax=0.5] {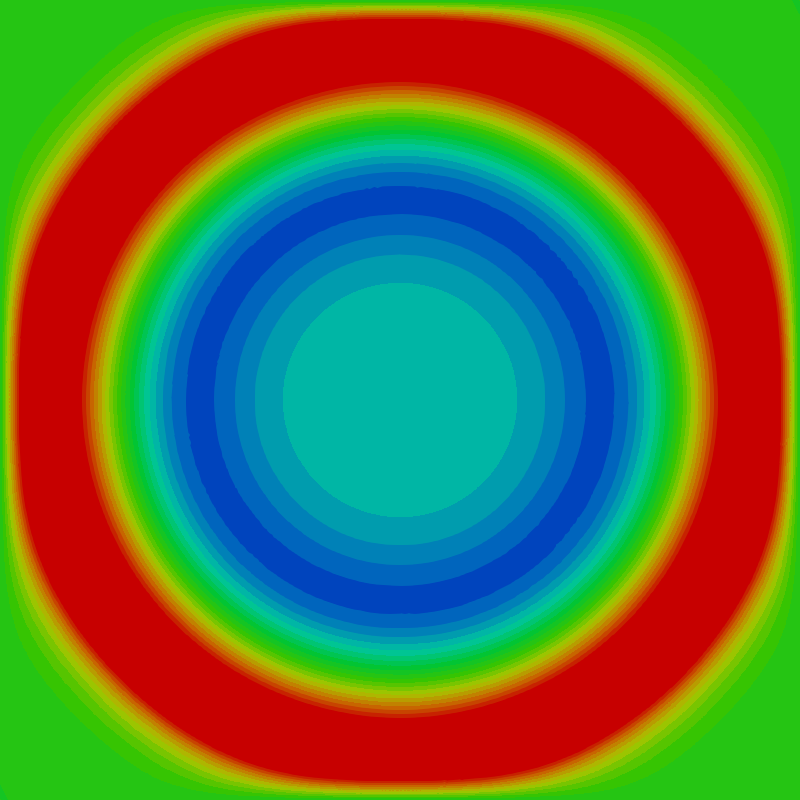};

            \nextgroupplot[title=\ref*{leg::p2d::grcbc_ff}\\GRCBC$_{0.01 \infty}$, yticklabel=\empty, xticklabel=\empty]
            \addplot graphics[xmin=-0.5,ymin=-0.5,xmax=0.5,ymax=0.5] {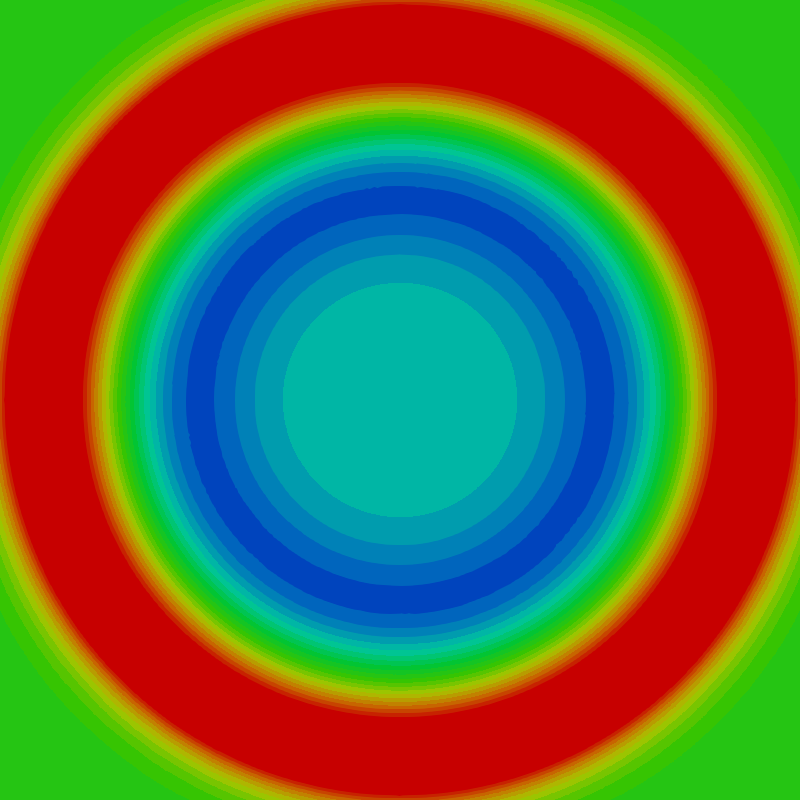};

            \nextgroupplot[title=\ref*{leg::p2d::grcbc_so}\\GRCBC$_{0.01p_{\infty}}$, yticklabel=\empty, xticklabel=\empty]
            \addplot graphics[xmin=-0.5,ymin=-0.5,xmax=0.5,ymax=0.5] {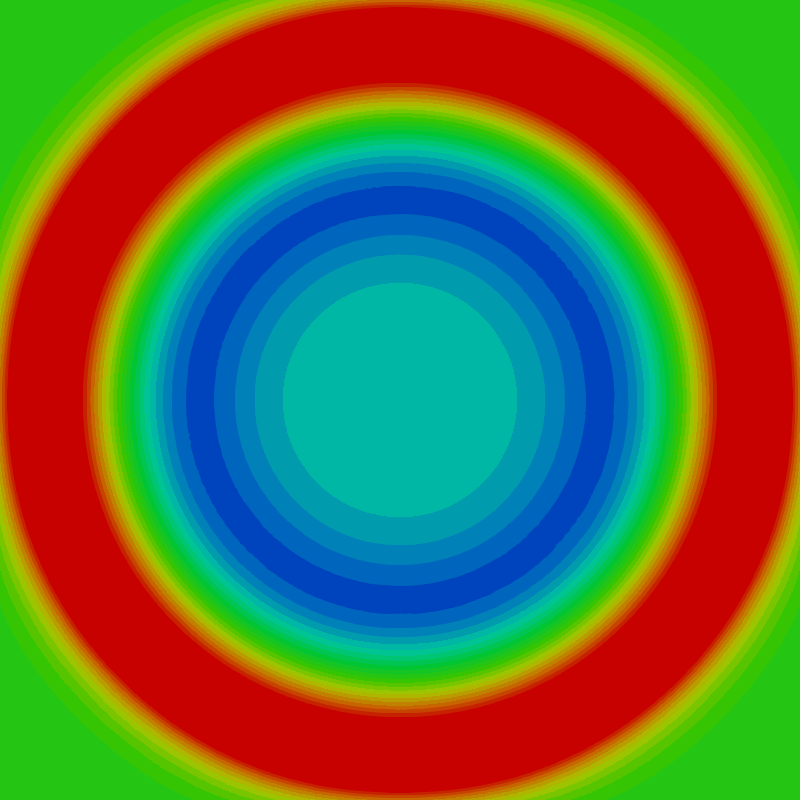};

            \nextgroupplot[title=\ref*{leg::p2d::nscbc_so}\\NSCBC$_{0.28p_{\infty}}$, yticklabel=\empty, xticklabel=\empty]
            \addplot graphics[xmin=-0.5,ymin=-0.5,xmax=0.5,ymax=0.5] {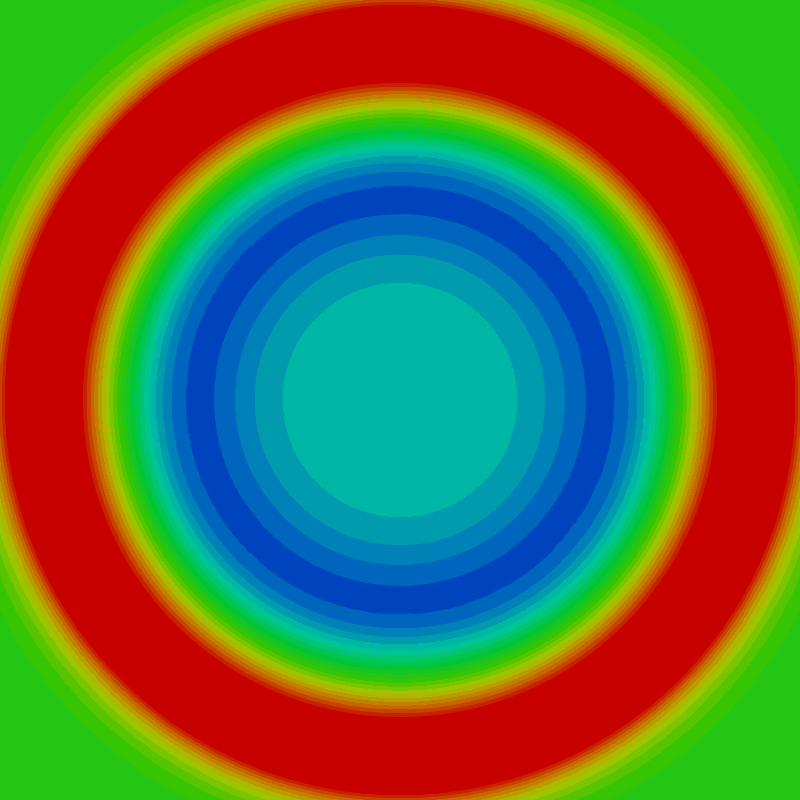};

            \nextgroupplot[title=\ref*{leg::p2d::grcbc_ma}\\GRCBC$_{0.01\rho \VEL_\infty}$, yticklabel=\empty, xticklabel=\empty]
            \addplot graphics[xmin=-0.5,ymin=-0.5,xmax=0.5,ymax=0.5] {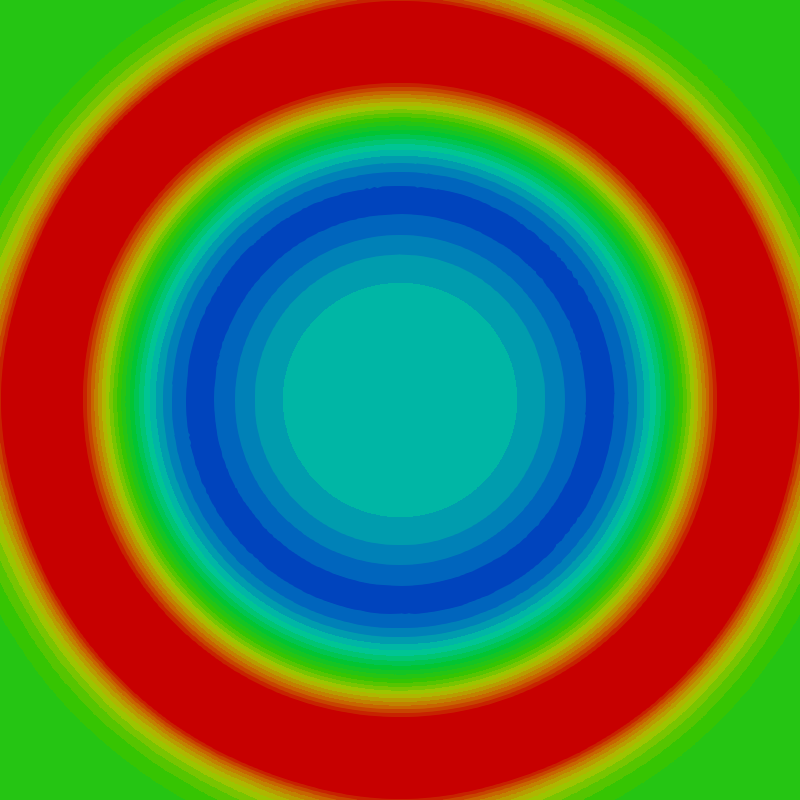};

            \nextgroupplot[title=\ref*{leg::p2d::grcbc_te}\\GRCBC$_{0.01 \theta_\infty}$, yticklabel=\empty, xticklabel=\empty]
            \addplot graphics[xmin=-0.5,ymin=-0.5,xmax=0.5,ymax=0.5] {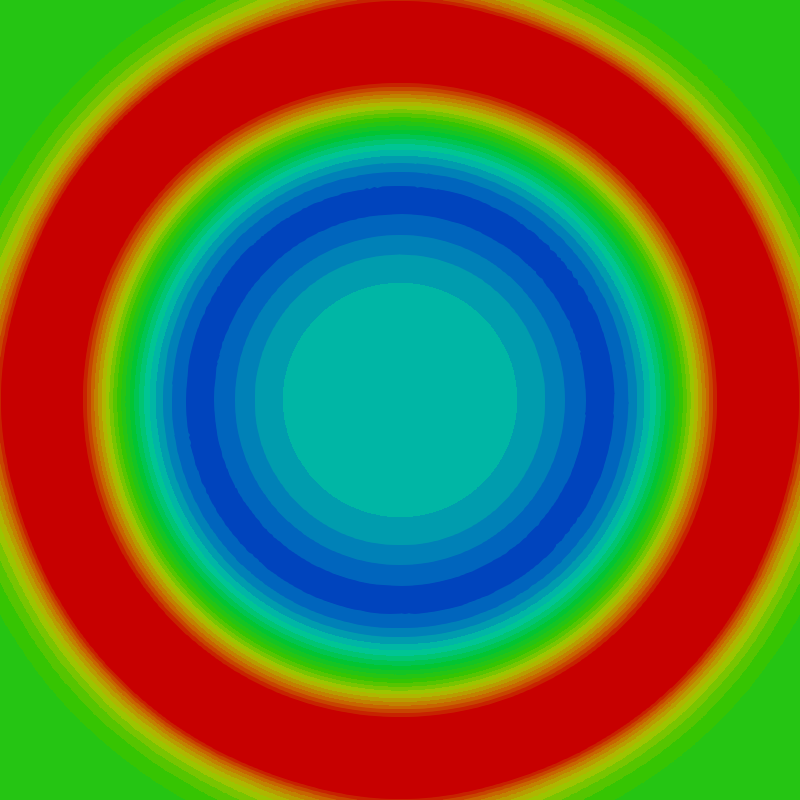};

            \coordinate (top) at (rel axis cs:1,1);

            \nextgroupplot[ylabel=$t_1$, xticklabel=\empty]
            \addplot graphics[xmin=-0.5,ymin=-0.5,xmax=0.5,ymax=0.5] {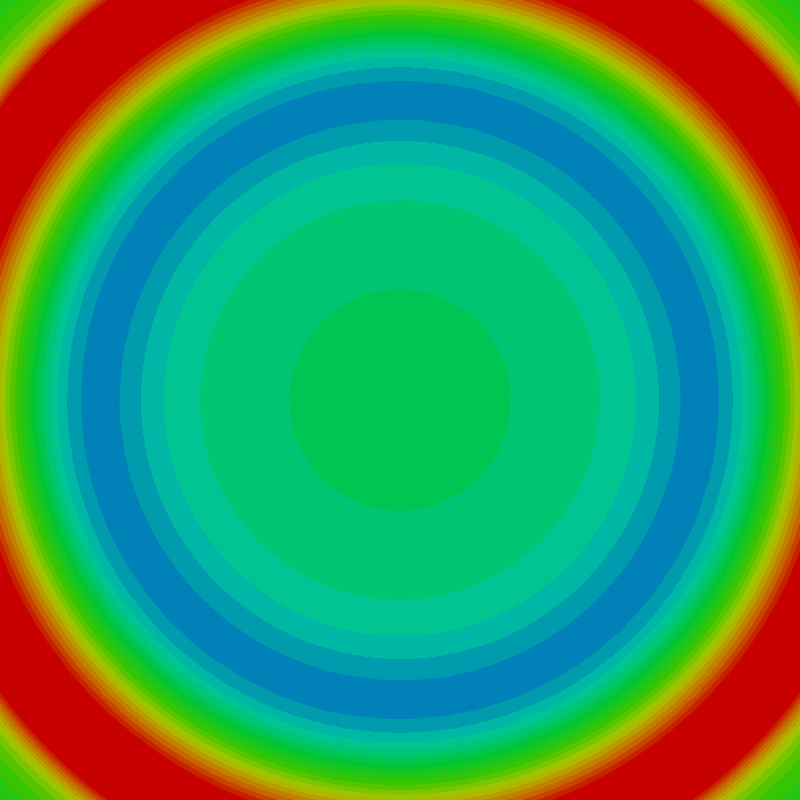};

            \nextgroupplot[yticklabel = \empty, xticklabel=\empty]
            \addplot graphics[xmin=-0.5,ymin=-0.5,xmax=0.5,ymax=0.5] {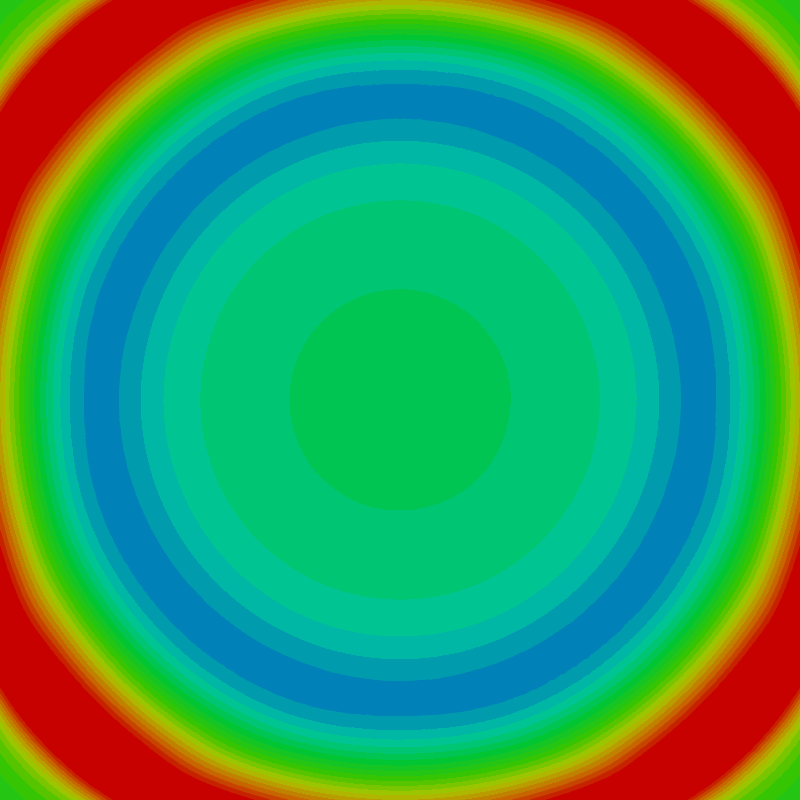};

            \nextgroupplot[yticklabel = \empty, xticklabel=\empty]
            \addplot graphics[xmin=-0.5,ymin=-0.5,xmax=0.5,ymax=0.5] {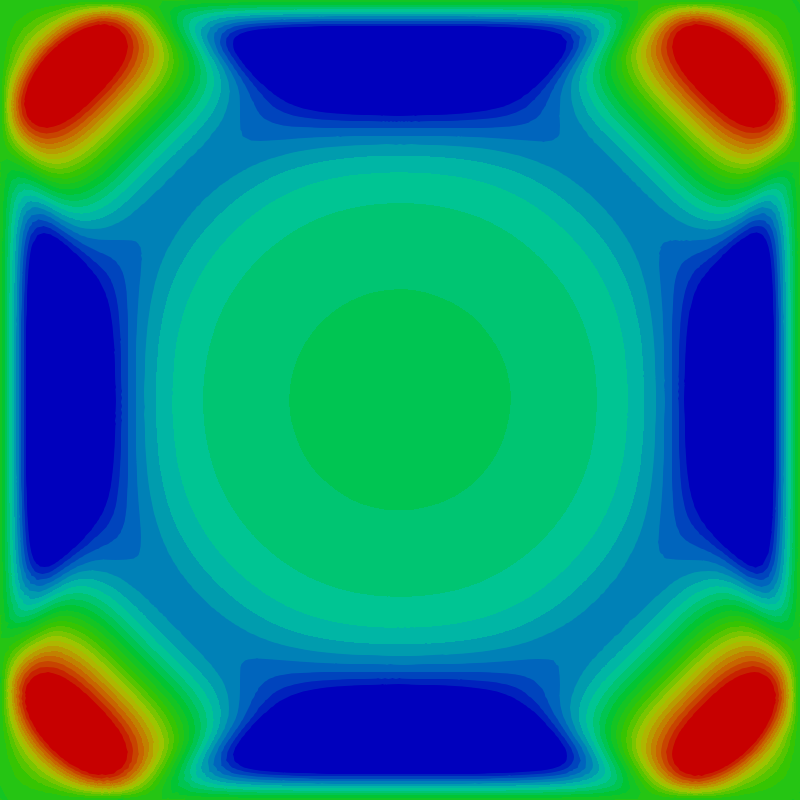};

            \nextgroupplot[yticklabel = \empty, xticklabel=\empty]
            \addplot graphics[xmin=-0.5,ymin=-0.5,xmax=0.5,ymax=0.5] {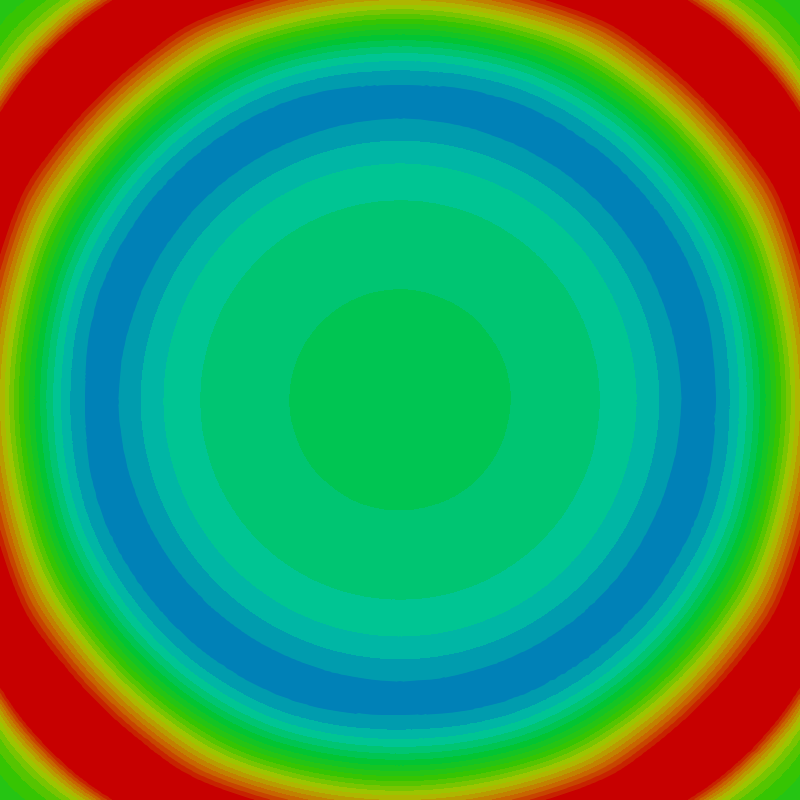};

            \nextgroupplot[yticklabel = \empty, xticklabel=\empty]
            \addplot graphics[xmin=-0.5,ymin=-0.5,xmax=0.5,ymax=0.5] {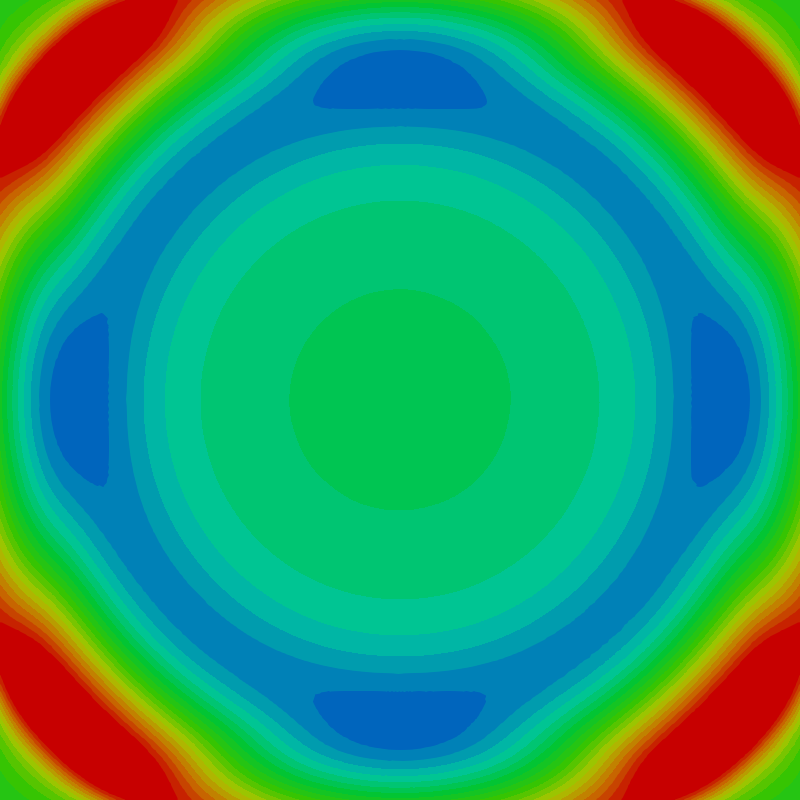};

            \nextgroupplot[yticklabel = \empty, xticklabel=\empty]
            \addplot graphics[xmin=-0.5,ymin=-0.5,xmax=0.5,ymax=0.5] {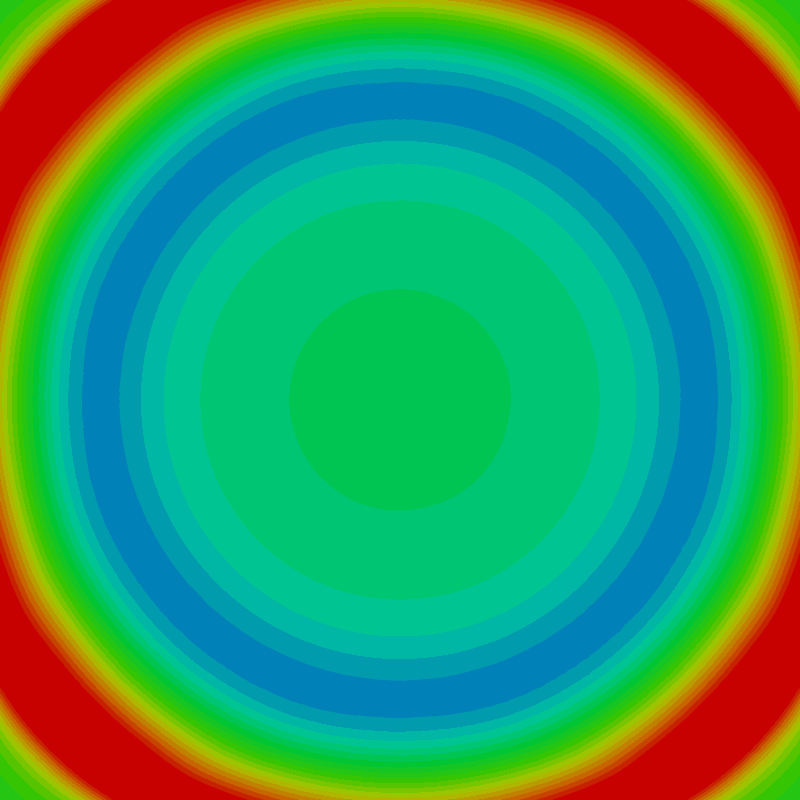};

            \nextgroupplot[yticklabel = \empty, xticklabel=\empty]
            \addplot graphics[xmin=-0.5,ymin=-0.5,xmax=0.5,ymax=0.5] {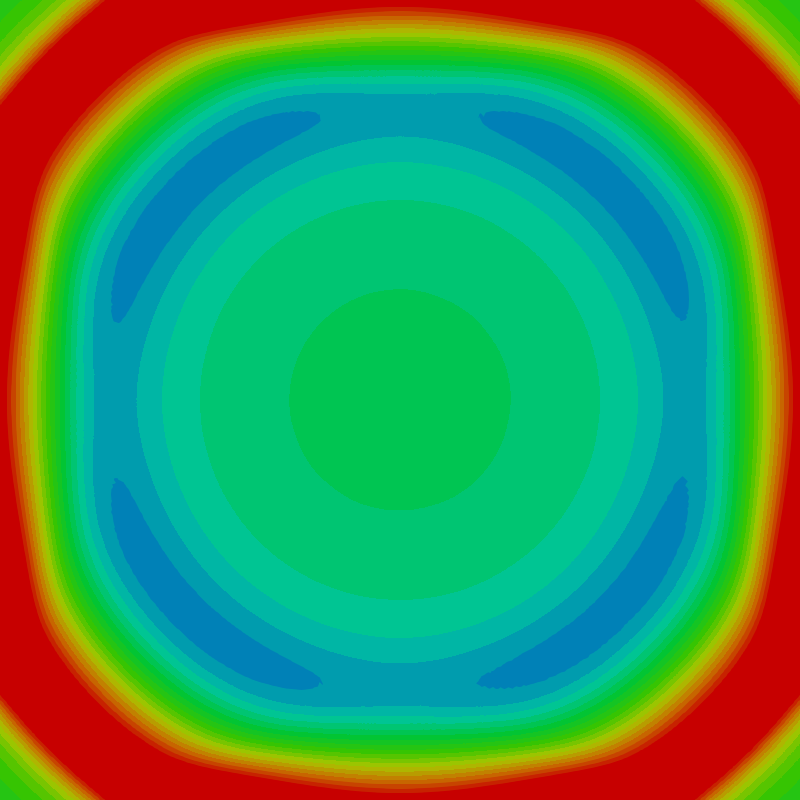};

            \nextgroupplot[yticklabel = \empty, xticklabel=\empty]
            \addplot graphics[xmin=-0.5,ymin=-0.5,xmax=0.5,ymax=0.5] {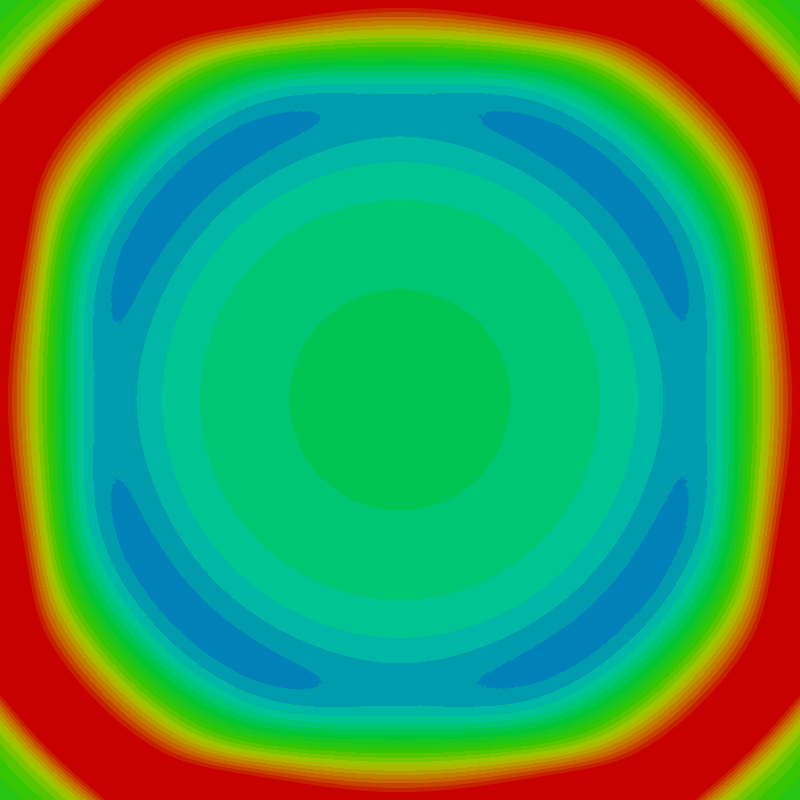};

            \nextgroupplot[ylabel=$t_2$]
            \addplot graphics[xmin=-0.5,ymin=-0.5,xmax=0.5,ymax=0.5] {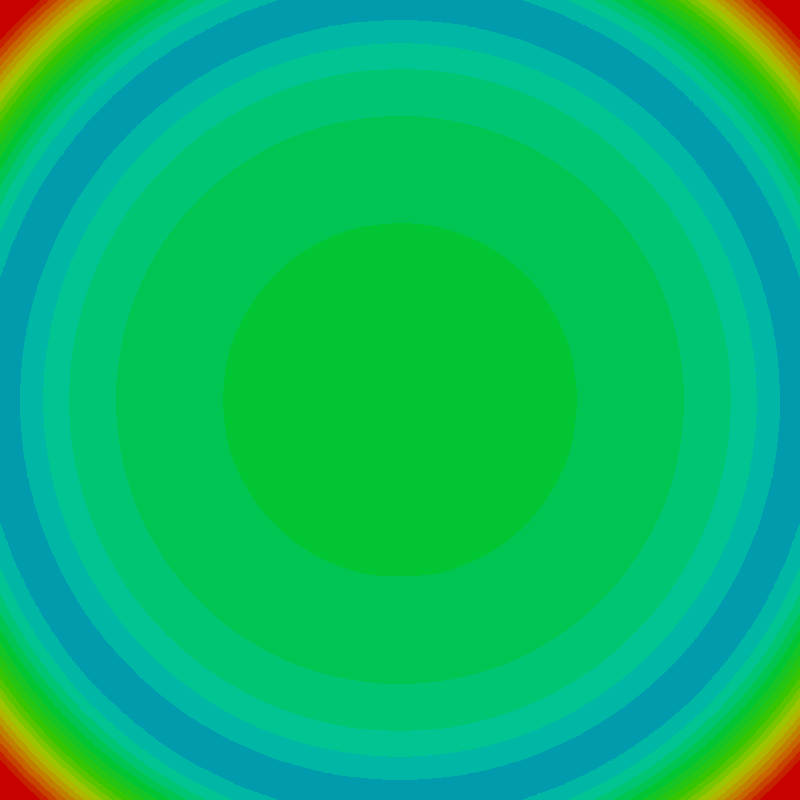};

            \nextgroupplot[yticklabel=\empty]
            \addplot graphics[xmin=-0.5,ymin=-0.5,xmax=0.5,ymax=0.5] {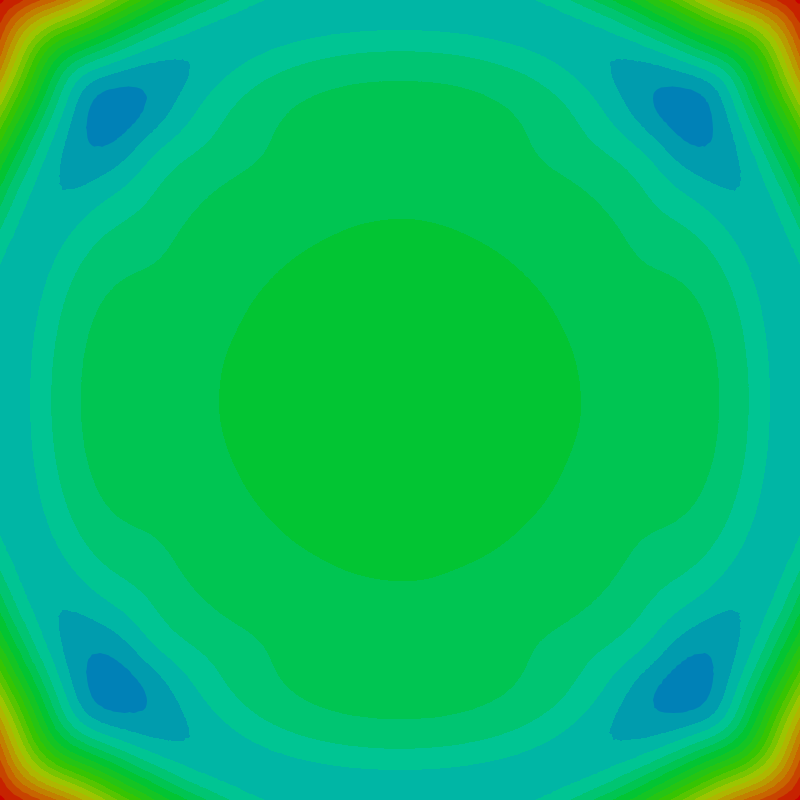};

            \nextgroupplot[yticklabel=\empty]
            \addplot graphics[xmin=-0.5,ymin=-0.5,xmax=0.5,ymax=0.5] {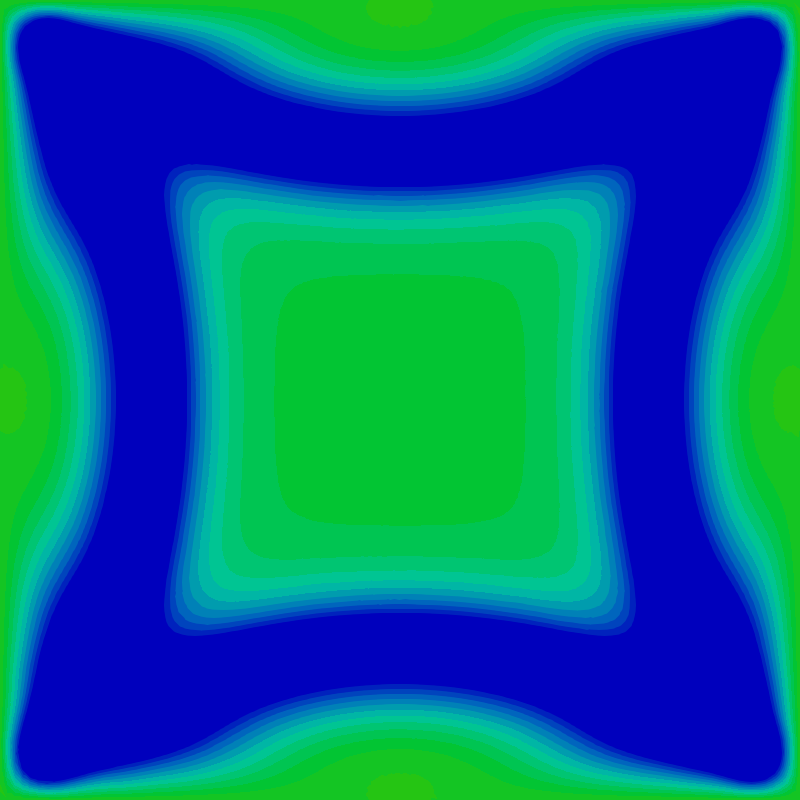};

            \nextgroupplot[yticklabel=\empty]
            \addplot graphics[xmin=-0.5,ymin=-0.5,xmax=0.5,ymax=0.5] {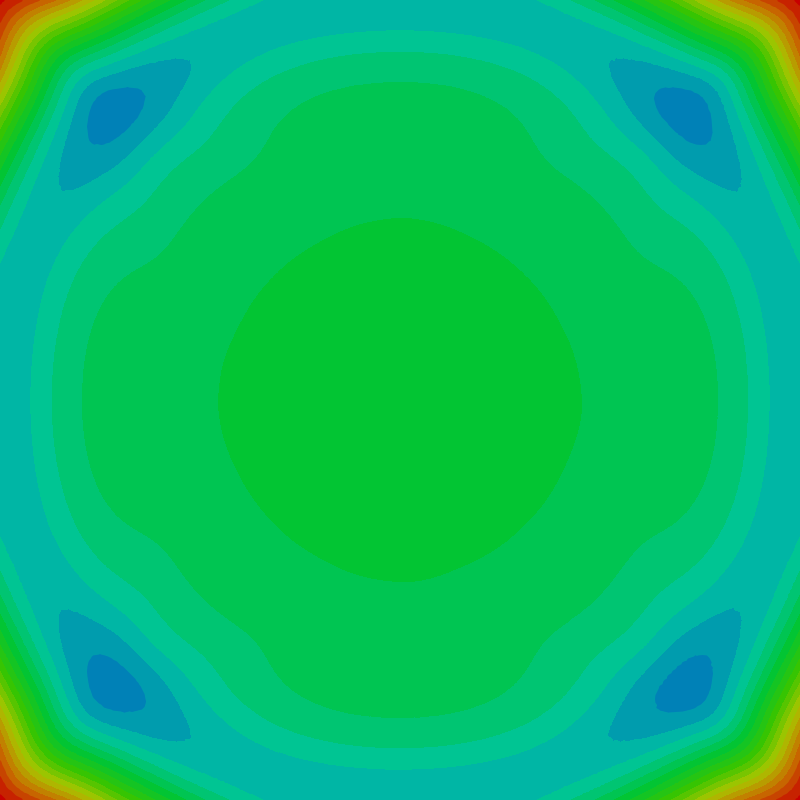};

            \nextgroupplot[yticklabel=\empty]
            \addplot graphics[xmin=-0.5,ymin=-0.5,xmax=0.5,ymax=0.5] {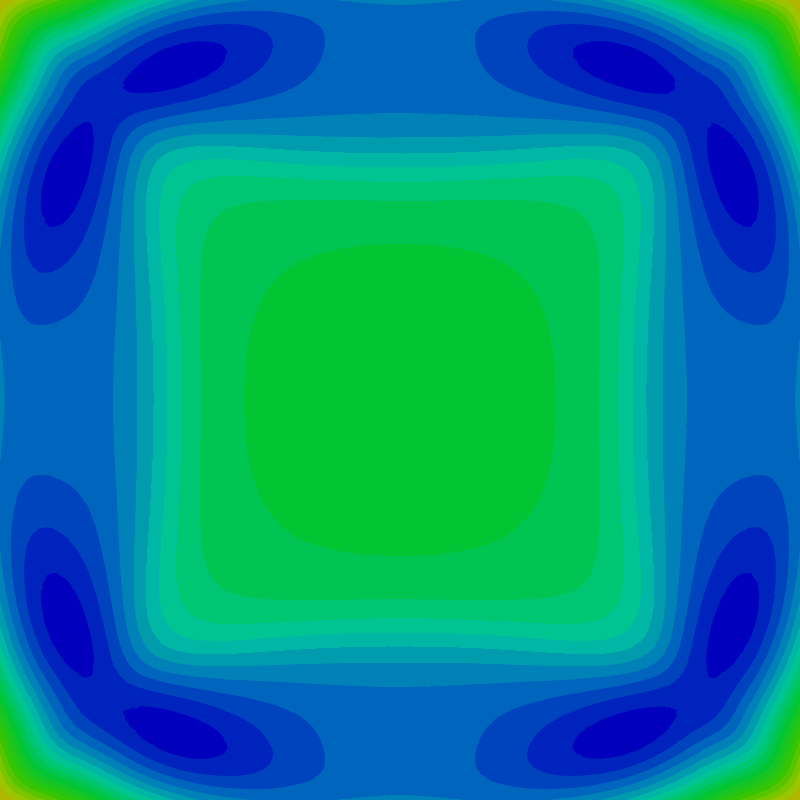};

            \nextgroupplot[yticklabel=\empty]
            \addplot graphics[xmin=-0.5,ymin=-0.5,xmax=0.5,ymax=0.5] {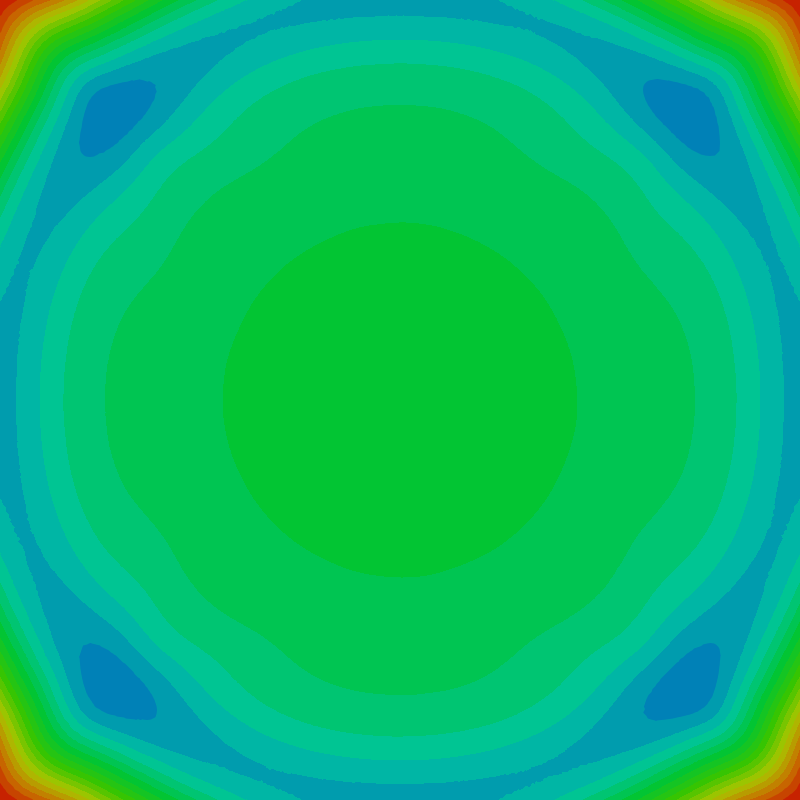};

            \nextgroupplot[yticklabel=\empty]
            \addplot graphics[xmin=-0.5,ymin=-0.5,xmax=0.5,ymax=0.5] {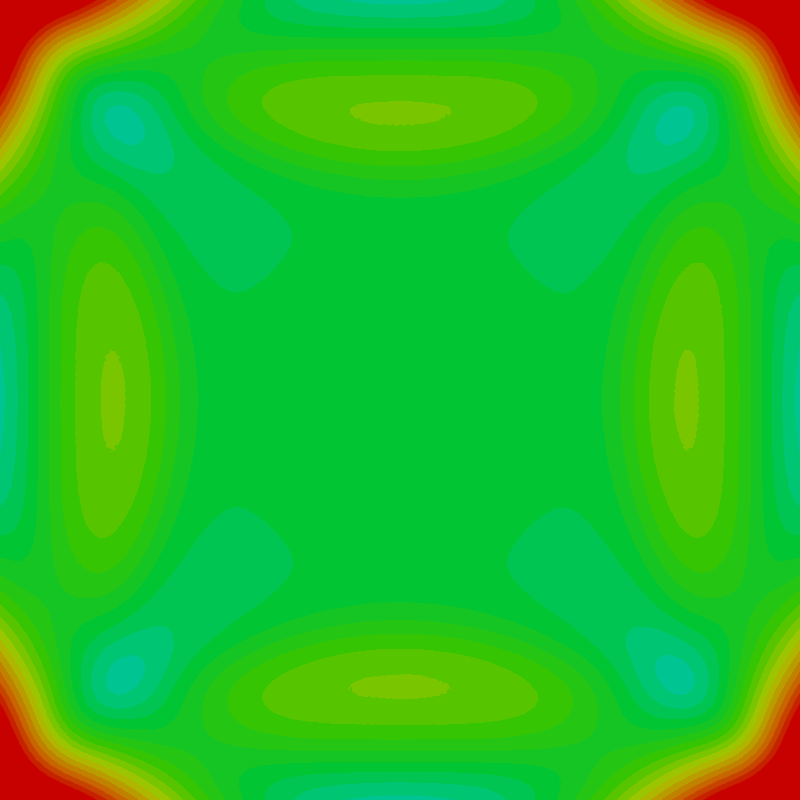};

            \nextgroupplot[yticklabel=\empty]
            \addplot graphics[xmin=-0.5,ymin=-0.5,xmax=0.5,ymax=0.5] {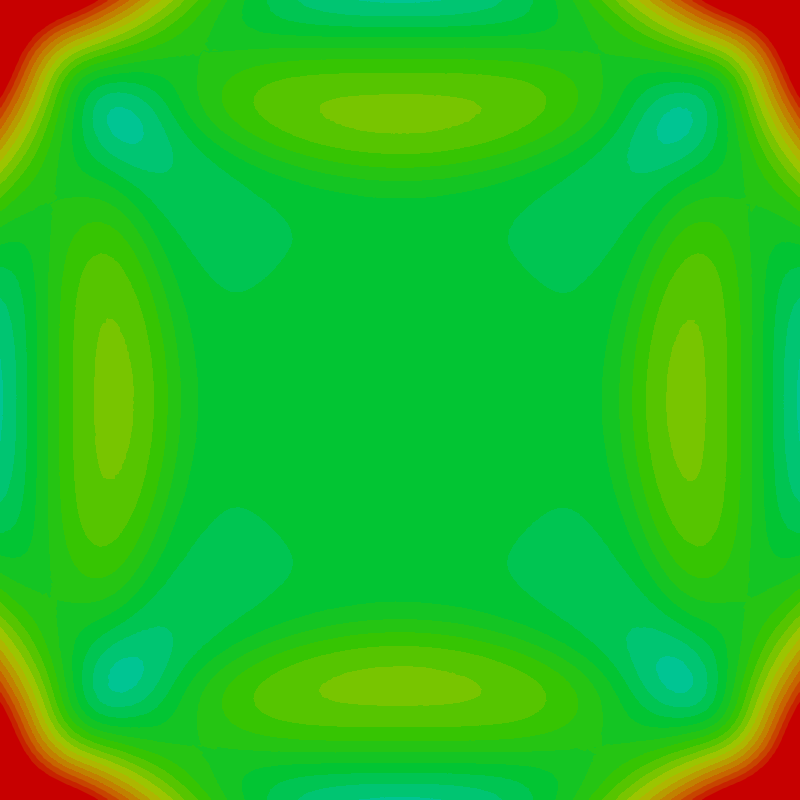};

            \coordinate (bottom) at (rel axis cs:1,0);

        \end{groupplot}

        \path (top)--(bottom) coordinate[midway] (group center);
        \node[right=0em,inner sep=0pt] at(group center -| current bounding box.east) {\pgfplotscolorbarfromname{pulse2d_colorbar}};

\end{tikzpicture}
    \caption[]{Pressure pulse with strength $\alpha = 0.5$: scaled fluctuating pressure $\Delta \tilde{p}$ for different boundary conditions at time instances $t_0 = 0.75$, $t_1 = 1$ and $t_2 = 1.25$.}
    \label{fig::p2d_0.5_snapshots}
\end{figure}

Similarly to the one-dimensional acoustic pulse, we compute a reference solution in a larger circular domain $\Omega_{r} = \left\lbrace (x, y): x^2 + y^2 \leq 40 R \right\rbrace$,
where the pressure pulse propagates undisturbed through the domain during the time interval $(0, t_{end}]$ without impinging on the boundaries.
The time-stepping is performed with a time step $\Delta t = 2\cdot 10^{-3}$ up to a dimensionless time $t_{end}=4$. Besides the scaled $L_2(\Omega)$-errors of the pressure~$e_p$,
in Figure~\ref{fig::p2d_0.5} we show the $L_2(\Omega)$-norm of the scaled fluctuating pressure
\begin{align}
    \label{eq::scaled_fluctuating_pressure}
    \Delta \tilde{p} & := \frac{p_h - p_\infty}{p_\infty},
\end{align}
to visualize the propagation of the reference solution. Note that we excluded the simulation results with no relaxation ($\hat{\mat{D}} = \mat{0}$), 
since they almost coincide with the far-field solution (\ref*{leg::p2d::ff} FF), in accordance to the observations drawn from Subsection~\ref{sec::experiments::planar_acoustic_pulse}. 
Similarly, also the cases with a large relaxation and target state $\vec{U}_{p_\infty}$ are omitted, as they produce the same 
outcome as the subsonic outflow boundary condition (\ref*{leg::p2d::so} SO). In Figure~\ref{fig::p2d_0.5_snapshots}, the
snapshots of the scaled fluctuating pressure $\Delta \tilde{p}$ are reported at time instances $t_0 = 0.75$, $t_1 = 1$ and $t_2 = 1.25$, during which the pressure pulse intersects the boundaries.


In Figure~\ref{fig::p2d_0.5}b we observe that at $t = t_{end}$ the reference solution (\ref*{leg::p2d::ref} Reference) has not yet recovered the target pressure $p_\infty$,
indicating that the decay of the initial non-linear pressure pulse is progressing slowly. 
Analogous to the findings from the one-dimensional acoustic pulse, any imposition of a hard boundary condition leads to a full/partial reflection of the pressure pulse 
depending on the relaxation amount in the relaxation matrix $\hat{\mat{D}}$, as can be seen in Figure~\ref{fig::p2d_0.5} and in Figure~\ref{fig::p2d_0.5_snapshots}. 
The subsonic outflow boundary condition (\ref*{leg::p2d::so} SO) and the characteristic
boundary conditions (\ref*{leg::p2d::grcbc_so} GRCBC$_{p_\infty}$, \ref*{leg::p2d::grcbc_ma} GRCBC$_{\rho u_\infty}$, \ref*{leg::p2d::grcbc_te} GRCBC$_{\theta_\infty}$) generate 
the most amount of perturbations, although the latter can be tuned appropriately to reduce the reflections by adjusting the relaxation matrix $\hat{\mat{D}}$.
Remarkably, the boundary condition (\ref*{leg::p2d::nscbc_so} NSCBC$_{0.28 p_\infty}$) with the pressure target state $\vec{U}_{p_\infty}$ and the optimal determined
relaxation parameter $\sigma_{opt} = 0.28$, see Subsection~\ref{sec::subsonic_outflow}, shows the smallest reflection of the pressure pulse. Boundary conditions using the far-field state $\vec{U}_\infty$ (\ref*{leg::p2d::ff} FF, \ref*{leg::p2d::grcbc_ff} GRCBC$_\infty$)
seem to be an appropriate choice for a pressure pulse, as they exhibit little reflection and are not sensitive to the alternating change of inflow and outflow, even though
in the interval $t_1 < t < 2$ they dampen the numerical solution slightly more than the reference solution. This overdamping is likely caused by the absence of transverse terms in the boundary conditions,
which are crucial for the accurate propagation of the pressure pulse within the domain. This effect is evident in the discrepancy between the reference solution and the simulation 
results at the corners of the domain (see snapshot $t_2$ in Figure~\ref{fig::p2d_0.5_snapshots}).

\subsection{Zero-circulation vortex}
\label{sec::experiments::zero_circulation_vortex}

\graphicspath{{./figures/zero_circulation_vortex/data}}
\pgfplotsset{table/search path={./figures/zero_circulation_vortex/data}}

The aim of this benchmark is to investigate the non-reflecting behavior of newly implemented boundary conditions in
presence of vortical structures in a low-Mach flow. Particularly, we expect a significant improvement in the accuracy of the 
solution on outflow boundaries, in combination with the
transverse relaxation parameter~$\beta$, emerging from a low-Mach asymptotic analysis as described in
Section~\ref{sec::compressible_equations}.

\begin{figure}[htb]
    \centering
    \tikzsetnextfilename{zero_circulation_vortex_domain}
\begin{tikzpicture}[scale=2, >=latex]
    \draw[thick] (-1, -0.5) rectangle (1, 0.5);
    \fill[amber!30] (-1, -0.5) rectangle (1, 0.5);

    \node at (-1.2, 0) {$\Gamma_{in}$};
    \node at (1.2, 0) {$\Gamma_{out}$};
    \node at (0, 0.63) {$\Gamma_{in}$};
    \node at (0, -0.63) {$\Gamma_{in}$};

    \shade[outer color=amber!30, inner color=matyellow] (0, 0) circle (0.1);
    
    \node at (1.15, -0.4) {$\Omega$};
    \draw[-{>[flex=0.85]}] (0.1, 0) arc (0:315:0.1);
    \draw[-{>[flex=0.85]}] (0.15, 0) arc (0:315:0.15);
    \draw[-{>[flex=0.85]}] (0.2, 0) arc (0:315:0.2);

    \draw[thick, ->] (-1, -0.5) -- (-0.75, -0.5);
    \draw[thick, ->] (-1, -0.25) --(-0.75, -0.25) ;
    \draw[thick, ->] (-1, 0) --    (-0.75, 0);
    \draw[thick, ->] (-1, 0.25) -- (-0.75, 0.25);
    \draw[thick, ->] (-1, 0.5) --  (-0.75, 0.5);
    \node at (-0.925, 0.63) {$\Ma_\infty$};

    \begin{scope}[scale=0.5, xshift=180]

        \fill[amber!30] (-2, -1) rectangle (2, 1);
        \node at (2.3, -0.8) {$\mesh$};

        \pgfplotstableread[col sep=comma]{mesh.csv}\mesh
        \pgfplotstablegetrowsof{\mesh}
        \pgfmathsetmacro{\rows}{\pgfplotsretval-1}

        \foreach \row in {0,...,\rows} {
                \pgfplotstablegetelem{\row}{x0}\of\mesh
                \let\xi=\pgfplotsretval
                \pgfplotstablegetelem{\row}{y0}\of\mesh
                \let\yi=\pgfplotsretval
                \pgfplotstablegetelem{\row}{x1}\of\mesh
                \let\xj=\pgfplotsretval
                \pgfplotstablegetelem{\row}{y1}\of\mesh
                \let\yj=\pgfplotsretval
                \pgfplotstablegetelem{\row}{x2}\of\mesh
                \let\xk=\pgfplotsretval
                \pgfplotstablegetelem{\row}{y2}\of\mesh
                \let\yk=\pgfplotsretval
                \draw (\xi,\yi) -- (\xj,\yj) -- (\xk,\yk) -- cycle;
            }
    \end{scope}

\end{tikzpicture}
    \caption[]{Representation of the domain $\Omega$ with corresponding boundaries $\Gamma$ and the triangulation $\mesh$ with mesh size $h = 0.15$ and local mesh size $h_{out} = 0.1$ at $\Gamma_{out}$.}
    \label{fig::vortex_domain}
\end{figure}
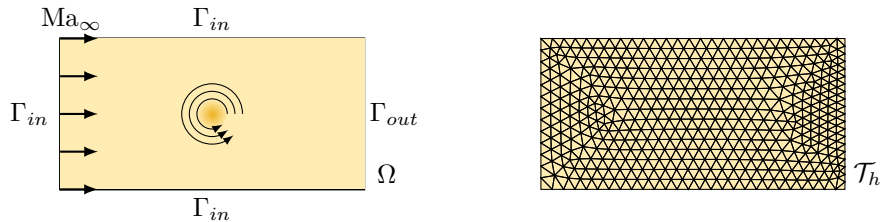

A cylindrical, zero-circulation vortex of radius $R=0.1$ is injected into an inviscid baseflow with Mach number $\Ma_\infty$. 
The domain consists of a rectangle $\Omega=(-20R, 20R)\times(-10R, 10R)$ and in the
initial configuration $t=0$ the vortex is located at the origin $\vec{x}_0 = (0, 0)$, as shown in Figure~\ref{fig::vortex_domain}.
The initial conditions are
\begin{align*}
    \rho_0 & = \rho_\infty \left(1 - \frac{\gamma - 1}{2} \Ma_\infty^2 \alpha^2  \exp\left(1 - r^2/R^2\right) \right)^{1/(\gamma - 1)},                                               \\
    \VEL_0 & = |\VEL_\infty| \left[\begin{pmatrix} 1 \\ 0 \end{pmatrix} + \begin{pmatrix} -y \\ x\end{pmatrix} \frac{\alpha}{R}  \exp\left(\frac{1}{2} \left(1 - r^2/R^2 \right)\right) \right], \\
    p_0    & = p_\infty \left(1 - \frac{\gamma - 1}{2} \Ma_\infty^2 \alpha^2 \exp\left(1 - r^2/R^2\right) \right)^{\gamma/(\gamma - 1)},
\end{align*}
where the polar coordinate $r$ is given as $r^2= (x-x_0)^2 + (y-y_0)^2$ and $\alpha$ denotes the vortex strength ratio of the vortex to the base flow, $\Ma_v = \alpha \Ma_\infty$, studied in the following subsections. We employ aerodynamic
scaling, see \eqref{eq::aerodynamic-scaling}, with dimensionless quantities $\rho_\infty=1$, $|\VEL_\infty| = 1$ and $p_\infty = 1/(\Ma_\infty^2 \gamma)$. The time-stepping is performed with a time step $\Delta t = 2\cdot 10^{-3}$
up to a dimensionless time $t_{end}=5$.

We expect that most reflections will be generated while the vortex intersects the outflow boundary~$\Gamma_{out}$. Dimensional analysis indicates that this intersection occurs at approximately $t \approx 2$. Consequently, the primary objective of this benchmark is to evaluate the non-reflecting behavior of the subsonic outflow boundary conditions.

At the subsonic inflows $\Gamma_{in}$ we expect upstream traveling acoustic waves. Given the good results of the far-field target state $\vec{U}_\infty$ for acoustic perturbations in previous benchmarks, we impose a GRCBC$_{0.01 \infty}$ as the closure condition.

Since the transverse relaxation parameter $\beta$ is derived from a
low-Mach asymptotic analysis, we consider an almost incompressible baseflow with a Mach number of
$\Ma_\infty = 0.03$ and the commonly used limiting case for incompressible flow at $\Ma_\infty = 0.3$.

An exact solution can be constructed by replacing $\vec{x} = \VEL_\infty t + \vec{x}_0$, since the vortex propagates at the speed and in the direction of the base flow $\VEL_\infty$ due to the local equilibrium.
In order to visualize the exact solution in our results, additionally to the scaled $L_2(\Omega)$-errors of the pressure $e_p$, we show the scaled fluctuating pressure $\Delta \tilde{p}$ 
defined in \eqref{eq::scaled_fluctuating_pressure}.

Simulations conducted with increased relaxation parameters $\hat{\mat{D}}$ confirm the observations from previous benchmarks, with higher relaxation increasingly reproducing
the solutions of the common boundary conditions (\ref*{leg::vw::ff} FF) and (\ref*{leg::vw::so} SO). Consequently, we only show results obtained with low relaxation parameters.

\begin{figure}[htb]
    \centering
    \input{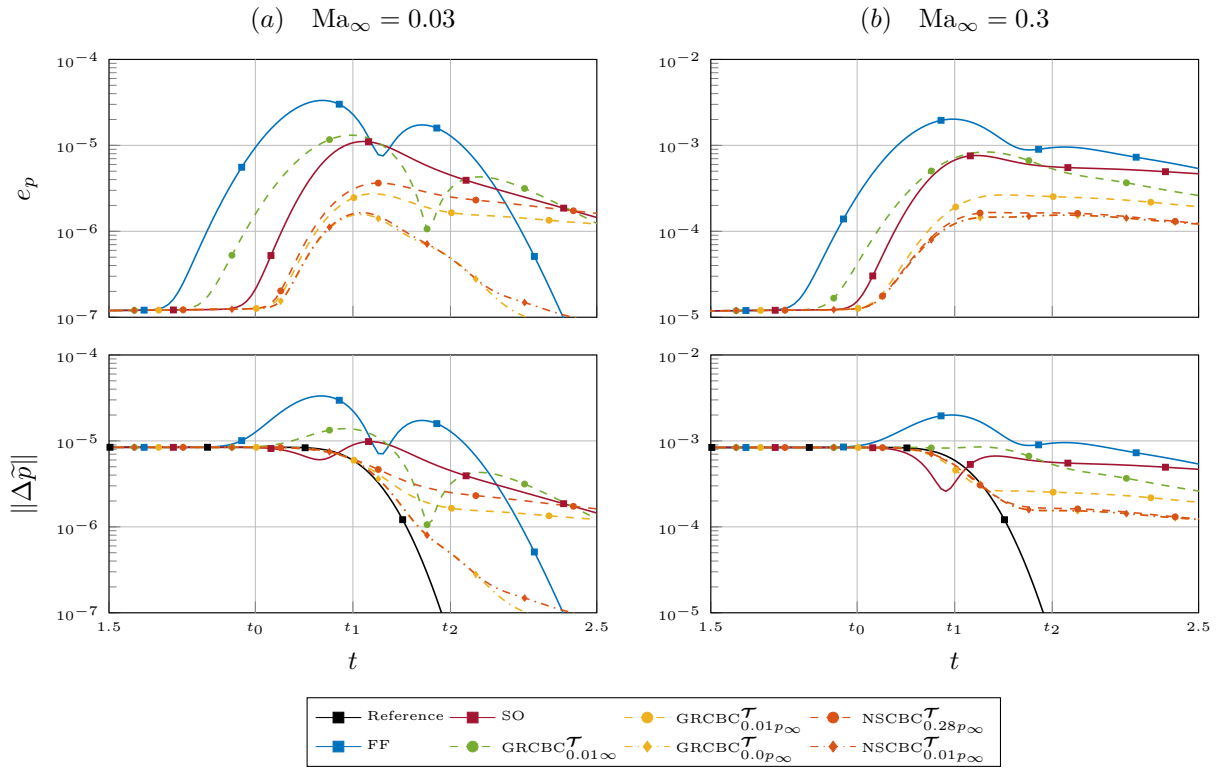}
    \caption[Relative $L_2(\Omega)$-errors of pressure and velocity fields for a weak vortex]
    {Zero-circulation vortex with vortex strength ratio $\alpha = 1/3$: scaled $L_2(\Omega)$-errors of pressure $e_p$ and scaled pressure fluctuations $\Delta \tilde{p}$ for Mach numbers $(a)$ $\Ma_\infty = 0.03$ and $(b)$ $\Ma_\infty = 0.3$.}
    \label{fig::low_vortex_strength_ratio}
\end{figure}

\subsubsection[Low vortex strength ratio]{Low vortex strength ratio $\alpha = 1/3$}
\label{sec::weak_vortex}

\begin{figure}[htb]
    \centering
\begin{tikzpicture}
    \pgfplotsset{every axis/.append style={
        width=1.5cm, 
        height=3cm, 
        mark repeat = 15,
        mark size = 1pt,
        mark options={solid},
        tick label style = {font=\tiny},
        scale only axis,
        xmin=10, 
        xmax=20, 
        ymin=-10,
        ymax=10,
        ytick = {-10, 0, 10},
        xtick = {10, 15, 20}, 
        yticklabels = {$-10$, $\frac{y}{R}$, $10$},
        xticklabels = {$10$, $\frac{x}{R}$, $20$}, 
        xtick pos=bottom,
        xmajorgrids=true, ymajorgrids=true,
        title style={font=\tiny, align=center},
        ylabel style={font=\small},
         },
        }


        \begin{groupplot}[group style={group size=6 by 3, x descriptions at=edge bottom, vertical sep=0.5cm, horizontal sep=0.5cm}, colorbar to name=weak_vortex_colorbar,   
         colorbar style={title=$\omega$, height=\pgfkeysvalueof{/pgfplots/parent axis height}, width=0.3cm, ytick = {-1.5, 0, 5, 10}, 
         ytick pos=right, yticklabel style={anchor=west, xshift=0.5em}}, point meta max = 10, point meta min = -1.5]
            
            \nextgroupplot[title={\ref*{leg::vw::ref}\\Reference}, ylabel=$t_0$, xticklabel=\empty]
            \addplot graphics[xmin=10,ymin=-10,xmax=20,ymax=10] {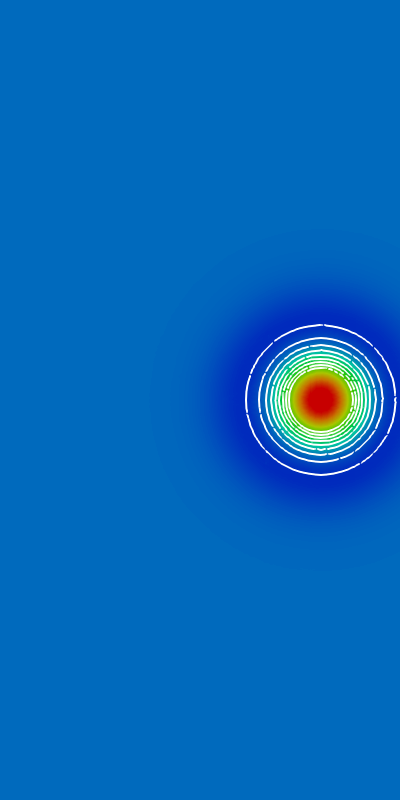};

            \nextgroupplot[title={\ref*{leg::vw::ff} \\ FF}, yticklabel=\empty, xticklabel=\empty]
            \addplot graphics[xmin=10,ymin=-10,xmax=20,ymax=10] {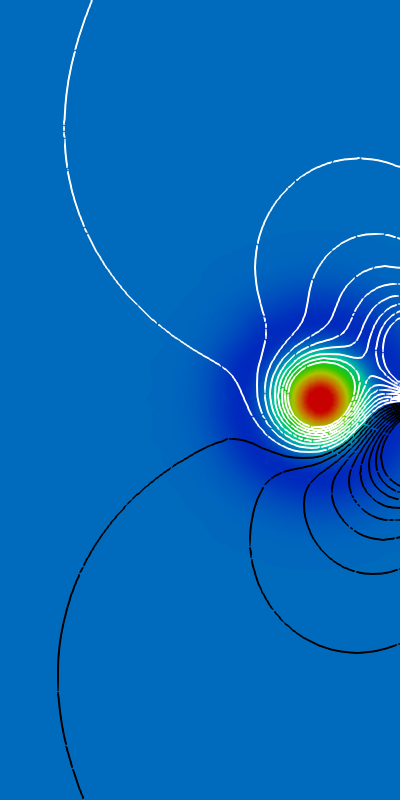};

            \nextgroupplot[title={\ref*{leg::vw::so} \\SO}, yticklabel=\empty, xticklabel=\empty]
            \addplot graphics[xmin=10,ymin=-10,xmax=20,ymax=10] {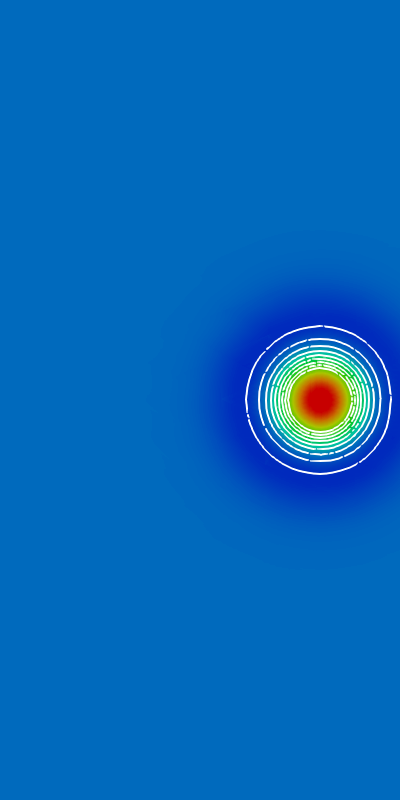};

            \nextgroupplot[title={\ref*{leg::vw::grcbc_ff} \\GRCBC$^{\TRA}_{0.01 \infty}$}, yticklabel=\empty, xticklabel=\empty]
            \addplot graphics[xmin=10,ymin=-10,xmax=20,ymax=10] {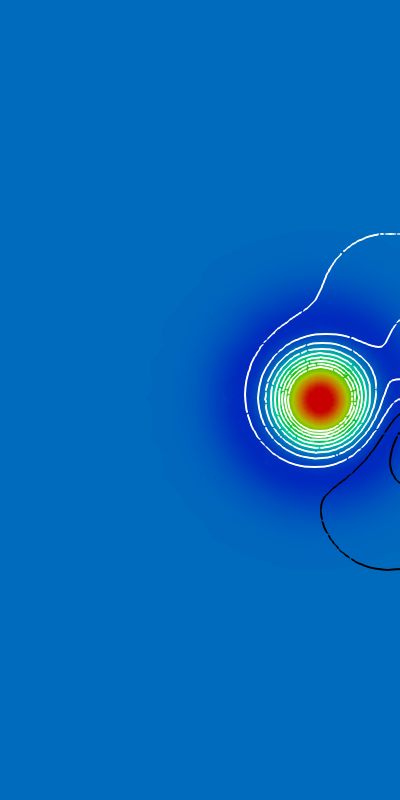};

            \nextgroupplot[title={\ref*{leg::vw::grcbc_so} \\GRCBC$^{\TRA}_{0.01p_{\infty}}$}, yticklabel=\empty, xticklabel=\empty]
            \addplot graphics[xmin=10,ymin=-10,xmax=20,ymax=10] {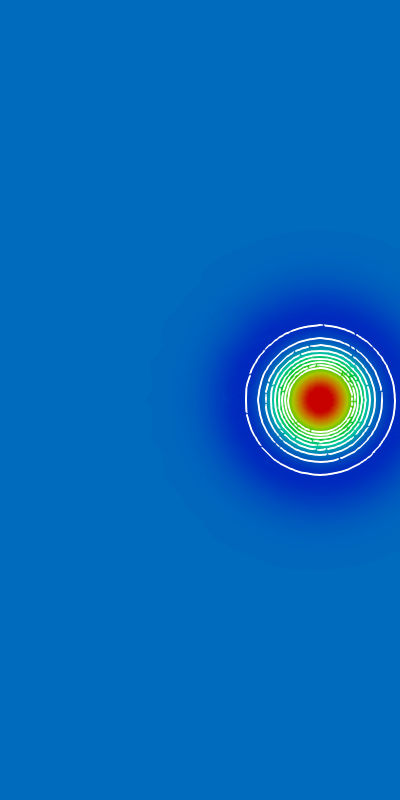};

            \nextgroupplot[title={\ref*{leg::vw::nscbc} \\NSCBC$^{\TRA}_{0.28p_{\infty}}$}, colorbar, yticklabel=\empty, xticklabel=\empty]
            \addplot graphics[xmin=10,ymin=-10,xmax=20,ymax=10] {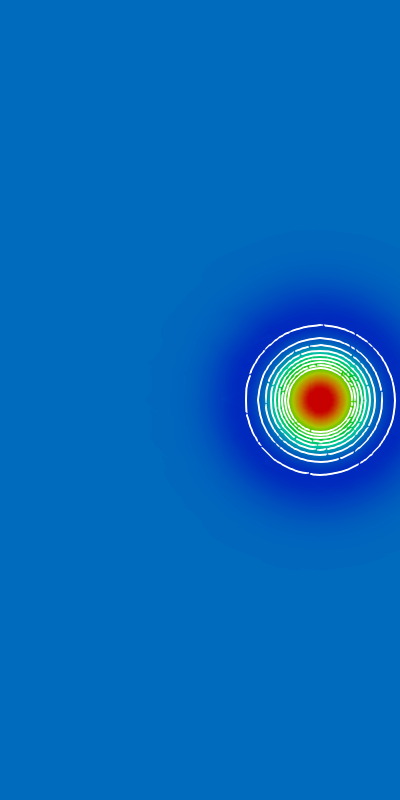};

            \coordinate (top) at (rel axis cs:1,0.75);

            \nextgroupplot[ylabel=$t_1$, xticklabel=\empty]
            \addplot graphics[xmin=10,ymin=-10,xmax=20,ymax=10] {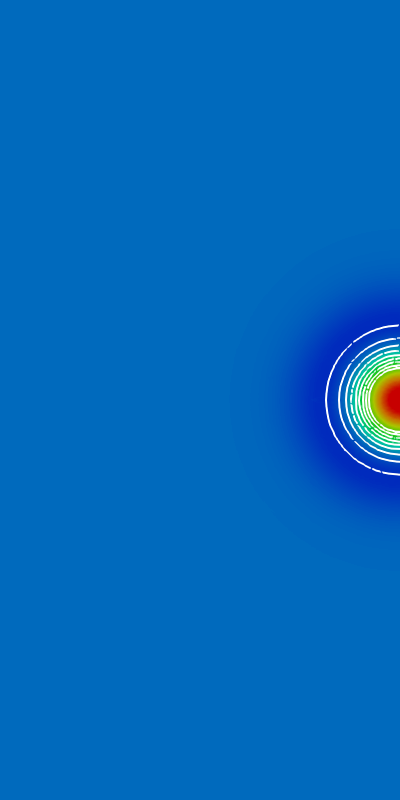};

            \nextgroupplot[xticklabel=\empty, yticklabel=\empty,]
            \addplot graphics[xmin=10,ymin=-10,xmax=20,ymax=10] {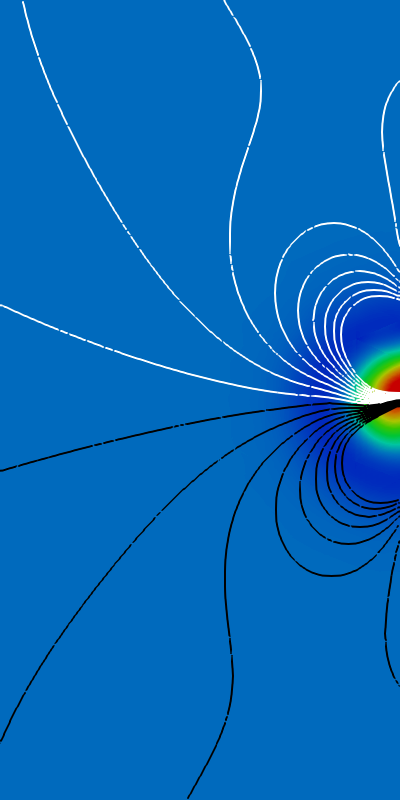};

            \nextgroupplot[yticklabel=\empty, xticklabel=\empty]
            \addplot graphics[xmin=10,ymin=-10,xmax=20,ymax=10] {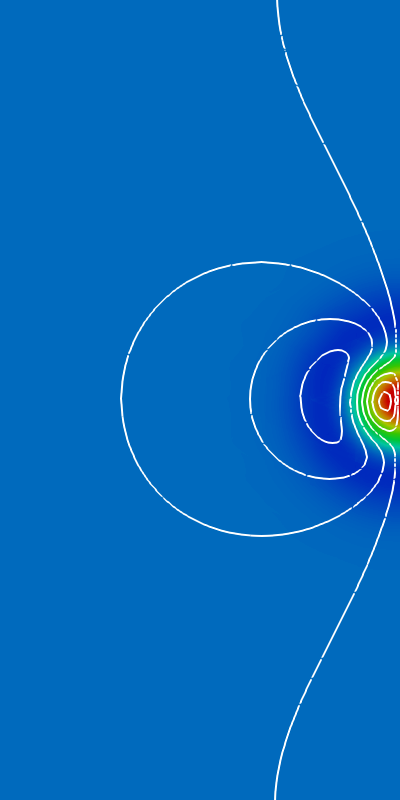};

            \nextgroupplot[yticklabel=\empty, xticklabel=\empty]
            \addplot graphics[xmin=10,ymin=-10,xmax=20,ymax=10] {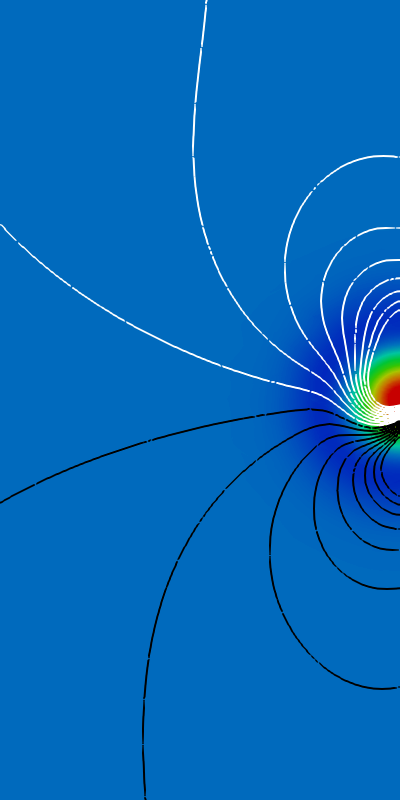};

            \nextgroupplot[yticklabel=\empty, xticklabel=\empty]
            \addplot graphics[xmin=10,ymin=-10,xmax=20,ymax=10] {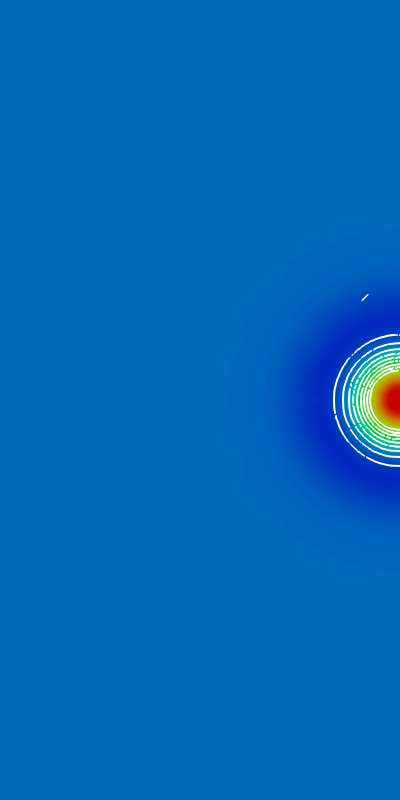};

            \nextgroupplot[yticklabel=\empty, xticklabel=\empty]
            \addplot graphics[xmin=10,ymin=-10,xmax=20,ymax=10] {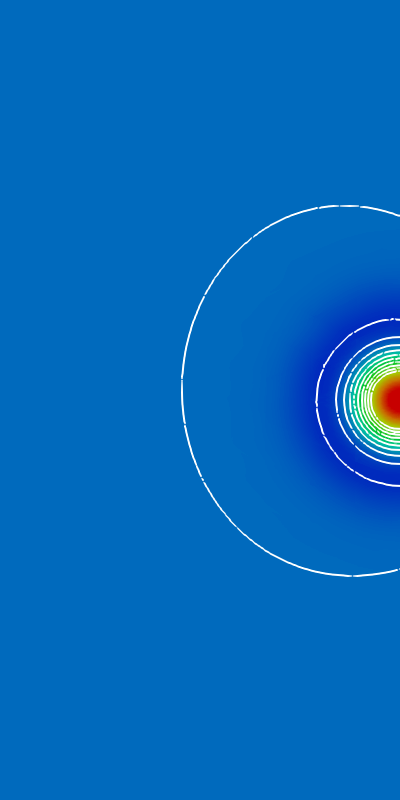};

            \coordinate (mid1) at (rel axis cs:1,0.75);
            \coordinate (mid2) at (rel axis cs:1,0.25);

            \nextgroupplot[ylabel=$t_2$]
            \addplot graphics[xmin=10,ymin=-10,xmax=20,ymax=10] {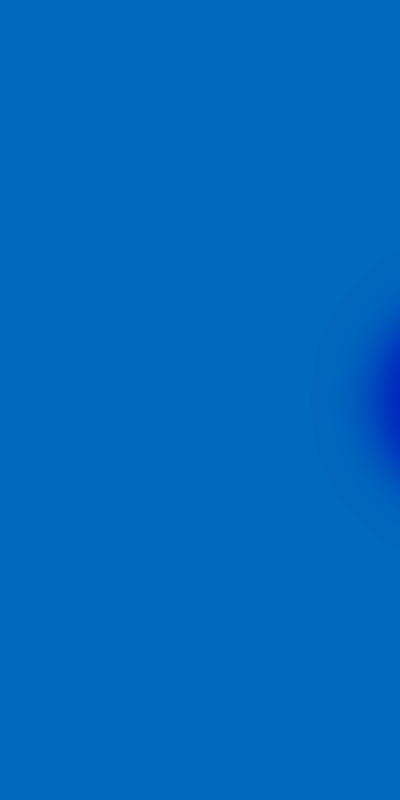};

            \nextgroupplot[yticklabel=\empty]
            \addplot graphics[xmin=10,ymin=-10,xmax=20,ymax=10] {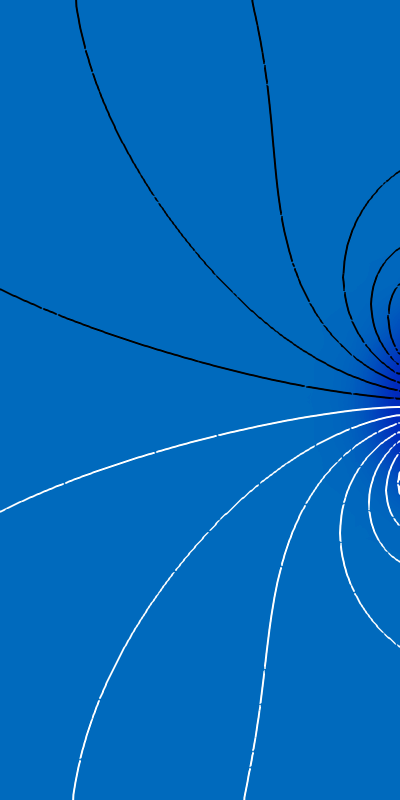};

            \nextgroupplot[yticklabel=\empty]
            \addplot graphics[xmin=10,ymin=-10,xmax=20,ymax=10] {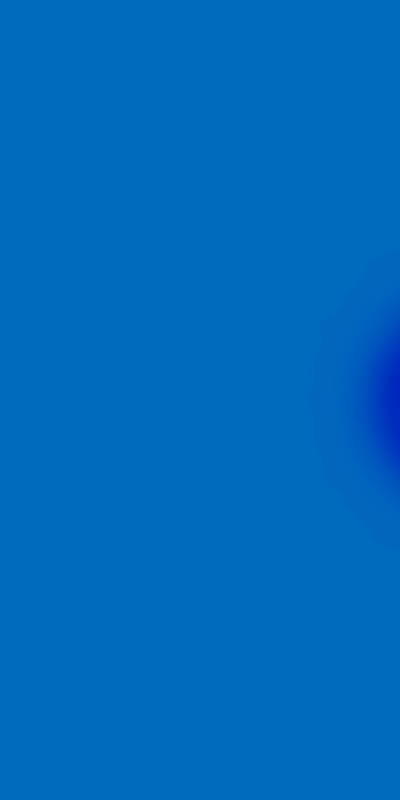};

            \nextgroupplot[yticklabel=\empty]
            \addplot graphics[xmin=10,ymin=-10,xmax=20,ymax=10] {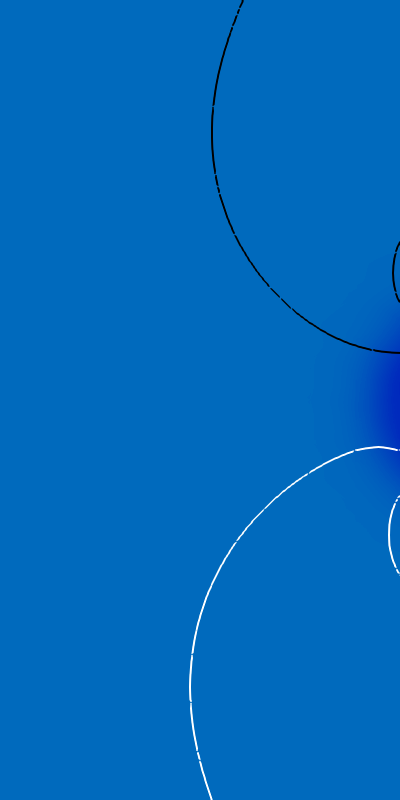};

            \nextgroupplot[yticklabel=\empty]
            \addplot graphics[xmin=10,ymin=-10,xmax=20,ymax=10] {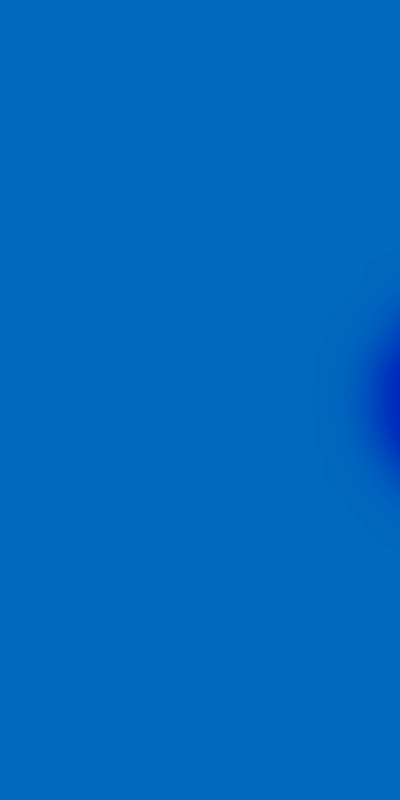};

            \nextgroupplot[yticklabel=\empty]
            \addplot graphics[xmin=10,ymin=-10,xmax=20,ymax=10] {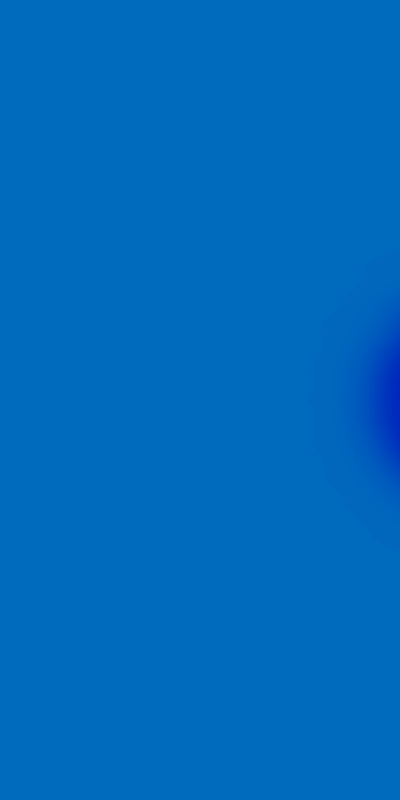};

            \coordinate (bottom) at (rel axis cs:1,0.25);

        \end{groupplot}

        \coordinate (right) at (current bounding box.east);

        \path (top)--(mid1) coordinate[midway] (center1);
        \path (mid2)--(bottom) coordinate[midway] (center2);

        \begin{axis}[%
            hide axis,
            scale only axis,
            colorbar,
            at={(center1 -| right)},
            anchor=west,
            xshift=0cm,
            colorbar style={title=$\omega$, 
            height=\pgfkeysvalueof{/pgfplots/parent axis height}, 
            width=0.3cm, 
            ytick = {-2.5, 0, 5, 10}, 
            yticklabels = {-2.5, 0, 5, 10}, 
            ytick pos=right, 
            yticklabel style={anchor=west, xshift=0em}}, 
            point meta max = 10, point meta min = -2.5
          ]
            \addplot [draw=none] coordinates {(0,0)};
        \end{axis}
        
        \begin{axis}[
            hide axis,
            scale only axis,
            colorbar,
            colormap name=xray,
            at={(center2 -| right)},
            anchor=west,
            xshift=0cm,
            point meta min=-1e-4,
            point meta max=1e-4,
            colorbar style={title=$\Delta \tilde{p}$, 
            height=\pgfkeysvalueof{/pgfplots/parent axis height}, 
            width=0.3cm, 
            ytick = {-1e-4, 0, 1e-4}, 
            yticklabels = {$-10^{-4}$, 0, $10^{-4}$}, 
            scaled y ticks=false,
            ytick pos=right, 
            yticklabel style={anchor=west, xshift=0em}}, 
          ]
            \addplot [draw=none] coordinates {(0,0)};
        \end{axis}

\end{tikzpicture}
    \caption[Snapshots weak vortex]{Zero-circulation vortex with vortex strength ratio $\alpha = 1/3$ and Mach number $\Ma_\infty = 0.03$: vorticity field $\omega$ and isocontour lines of the scaled fluctuating pressure~$\Delta \tilde{p}$ for different boundary conditions at time instances
        $t_0 = 1.8$, $t_1 = 2.0$ and $t_2 = 2.2$}
    \label{fig::vortex_weak_snapshots}
\end{figure}
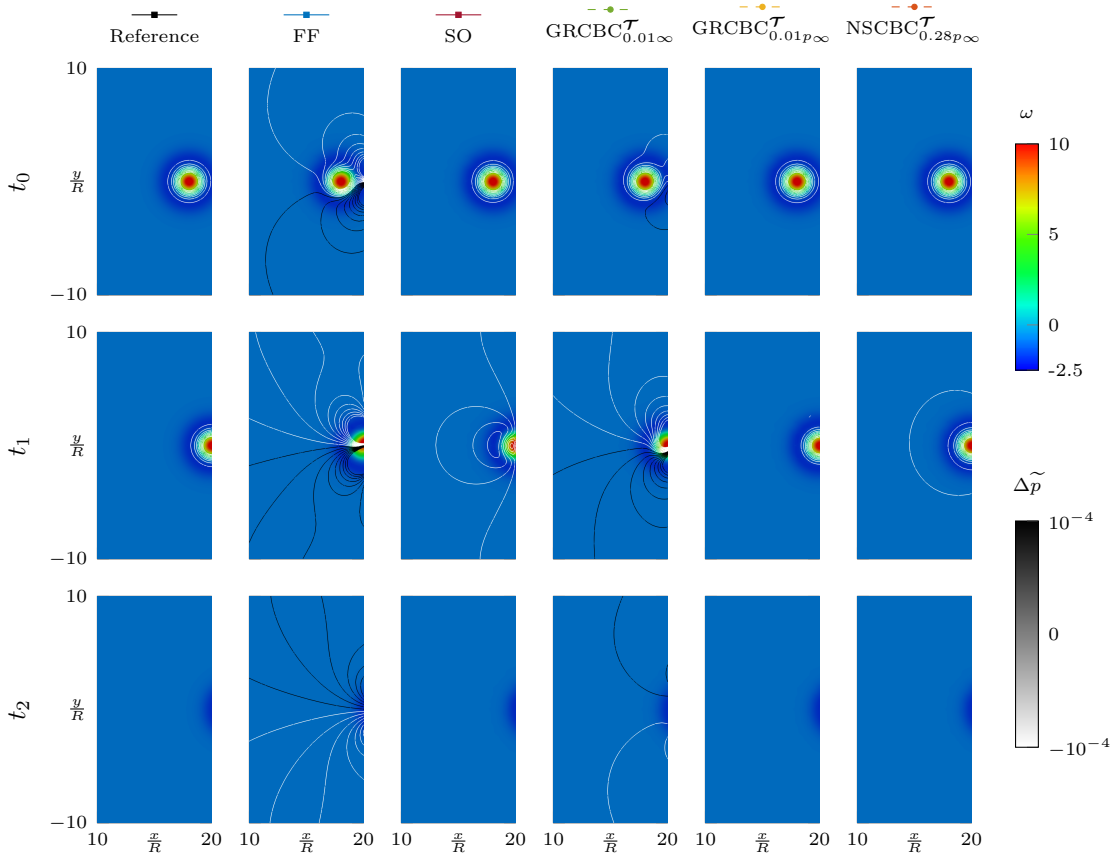

The first test consists of a vortex with a relatively low strength ratio $\alpha = 1/3$, in the sense that
the subsonic outflow does not alternate to a subsonic inflow during the intersection of the vortex. Numerical results are shown in Figure~\ref{fig::low_vortex_strength_ratio},
in which the scaled $L_2(\Omega)$-errors of the pressure $e_p$ and the scaled pressure fluctuation $\Delta \tilde{p}$ are shown for the two distinct baseflow
Mach numbers $\Ma_\infty$. In order to visually assess the performance of the boundary conditions,
snapshots of the vorticity field and isocontour lines of the scaled fluctuating pressure $\Delta \tilde{p}$ are shown for three time instances
$t_0 = 1.8$, $t_1 = 2.0$ and $t_2 = 2.2$ and Mach number $\Ma_\infty = 0.03$ in Figure~\ref{fig::vortex_weak_snapshots}.

As explained and observed in the earlier benchmarks, the subsonic outflow boundary condition (\ref*{leg::vw::so} SO) causes considerable reflections at the outflow boundary $\Gamma_{out}$.
Notably, the reflections are less intense than those generated by the far-field boundary condition~(\ref*{leg::vw::ff} FF), which indicates
that the subsonic outflow target state $\vec{U}_{p_\infty}$ may be more appropriate in this scenario.

In comparison to the far-field boundary condition (\ref*{leg::vw::ff} FF), the CBCs with target state $\vec{U}_\infty$ (\ref*{leg::vw::grcbc_ff} GRCBC$^{\TRA}_\infty$)
show an improvement as they produce less reflections due to the inclusion of the transverse relaxation. Nevertheless, the imposition of the far-field target state
$\vec{U}_\infty$ on a subsonic outflow with exiting vortical structures is inappropriate and leads to a significant error in the solution, which is consistent with the
findings of \cite{pirozzoliGeneralizedCharacteristicRelaxation2013}.

A notable improvement is achieved by the CBCs with target state $\vec{U}_{p_\infty}$ (\ref*{leg::vw::grcbc_so}~GRCBC$^{\TRA}_{p_{\infty}}$, \ref*{leg::vw::nscbc} NSCBC$^{\TRA}_{p_{\infty}}$) along with an additional transverse
relaxation parameter $\beta$. In the plots of the scaled pressure fluctuation norm $\| \Delta \tilde{p} \|$ in Figure~\ref{fig::low_vortex_strength_ratio}, it is evident that these CBCs follow the exact
solution closely during the crossover of the vortex, which is further supported by the snapshots of the vorticity field
and the isocontour lines of the scaled fluctuating pressure $\Delta \tilde{p}$ in Figure~\ref{fig::vortex_weak_snapshots}.

These observations hold for both baseflow Mach numbers $\Ma_\infty$, although for increasing Mach number the deviations from the exact solution become more pronounced, indicating
the eventual limit of the low-Mach asymptotic analysis for the relaxation parameter $\beta$.

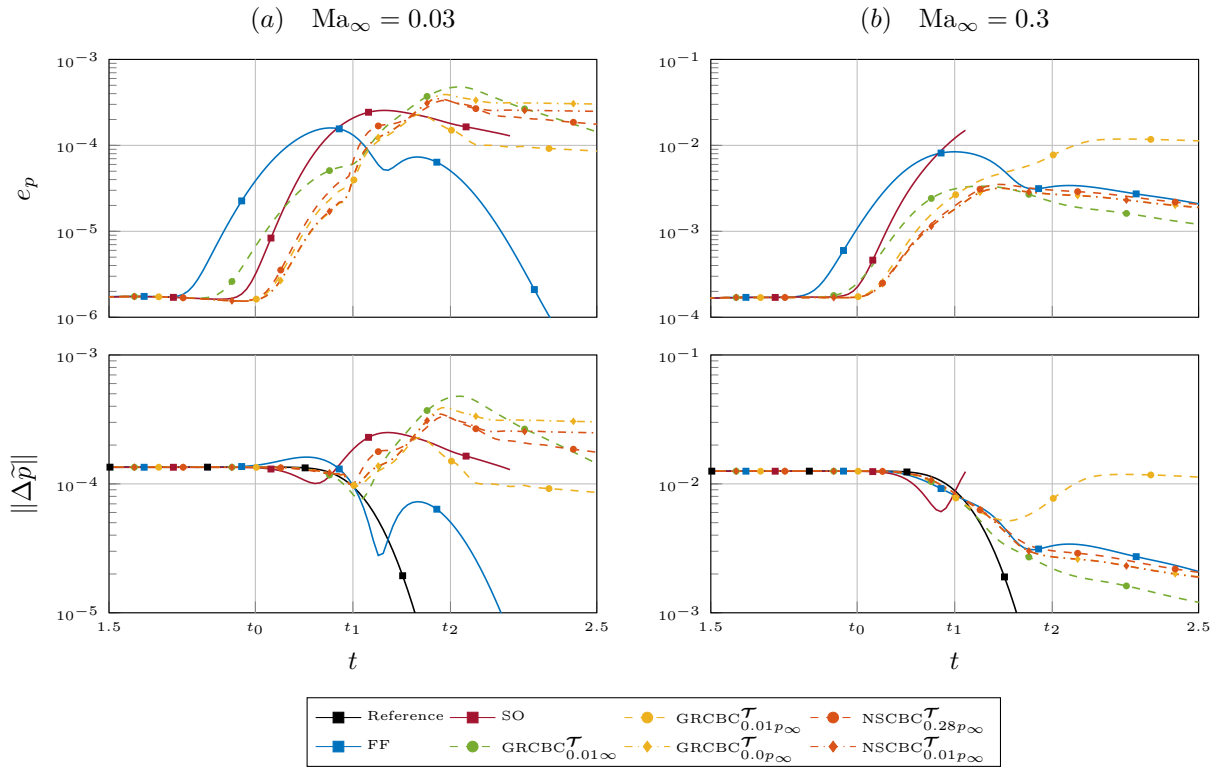
\begin{figure}[htb]
    \centering
    \tikzsetnextfilename{high_vortex_strength_ratio}
\begin{tikzpicture}
    \pgfplotsset{every axis/.append style={
                width=8cm,
                height=5cm,
                mark repeat = 20,
                mark size = 1pt,
                mark options={solid},
                xmin=1.5,
                xmax=2.5,
                max space between ticks=20,
                semithick,
                xtick pos=bottom,
                xmajorgrids=true, ymajorgrids=true,
                axis line style={thin},
                tick label style = {font=\tiny},
                ytick pos=left,
                xtick = {1.5, 2.5}},
        legend style={font=\tiny},
        legend cell align={left},
        legend columns=2,
        transpose legend,
        extra x ticks={1.8, 2, 2.2},
        extra x tick style={grid=major},
        extra x tick labels={$t_0$, $t_1$, $t_2$}
    }

    \begin{groupplot}[group style={group size=2 by 2, horizontal sep=1.5cm, vertical sep=0.5cm}]

        \nextgroupplot[ymode=log, title={$(a) \quad \Ma_\infty=0.03$}, ylabel=$e_p$, xticklabels=\empty, extra x tick labels=\empty, ymax=1e-3, ymin=1e-6]

        \addplot[black,    mark phase = 0, mark=square*, mark options={solid}] table [y=e_p, x=t]{exact_Ma0.03_Mat0.04_dt0.002.dat};
        \addplot[matblue,  mark phase = 8, mark=square*, mark options={solid}]  table [y=e_p, x=t]{standard_farfield_Ma0.03_Mat0.04_dt0.002.dat};
        \addplot[matred, mark phase = 14, mark=square*, mark options={solid}] table [y=e_p, x=t]{standard_outflow_Ma0.03_Mat0.04_dt0.002.dat};

        \addplot[matgreen, mark phase = 6, mark=*, dashed, mark options={solid}]          table [y=e_p_0.01, x=t]{grcbc_farfield_Ma0.03_Mat0.04_dt0.002.dat};

        \addplot[matyellow, mark phase = 11, mark=*, dashed, mark options={solid}]        table [y=e_p_0.01, x=t]{grcbc_outflow_Ma0.03_Mat0.04_dt0.002.dat};
        \addplot[matyellow, mark phase = 16, mark=diamond*, dash dot, mark options={solid}] table [y=e_p_0.0, x=t]{grcbc_outflow_Ma0.03_Mat0.04_dt0.002.dat};

        \addplot[matorange, mark phase = 16, mark=*, dashed, mark options={solid}]   table [y=e_p_0.28, x=t]      {nscbc_outflow_Ma0.03_Mat0.04_dt0.002.dat};
        \addplot[matorange, mark phase = 6, mark=diamond*, dash dot, mark options={solid}] table [y=e_p_0.01, x=t]{nscbc_outflow_Ma0.03_Mat0.04_dt0.002.dat};

        \nextgroupplot[ymode=log, title={$(b) \quad \Ma_\infty=0.3$}, ymax=1e-1, ymin=1e-4, xticklabels=\empty, extra x tick labels=\empty]

        \addplot[black,     mark phase = 0,  mark=square*, mark options={solid}] table [y=e_p, x=t]{exact_Ma0.3_Mat0.4_dt0.002.dat};
        \addplot[matblue,   mark phase = 8,  mark=square*, mark options={solid}]  table [y=e_p, x=t]{standard_farfield_Ma0.3_Mat0.4_dt0.002.dat};
        \addplot[matred,  mark phase = 14, mark=square*, mark options={solid}] table [y=e_p, x=t]{standard_outflow_Ma0.3_Mat0.4_dt0.002.dat};

        \addplot[matgreen, mark phase = 6, mark=*, dashed, mark options={solid}]          table [y=e_p_0.01, x=t]{grcbc_farfield_Ma0.3_Mat0.4_dt0.002.dat};

        \addplot[matyellow, mark phase = 11, mark=*, dashed, mark options={solid}]        table [y=e_p_0.01, x=t]{grcbc_outflow_Ma0.3_Mat0.4_dt0.002.dat};
        \addplot[matyellow, mark phase = 16, mark=diamond*, dash dot, mark options={solid}] table [y=e_p_0.0, x=t]{grcbc_outflow_Ma0.3_Mat0.4_dt0.002.dat};

        \addplot[matorange, mark phase = 16, mark=*, dashed, mark options={solid}]   table [y=e_p_0.28, x=t]      {nscbc_outflow_Ma0.3_Mat0.4_dt0.002.dat};
        \addplot[matorange, mark phase = 6, mark=diamond*, dash dot, mark options={solid}] table [y=e_p_0.01, x=t]{nscbc_outflow_Ma0.3_Mat0.4_dt0.002.dat};

        \nextgroupplot[ymode=log, ylabel=$\| \Delta \tilde{p} \|$, ymax=1e-3, ymin=1e-5, xlabel=$t$, legend to name=vortex_strong_legend]

        \addplot[black,     mark phase = 0,  mark=square*, mark options={solid}] table [y=p, x=t]{exact_Ma0.03_Mat0.04_dt0.002.dat}; \addlegendentry{Reference}; \label{leg::vs::ref};
        \addplot[matblue,   mark phase = 8,  mark=square*, mark options={solid}]  table [y=p, x=t]{standard_farfield_Ma0.03_Mat0.04_dt0.002.dat};      \addlegendentry{FF}; \label{leg::vs::ff};
        \addplot[matred,  mark phase = 14, mark=square*, mark options={solid}] table [y=p, x=t]{standard_outflow_Ma0.03_Mat0.04_dt0.002.dat}; \addlegendentry{SO};  \label{leg::vs::so};

        \addplot[matgreen, mark phase = 6, mark=*, dashed, mark options={solid}]          table [y=p_0.01, x=t]{grcbc_farfield_Ma0.03_Mat0.04_dt0.002.dat};     \addlegendentry{GRCBC$_{0.01 \infty}^{\TRA}$}; \label{leg::vs::grcbc_ff};

        \addplot[matyellow, mark phase = 11, mark=*, dashed, mark options={solid}]        table [y=p_0.01, x=t]{grcbc_outflow_Ma0.03_Mat0.04_dt0.002.dat};   \addlegendentry{GRCBC$_{0.01p_{\infty}}^{\TRA}$}; \label{leg::vs::grcbc_so};
        \addplot[matyellow, mark phase = 16, mark=diamond*, dash dot, mark options={solid}] table [y=p_0.0, x=t]{grcbc_outflow_Ma0.03_Mat0.04_dt0.002.dat};   \addlegendentry{GRCBC$_{0.0p_{\infty}}^{\TRA}$};

        \addplot[matorange, mark phase = 16, mark=*, dashed, mark options={solid}]   table [y=p_0.28, x=t]      {nscbc_outflow_Ma0.03_Mat0.04_dt0.002.dat};    \addlegendentry{NSCBC$_{0.28p_{\infty}}^{\TRA}$}; \label{leg::vs::nscbc};
        \addplot[matorange, mark phase = 6, mark=diamond*, dash dot, mark options={solid}] table [y=p_0.01, x=t]{nscbc_outflow_Ma0.03_Mat0.04_dt0.002.dat};    \addlegendentry{NSCBC$_{0.01p_{\infty}}^{\TRA}$};
        \coordinate (left) at (rel axis cs:0,0);

        \nextgroupplot[ymode=log, ymax=1e-1, ymin=1e-3, xlabel=$t$]

        \addplot[black,     mark phase = 0,  mark=square*, mark options={solid}] table [y=p, x=t]{exact_Ma0.3_Mat0.4_dt0.002.dat};
        \addplot[matblue,   mark phase = 8,  mark=square*, mark options={solid}]  table [y=p, x=t]{standard_farfield_Ma0.3_Mat0.4_dt0.002.dat};
        \addplot[matred,  mark phase = 14, mark=square*, mark options={solid}] table [y=p, x=t]{standard_outflow_Ma0.3_Mat0.4_dt0.002.dat};

        \addplot[matgreen, mark phase = 6, mark=*, dashed, mark options={solid}]          table [y=p_0.01, x=t]{grcbc_farfield_Ma0.3_Mat0.4_dt0.002.dat};

        \addplot[matyellow, mark phase = 11, mark=*, dashed, mark options={solid}]        table [y=p_0.01, x=t]{grcbc_outflow_Ma0.3_Mat0.4_dt0.002.dat};
        \addplot[matyellow, mark phase = 16, mark=diamond*, dash dot, mark options={solid}] table [y=p_0.0, x=t]{grcbc_outflow_Ma0.3_Mat0.4_dt0.002.dat};

        \addplot[matorange, mark phase = 16, mark=*, dashed, mark options={solid}]   table [y=p_0.28, x=t]      {nscbc_outflow_Ma0.3_Mat0.4_dt0.002.dat};
        \addplot[matorange, mark phase = 6, mark=diamond*, dash dot, mark options={solid}] table [y=p_0.01, x=t]{nscbc_outflow_Ma0.3_Mat0.4_dt0.002.dat};
        \coordinate (right) at (rel axis cs:1,0);

    \end{groupplot}

    \path (left)--(right) coordinate[midway] (group center);
    \node[inner sep=0pt,yshift=-4.5em] at(group center) {\pgfplotslegendfromname{vortex_strong_legend}};

\end{tikzpicture}
    \caption[Relative $L_2(\Omega)$-errors of pressure and velocity fields for a strong vortex]
    {Zero-circulation vortex with vortex strength ratio $\alpha = 4/3$: scaled $L_2(\Omega)$-errors of pressure $e_p$ and scaled pressure fluctuations $\Delta \tilde{p}$ for Mach numbers $(a)$ $\Ma_\infty = 0.03$ and $(b)$ $\Ma_\infty = 0.3$.}
    \label{fig::high_vortex_strength_ratio}
\end{figure}

\subsubsection[High vortex strength ratio]{High vortex strength ratio $\alpha = 4/3$}
The second vortex is characterized by a higher vortex strength ratio $\alpha = 4/3$, causing flow inversion to
occur at the outflow boundary $\Gamma_{out}$. Results are displayed in Figure~\ref{fig::high_vortex_strength_ratio}, where
scaled $L_2(\Omega)$-errors of the pressure, $e_p$, and the scaled pressure fluctuation, $\Delta \tilde{p}$, are shown for the two distinct baseflow
Mach numbers $\Ma_\infty$.

Since the outflow target state $\vec{U}_{p_\infty}$ is not an appropriate state to handle a subsonic inflow,
the subsonic outflow condition (\ref*{leg::vs::so} SO) is not well-defined anymore and diverges, due to the flow inversion at the outflow boundary.
However, the CBCs  (\ref*{leg::vs::grcbc_so} GRCBC$^{\TRA}_{p_{\infty}}$, \ref*{leg::vs::nscbc} NSCBC$^{\TRA}_{p_{\infty}}$), despite sharing the same target state $\vec{U}_{p_\infty}$,
effectively manage the flow inversion, although producing relatively large reflections.

Boundary conditions (\ref*{leg::vs::ff} FF, \ref*{leg::vs::grcbc_ff} GRCBC$^{\TRA}_{\infty}$) with the far-field target state $\vec{U}_{\infty}$ are by construction capable of
handling flow inversion and seem to dampen the reflections more effectively than the CBCs with target state $\vec{U}_{p_\infty}$.
Yet, the imposition of the far-field target state $\vec{U}_\infty$ on a subsonic outflow boundary in the presence of vortices remains inappropriate.

For this vortex strength ratio, the effect of the transverse terms $\TRA$ seems to be marginal. Although the relaxation
parameter $\beta$ causes a later deviation from the exact solution (see interval $t_0 < t_1$) in comparison to (\ref*{leg::vs::ff} FF) and (\ref*{leg::vs::so} SO), the incoming disturbances due to the flow inversion
seem to be too strong to be effectively damped by the transverse relaxation.

\subsection{Laminar flow around a cylinder}
\label{sec::experiments::cylinder_flow}

\graphicspath{{./figures/flow_around_a_cylinder/data}}
\pgfplotsset{table/search path={./figures/flow_around_a_cylinder/data}}

The final benchmark involves a two-dimensional circular cylinder within a laminar viscous flow, incorporating all features from the previous benchmarks.
Vortices, shed from the cylinder, generate alternating acoustic waves, resulting in a dipole sound, requiring the boundary conditions to be non-reflecting on the whole boundary.
Additionally, the outflow boundary closure must be capable of handling the vortices crossing the artificial boundary without generating reflections.
\begin{figure}[h]
    \centering
    \tikzsetnextfilename{flow_around_a_cylinder_domain}
\begin{tikzpicture}[scale=0.75, >=latex]
    \draw[thick] (-1.5, -2.0) rectangle (2.5, 2.0);
    \fill[amber!30] (-1.5, -2.0) rectangle (2.5, 2.0);
    \draw[thick] (0, 0) circle (0.1);
    \fill[white!30] (0, 0) circle (0.1);

    \node at (-2, 0) {$\Gamma_{in}$};
    \node at (0.5, 2.4) {$\Gamma_{in}$};
    \node at (0.5, -2.4) {$\Gamma_{in}$};
    \node at (3.1, 0) {$\Gamma_{out}$};

    \draw[thick, ->] (-1.5, 2.0) --  (-0.6, 2.0);
    \draw[thick, ->] (-1.5, 1.5) --  (-0.6, 1.5);
    \draw[thick, ->] (-1.5, 0.5) --  (-0.6, 0.5);
    \draw[thick, ->] (-1.5, 1) --    (-0.6, 1) ;
    \draw[thick, ->] (-1.5, 0) --    (-0.6, 0);
    \draw[thick, ->] (-1.5, -0.5) -- (-0.6, -0.5);
    \draw[thick, ->] (-1.5, -1) --   (-0.6, -1);
    \draw[thick, ->] (-1.5, -1.5) --  (-0.6, -1.5);
    \draw[thick, ->] (-1.5, -2.0) --  (-0.6, -2.0);

    \node at (-1.1, 2.4) {$\Ma_\infty$};
    \node at (2.9, -1.7) {$\Omega$};

    \draw[densely dashed, thin] (0,0) circle (0.3);

    \begin{scope}[xshift=200, scale=5]

        \draw[loosely dashed, thick] (0,0) circle (0.4);
        \fill[amber!30] (0, 0) circle (0.4);

        \draw[thick] (0, 0) circle (0.1);
        \fill[white!30] (0, 0) circle (0.1);

        \draw[solid, ->] (0, 0) -- (0.2, 0) node[anchor=west] {$x$};
        \draw[solid, ->] (0, 0) -- (0, 0.2) node[anchor=south] {$y$};
        \draw[solid, ->] (0, 0) -- (0.1414, 0.1414) node[anchor=south west] {$r$};
        \draw[solid, ->] (0.1414,0) arc[start angle=0, end angle=45, radius=0.1414] node [midway, right] {$\varphi$};

        \node at (-0.0, -0.03) {$\Gamma_{cyl}$};
        \draw[solid] (0, 0) -- (-0.0707, 0.0707);
        \draw[solid, ->] (-0.15, 0.15) -- (-0.0707, 0.0707);
        \node at (-0.2, 0.2) {$R$};

    \end{scope}

    \begin{scope}[scale=0.2, xshift=1800]
        \fill[amber!30] (-7.5, -10.0) rectangle (12.5, 10.0);
        \fill[white!30] (0, 0) circle (0.5);
        \node at (14.5, -8.5) {$\mesh$};

        \pgfplotstableread[col sep=comma]{mesh.csv}\mesh
        \pgfplotstablegetrowsof{\mesh}
        \pgfmathsetmacro{\rows}{\pgfplotsretval-1}

        \foreach \row in {0,...,\rows} {
                \pgfplotstablegetelem{\row}{x0}\of\mesh
                \let\xi=\pgfplotsretval
                \pgfplotstablegetelem{\row}{y0}\of\mesh
                \let\yi=\pgfplotsretval
                \pgfplotstablegetelem{\row}{x1}\of\mesh
                \let\xj=\pgfplotsretval
                \pgfplotstablegetelem{\row}{y1}\of\mesh
                \let\yj=\pgfplotsretval
                \pgfplotstablegetelem{\row}{x2}\of\mesh
                \let\xk=\pgfplotsretval
                \pgfplotstablegetelem{\row}{y2}\of\mesh
                \let\yk=\pgfplotsretval
                \draw (\xi,\yi) -- (\xj,\yj) -- (\xk,\yk) -- cycle;
            }

    \end{scope}

\end{tikzpicture}
    \caption{Representation of the domain $\Omega$ with corresponding boundaries $\Gamma$ and the triangulation $\mesh$ with region dependent mesh size $h$.}
    \label{fig::cylinder_domain}
\end{figure}
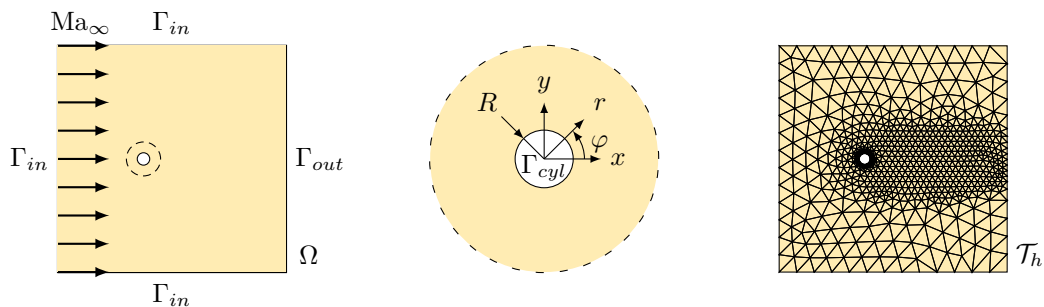

Let $\Omega~=~\left\lbrace (x, y) \in (-15R, 25R) \times (-20R, 20R) \, | \, x^2 + y^2 \geq R^2 \right\rbrace$ be the computational domain 
with a circular cylinder of radius $R=0.5$, as shown in Figure~\ref{fig::cylinder_domain}. The polar coordinates $r, \varphi$ in
counterclockwise direction are given as $r^2=x^2 + y^2$ and $\tan{\varphi} = y/x$.

An adiabatic wall boundary condition, as described in \eqref{eq::adiabatic_wall},
is enforced on the cylinder $\Gamma_{cyl}$. At the inflow boundaries $\Gamma_{in}$, we anticipate the presence of oblique impinging acoustic waves;
therefore, we apply the characteristic boundary condition (GRCBC$_{0.01 \infty}$), using the far-field target state $\vec{U}_\infty$ based on its strong performance in prior benchmarks.
The subsonic outflow boundary $\Gamma_{out}$, where vortices shed from the cylinder interact with sound waves,
poses the greatest challenge to manage, particularly due to its proximity to the cylinder. This boundary is the primary focus of the investigation of this test case.

We consider a Mach number of $\Ma_\infty = 0.2$ and a Reynolds number of $\Re_\infty = 150$. In this setting the dimensionless vortex shedding frequency
is known to be $0.18 \leq St \leq 0.185$ \cite{inoueSoundGenerationTwodimensional2002}, where $St = f D / u_\infty$ is the Strouhal number, $f$ the shedding frequency, $D=2R$ the cylinder diameter and $u_\infty$ the free stream velocity.
Every numerical investigation was initiated with a uniform flow at $t=0$ and advanced over time
$t=250$ with time step $\Delta t = 5\cdot 10^{-2}$, until the initial disturbances exited the domain. For an aerodynamic scaling \eqref{eq::aerodynamic-scaling}, the time step $\Delta t$ corresponds to a Nyquist frequency of $St_{ny} = 10$,
ensuring that the shedding frequency is well resolved. The simulation was then extended for a prolonged period, $\overline{t} = \left\lbrace t \, | \, 250 < t \leq 750 \right\rbrace$,
during which statistical data was collected over multiple shedding cycles.

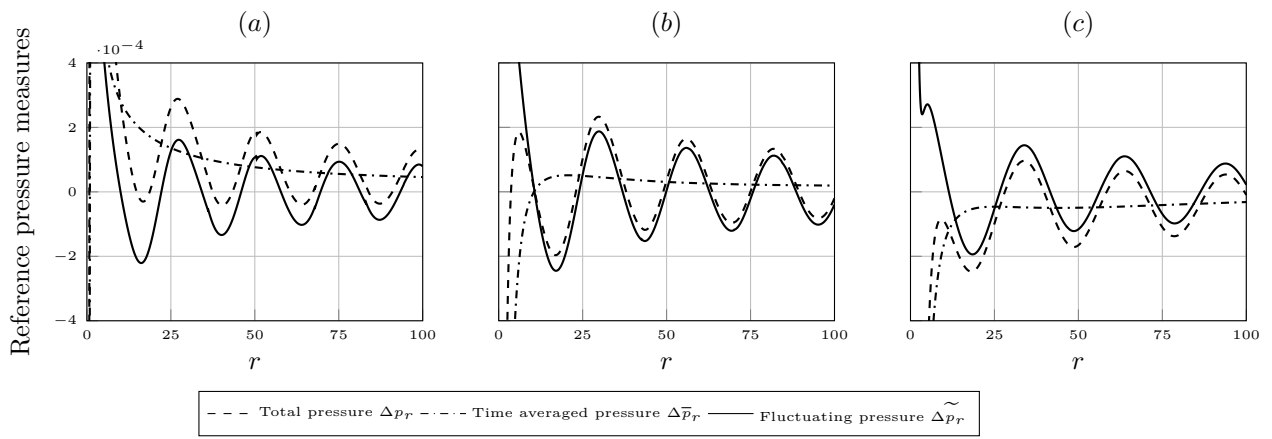
\begin{figure}[h]
    \tikzsetnextfilename{flow_around_a_cylinder_radial_reference_pressure}
\begin{tikzpicture}[]
    \pgfplotsset{every axis/.append style={
        width=6cm, 
        height=5cm, 
        no markers,  
        ymin=-4e-4, 
        ymax=4e-4,
        xmin=0,
        xmax=100,
        xtick={0, 25, 50, 75, 100},
        ytick={-4e-4, -2e-4, 0, 2e-4, 4e-4},
        xlabel=$r$,
        max space between ticks=20,
        thick,
        xtick pos=bottom,
        ytick pos=left,
        xmajorgrids=true, ymajorgrids=true,
        axis line style={thin}, 
        tick label style = {font=\tiny}},
        legend cell align={left},
        legend columns=3,
        legend style={font=\tiny},
    }

    \begin{groupplot}[group style={group size=3 by 1, x descriptions at=edge bottom, vertical sep=0.5cm, horizontal sep=1cm}]
        \nextgroupplot[title={$(a)$}, ylabel={Reference pressure measures}, legend to name=p_radial]
            \addplot[black, dashed] table [y=p, x=r]{reference_angle130.dat}; \addlegendentry{Total pressure $\Delta p_r$};
            \addplot[black, dash dot] table [y=p_mean, x=r]{reference_angle130.dat}; \addlegendentry{Time averaged pressure $\Delta \overline{p}_r$};
            \addplot[black] table [y=p_acou, x=r]{reference_angle130.dat}; \addlegendentry{Fluctuating pressure $\Delta \tilde{p}_r$};
            \coordinate (left) at (rel axis cs:0,0);

        \nextgroupplot[title={$(b)$}, yticklabels=\empty, ytick scale label code/.code={}]
            \addplot[black, dashed] table [y=p, x=r]{reference_angle101.5.dat};    
            \addplot[black, dash dot] table [y=p_mean, x=r]{reference_angle101.5.dat};
            \addplot[black] table [y=p_acou, x=r]{reference_angle101.5.dat};

        \nextgroupplot[title={$(c)$}, yticklabels=\empty, ytick scale label code/.code={}]
            \addplot[black, dashed] table [y=p, x=r]{reference_angle60.dat};    
            \addplot[black, dash dot] table [y=p_mean, x=r]{reference_angle60.dat};
            \addplot[black] table [y=p_acou, x=r]{reference_angle60.dat};
            \coordinate (right) at (rel axis cs:1,0);

    \end{groupplot}
         
    \path (left)--(right) coordinate[midway] (group center);
    \node[right=7em,yshift=-3.5em] at(group center -| current bounding box.west) {\pgfplotslegendfromname{p_radial}};

\end{tikzpicture}
    \caption[reference]{Laminar flow around a cylinder: plot of the reference solution pressure measures at angle $(a)$ $\varphi = 130^\circ$, $(b)$ $\varphi = 101.5^\circ$  and $(c)$ $\varphi = 60^\circ$  over the radial coordinate $r$ at time $t=731.55$.}
    \label{fig::inoue_reference}
\end{figure}

A reference solution has been generated in a reference domain $\Omega_{r}~=~\left\lbrace (x, y) | \, x^2 + y^2 \leq 1500^2 \right\rbrace$
by adding additional sponge layers to properly dampen all incoming disturbances. In Figure~\ref{fig::inoue_reference}, the total pressure $\Delta p$, the time-averaged
pressure $\Delta \overline{p}$ over the time interval $\overline{t}$ and the fluctuating pressure $\Delta \tilde{p}$ of the reference solution,
respectively defined as
\begin{align*}
    \Delta p             := p - p_\infty,                                                \qquad
    \Delta \overline{p}  := \frac{1}{|\overline{t}|} \int_{\overline{t}} \Delta p \, dt, \qquad
    \Delta \tilde{p}     := \Delta p -  \Delta \overline{p},
\end{align*}
are shown for different angles $\varphi$ over the polar coordinate $r$. These reference results closely match the simulation results of a Direct Numerical Solution (DNS) presented
in \cite[Figure 8]{inoueSoundGenerationTwodimensional2002}.

\begin{figure}[htb]
    \tikzsetnextfilename{flow_around_a_cylinder_polar_root_mean_square_pressure}
\begin{tikzpicture}
    \pgfplotsset{compat=newest,
    every axis/.append style={
    width=7cm,
    height=7cm,
    mark repeat = 60,
    mark size = 1pt,
    mark options={solid},
    xtick={0, 45, 90, 135, 180, 225, 270, 315},
    xticklabel=$\pgfmathprintnumber{\tick}^\circ$,
    ytick pos = both,
    thick,
    major tick length=0.1cm,
    axis line style={thin},
    tick label style = {font=\tiny},
    y axis line style={
            overlay,
            rotate around={90:(current axis.origin)},
            yshift=3cm,
            xshift=1cm,
        },
    ytick style={rotate around={90:(current axis.origin)}, xshift=1cm, yshift=3cm},
    yticklabel style={
            rotate=-90,
            rotate around={90:(current axis.origin)},
            xshift=-3cm,
            yshift=1cm,
            left,
        },
    ylabel = {$\Delta \tilde{p}_{\text{rms}}$},
    legend style={font=\tiny},
    legend columns=-1,
    }
    }

    \begin{polaraxis}[at={(0 cm, 0 cm)}, title={$(a) \quad r_0$}, ymin=0, ymax=1e-3, ytick={0, 0.25e-3, 0.5e-3, 0.75e-3, 1e-3},
            ylabel style={xshift=-5cm, yshift=-1.7cm, rotate=90, anchor=south}, legend to name=prms_polar_legend]
        \addplot+[black,    mark phase=0,  mark=square*, mark options={solid}] table [y=p_rms, x=theta] {reference_R4.dat};                              \addlegendentry{Reference}; \label{leg::cyl::ref};
        \addplot+[matblue,  mark phase=11, mark=square*, mark options={solid}] table [y=p_rms, x=theta] {farfield_inflow_and_outflow_R4.dat};                \addlegendentry{FF}; \label{leg::cyl::ff};
        \addplot+[matred, mark phase=21, mark=square*, mark options={solid}] table [y=p_rms, x=theta] {farfield_inflow_and_pressure_outflow_R4.dat};          \addlegendentry{SO}; \label{leg::cyl::so};
        \addplot+[matgreen,   mark phase=31, mark=*, dashed, mark options={solid}] table [y=p_rms, x=theta]{grcbc_farfield_inflow_and_outflow_R4.dat};        \addlegendentry{GRCBC$_{0.01\infty}^{\TRA \VIS}$}; \label{leg::cyl::grcbc_ff};
        \addplot+[matyellow,mark phase=41, mark=*, dashed, mark options={solid}] table [y=p_rms, x=theta]{grcbc_farfield_inflow_and_pressure_outflow_R4.dat};    \addlegendentry{GRCBC$_{0.01p_{\infty}}^{\TRA \VIS}$}; \label{leg::cyl::grcbc_so};
        \addplot+[matorange, mark phase=51, mark=*, dashed, mark options={solid}] table [y=p_rms, x=theta]{nscbc_farfield_inflow_and_pressure_outflow_R4.dat};    \addlegendentry{NSCBC$_{0.28p_{\infty}}^{\TRA \VIS}$}; \label{leg::cyl::nscbc_so};
    \end{polaraxis}

    \begin{polaraxis}[at={(8 cm, 0 cm)}, title={$(b) \quad r_1$}, ymin=0, ymax=5e-4, ytick={0, 1.25e-4, 2.5e-4, 3.75e-4, 5e-4},
            ylabel style={xshift=-5cm, yshift=-1.7cm, rotate=90, anchor=south},]
        \addplot+[black,    mark phase=0,  mark=square*, mark options={solid}] table [y=p_rms, x=theta]{reference_R7.5.dat};
        \addplot+[matblue,  mark phase=11, mark=square*, mark options={solid}] table [y=p_rms, x=theta]{farfield_inflow_and_outflow_R7.5.dat};
        \addplot+[matred, mark phase=21, mark=square*, mark options={solid}] table [y=p_rms, x=theta]{farfield_inflow_and_pressure_outflow_R7.5.dat};
        \addplot+[matgreen,   mark phase=31, mark=*, dashed, mark options={solid}] table [y=p_rms, x=theta]{grcbc_farfield_inflow_and_outflow_R7.5.dat};
        \addplot+[matyellow,mark phase=41, mark=*, dashed, mark options={solid}] table [y=p_rms, x=theta]{grcbc_farfield_inflow_and_pressure_outflow_R7.5.dat};
        \addplot+[matorange, mark phase=51, mark=*, dashed, mark options={solid}] table [y=p_rms, x=theta]{nscbc_farfield_inflow_and_pressure_outflow_R7.5.dat};
    \end{polaraxis}
    \node[yshift=-1em] at(current bounding box.south) {\pgfplotslegendfromname{prms_polar_legend}};
\end{tikzpicture}
    \caption{Laminar flow around a cylinder: polar plot of the scaled root-mean-square pressure fluctuation $\Delta \tilde{p}_{\text{rms}}$ for 
    the shifted radii $(a)$ $r_0$ and $(b)$ $r_1$.}
    \label{fig::polar_p_rms}
\end{figure}
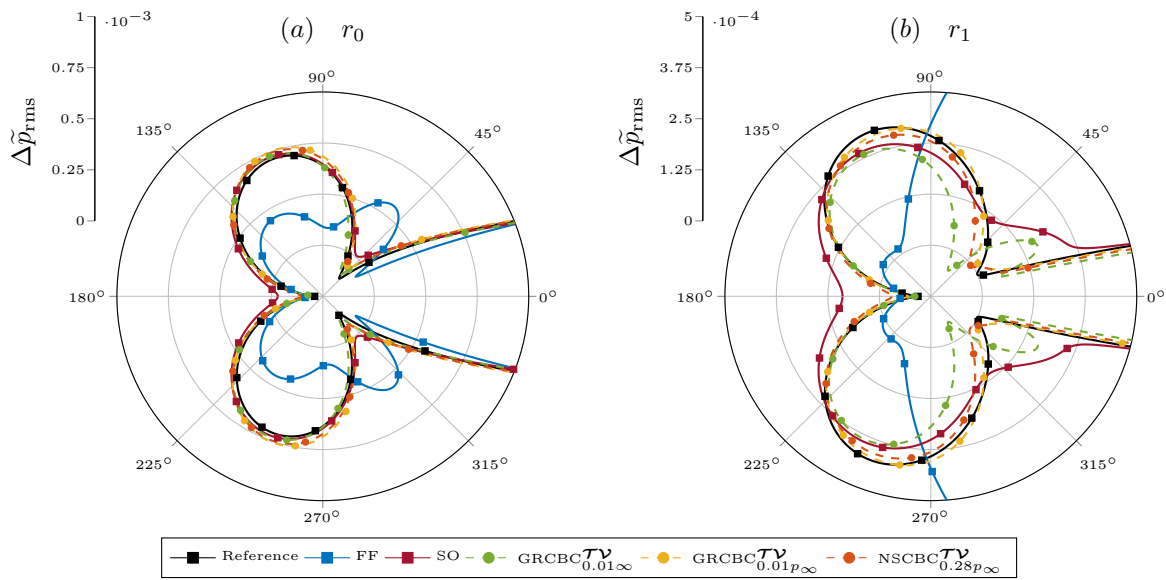

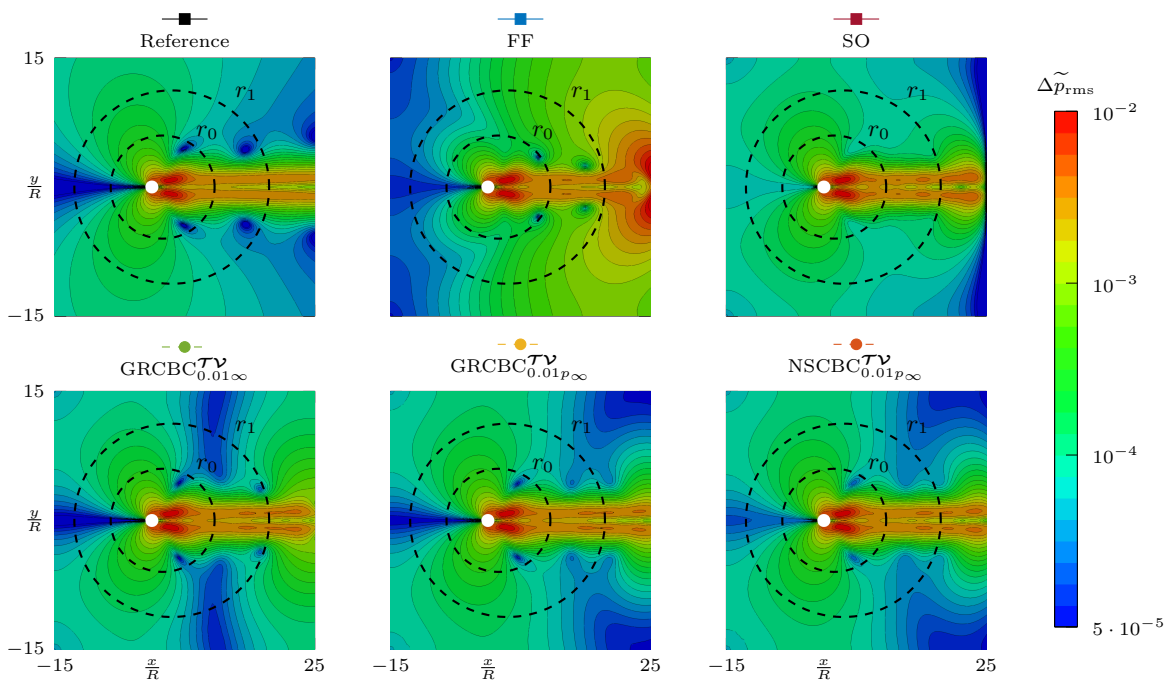
\begin{figure}[h]
    \centering
\begin{tikzpicture}
    \pgfplotsset{every axis/.append style={
    width=5cm,
    height=5cm,
    xmin=-15, xmax=25,
    ymin=-20, ymax=20,
    axis line style={black}, 
    tick style={black}, 
    tick label style = {font=\tiny},
    xlabel style={at={(axis description cs:0.5,-0.1)},anchor=north}, 
    ylabel style={at={(axis descriptioncs:-0.1,.5)},anchor=south}, 
    xtick={-15,0,25}, 
    ytick={-20,0,20},
    yticklabels = {$-15$, $\frac{y}{R}$, $15$},
    xticklabels = {$-15$, $\frac{x}{R}$, $25$}, 
    title style={font=\tiny, align=center, yshift=-0.5em},
    }}
            
    \begin{groupplot}[group style={group size=3 by 2, horizontal sep=1cm, vertical sep=1.0cm}, 
        colorbar to name=cylinder_prms_colorbar, 
        colormap name=netgen,  
        colorbar sampled,
        colorbar style={samples=25, title=$\Delta \tilde{p}_{\text{rms}}$, 
        height=2*\pgfkeysvalueof{/pgfplots/parent axis height}, 
        width=0.3cm, 
        ytick = {0, 10, 20, 30}, 
        yticklabels = {$5\cdot10^{-5}$, $10^{-4}$, $10^{-3}$, $10^{-2}$}, 
        ytick pos=right, 
        scaled y ticks=false, 
        yticklabel style={anchor=west, xshift=0.5em},
        title style={font=\tiny, align=center, yshift=0.5em},}, 
        point meta min = 0, 
        point meta max = 30]

        \nextgroupplot[title={\ref*{leg::cyl::ref} \\ Reference}, xticklabels=\empty, colorbar]
        \addplot graphics[xmin=-15,ymin=-20,xmax=25, ymax=20] {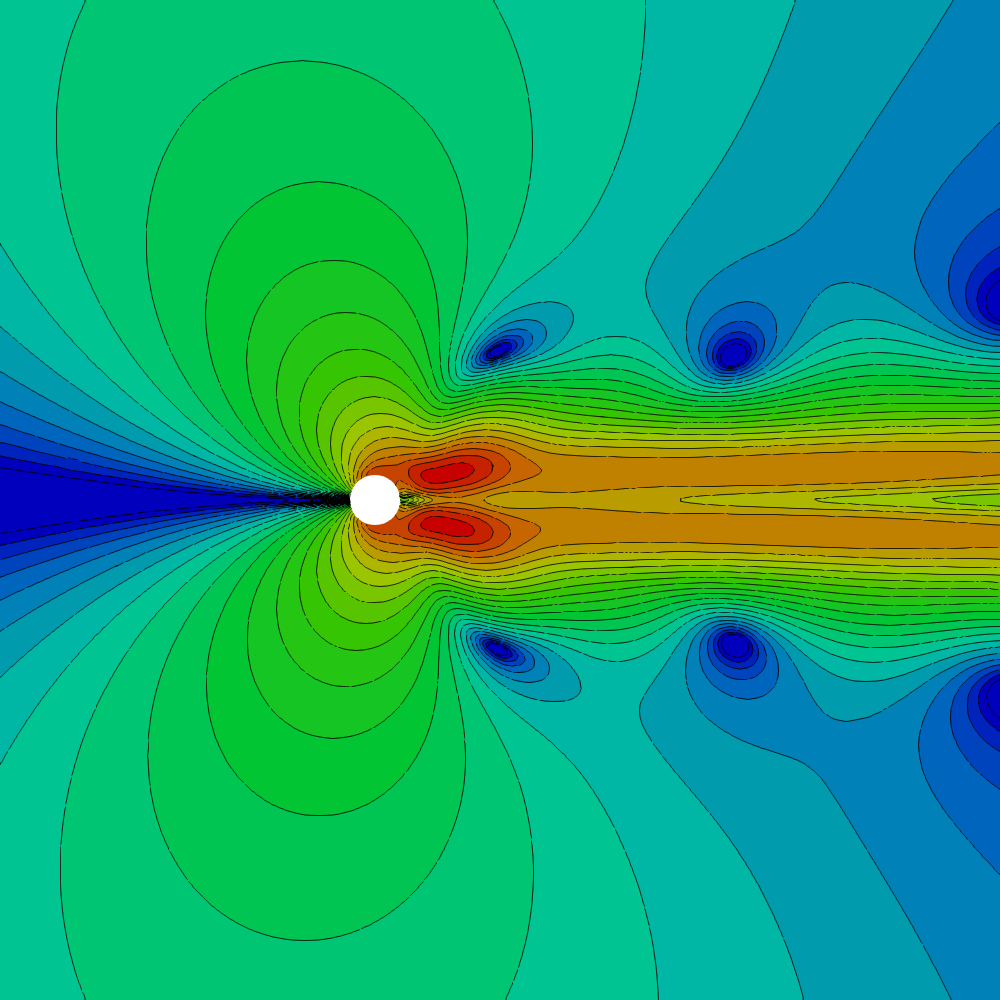};
        \addplot[domain=0:360, samples=100, dashed, thick, forget plot] ({15*cos(x) + 3}, {15*sin(x)})  node[] at (14.5, 14.5){\footnotesize $r_1$};
        \addplot[domain=0:360, samples=100, dashed, thick, forget plot] ({8*cos(x) + 1.6}, {8*sin(x)}) node[] at (8.5, 8.5){\footnotesize $r_0$};

        \nextgroupplot[title={\ref*{leg::cyl::ff} \\ FF}, yticklabels=\empty, xticklabels=\empty]
        \addplot graphics[xmin=-15,ymin=-20,xmax=25,ymax=20] {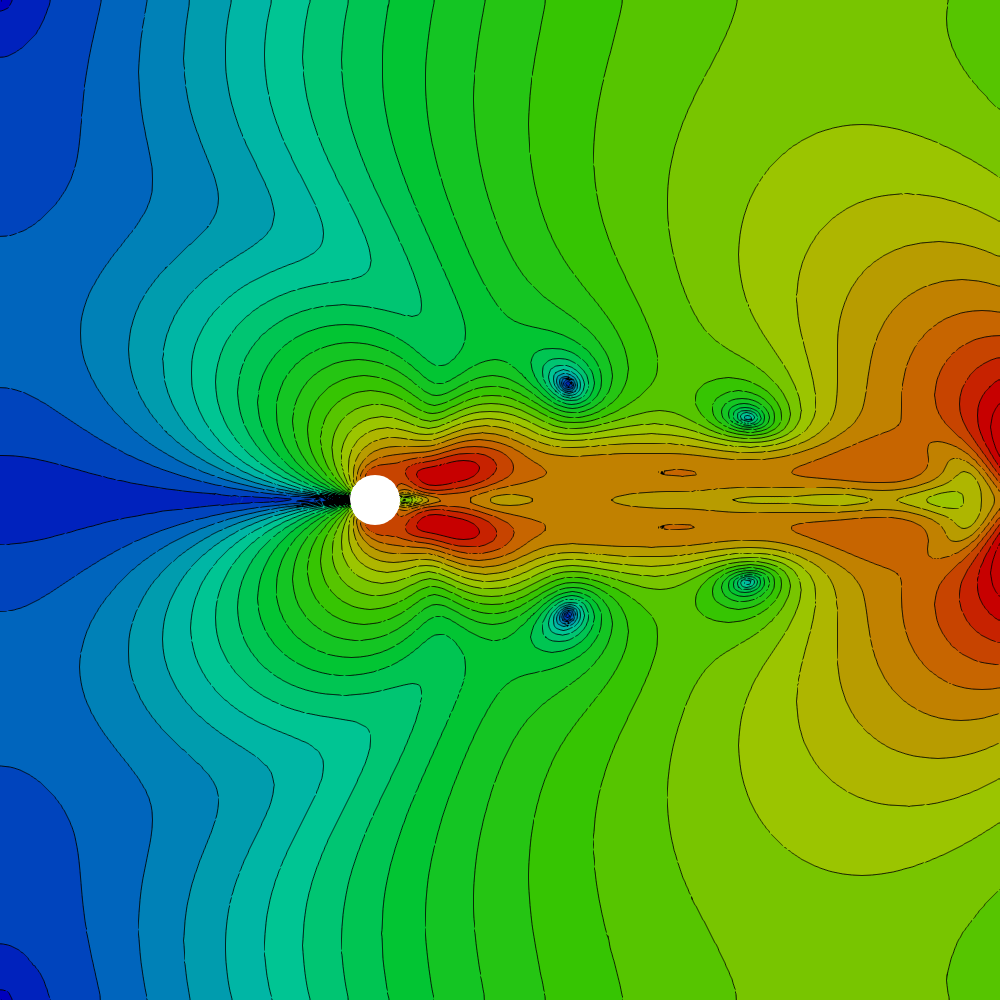};
        \addplot[domain=0:360, samples=100, dashed, thick, forget plot] ({15*cos(x) + 3}, {15*sin(x)})  node[] at (14.5, 14.5){\footnotesize $r_1$};
        \addplot[domain=0:360, samples=100, dashed, thick, forget plot] ({8*cos(x) + 1.6}, {8*sin(x)}) node[] at (8.5, 8.5){\footnotesize $r_0$};

        \nextgroupplot[title={\ref*{leg::cyl::so} \\ SO}, yticklabels=\empty, xticklabels=\empty]
        \addplot graphics[xmin=-15,ymin=-20,xmax=25,ymax=20] {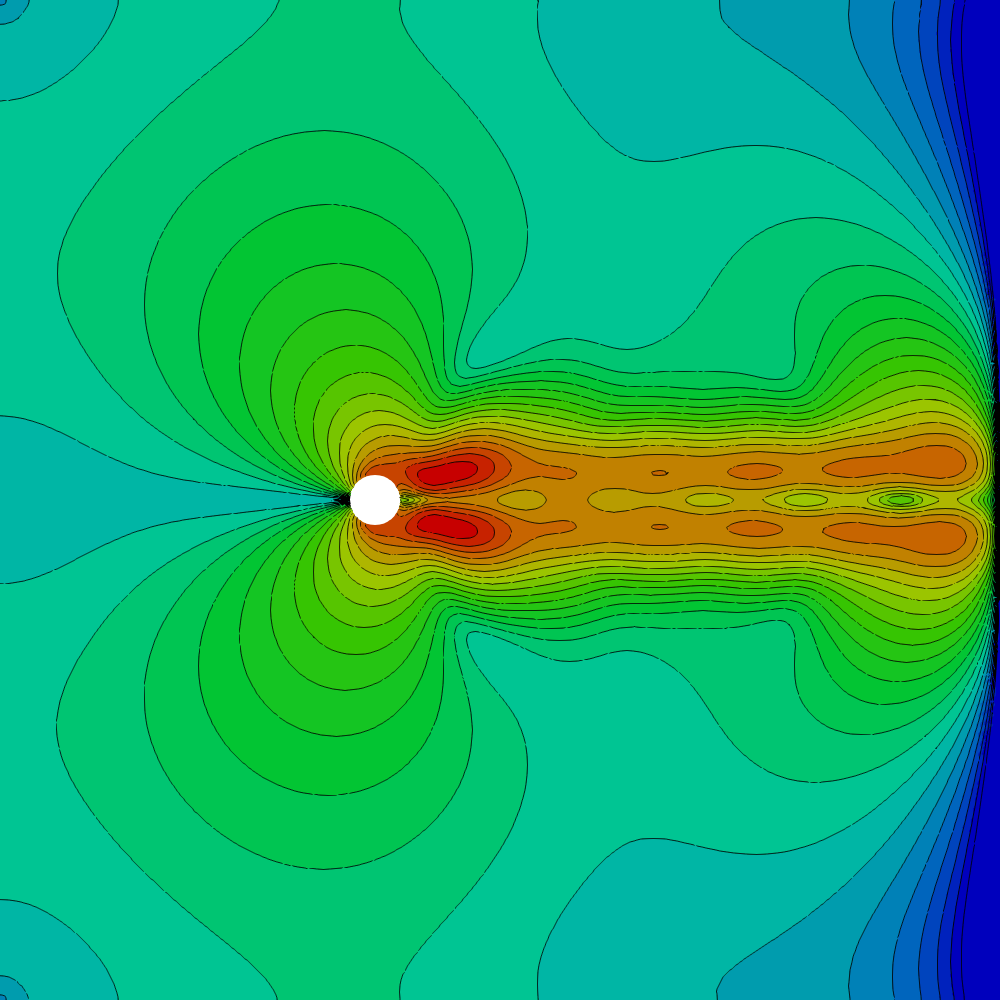};
        \addplot[domain=0:360, samples=100, dashed, thick, forget plot] ({15*cos(x) + 3}, {15*sin(x)})  node[] at (14.5, 14.5){\footnotesize $r_1$};
        \addplot[domain=0:360, samples=100, dashed, thick, forget plot] ({8*cos(x) + 1.6}, {8*sin(x)}) node[] at (8.5, 8.5){\footnotesize $r_0$};

        \coordinate (top) at (rel axis cs:1,1);

        \nextgroupplot[title={\ref*{leg::cyl::grcbc_ff} \\ GRCBC$_{0.01\infty}^{\TRA \VIS}$}]
        \addplot graphics[xmin=-15,ymin=-20,xmax=25,ymax=20] {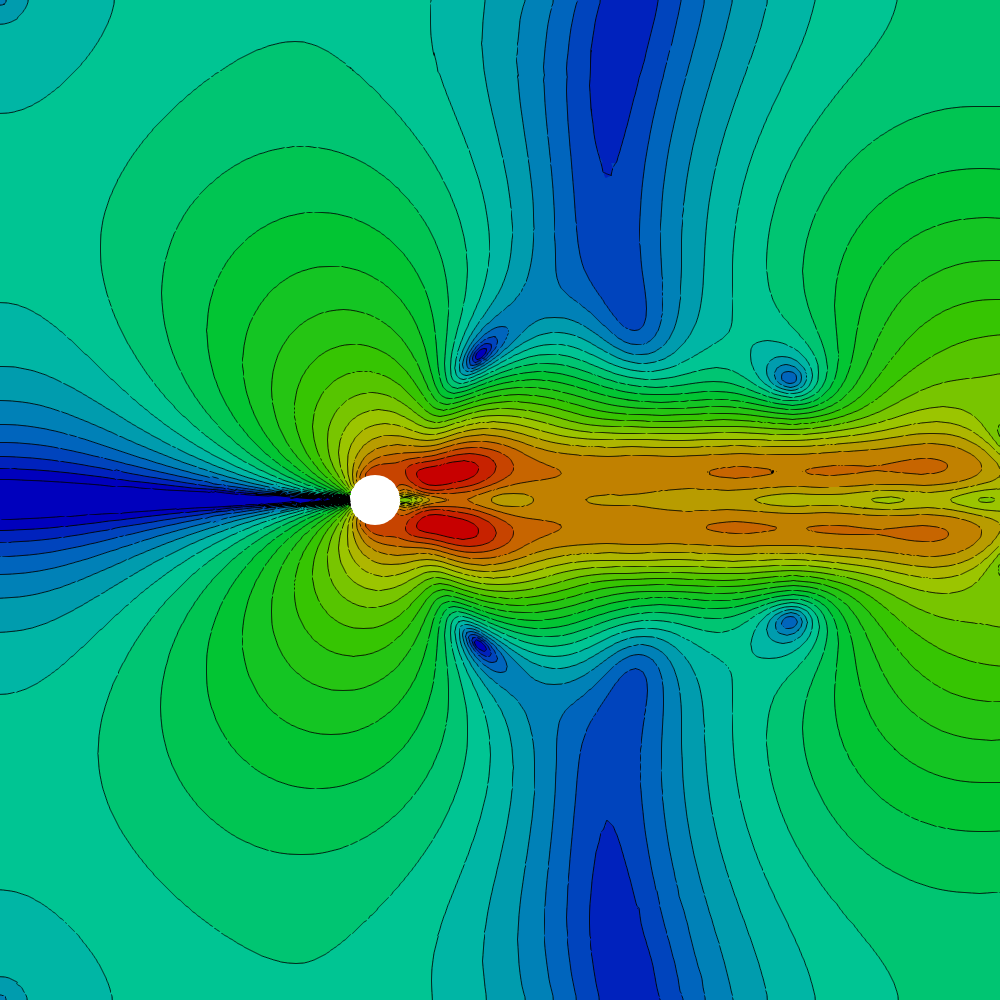};
        \addplot[domain=0:360, samples=100, dashed, thick, forget plot] ({15*cos(x) + 3}, {15*sin(x)})  node[] at (14.5, 14.5){\footnotesize $r_1$};
        \addplot[domain=0:360, samples=100, dashed, thick, forget plot] ({8*cos(x) + 1.6}, {8*sin(x)}) node[] at (8.5, 8.5){\footnotesize $r_0$};

        \nextgroupplot[title={\ref*{leg::cyl::grcbc_so} \\ GRCBC$_{0.01p_\infty}^{\TRA \VIS}$}, yticklabels=\empty]
        \addplot graphics[xmin=-15,ymin=-20,xmax=25,ymax=20] {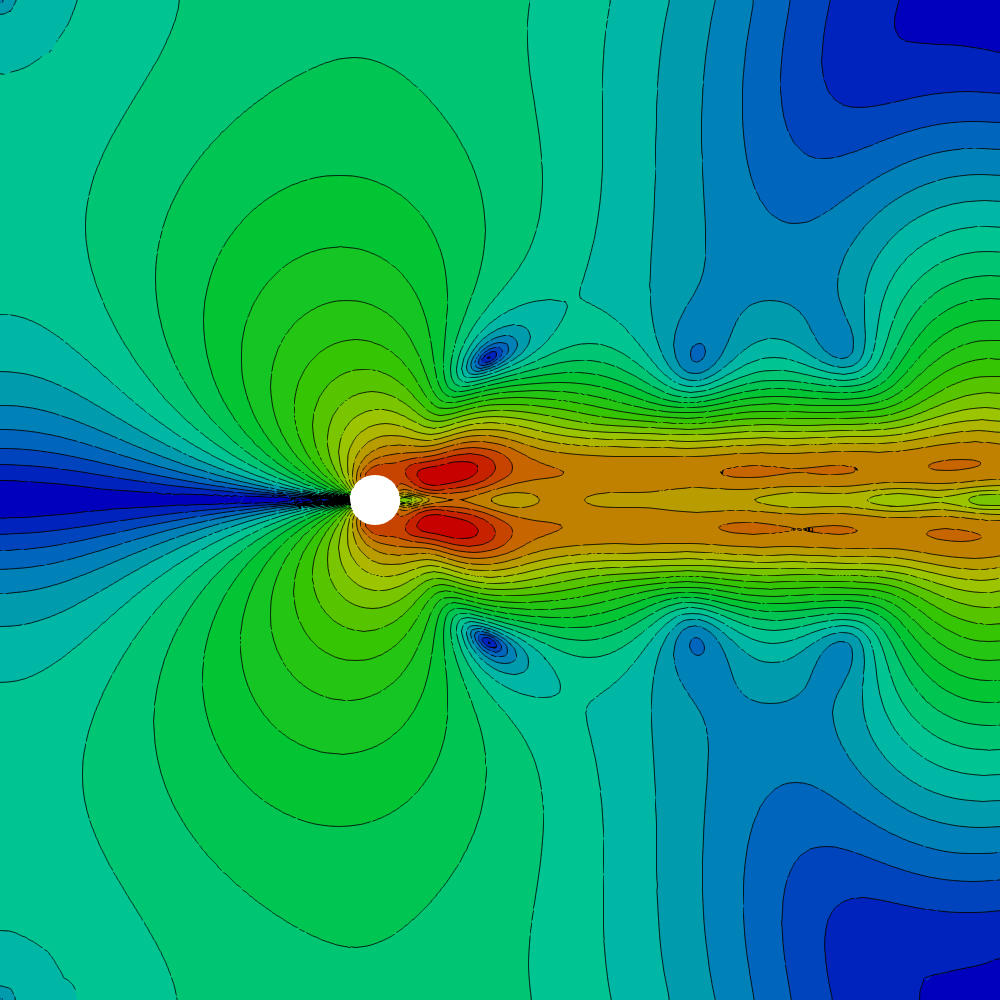};
        \addplot[domain=0:360, samples=100, dashed, thick, forget plot] ({15*cos(x) + 3}, {15*sin(x)})  node[] at (14.5, 14.5){\footnotesize $r_1$};
        \addplot[domain=0:360, samples=100, dashed, thick, forget plot] ({8*cos(x) + 1.6}, {8*sin(x)}) node[] at (8.5, 8.5){\footnotesize $r_0$};

        \nextgroupplot[title={\ref*{leg::cyl::nscbc_so} \\ NSCBC$_{0.01p_\infty}^{\TRA \VIS}$}, yticklabels=\empty]
        \addplot graphics[xmin=-15,ymin=-20,xmax=25, ymax=20] {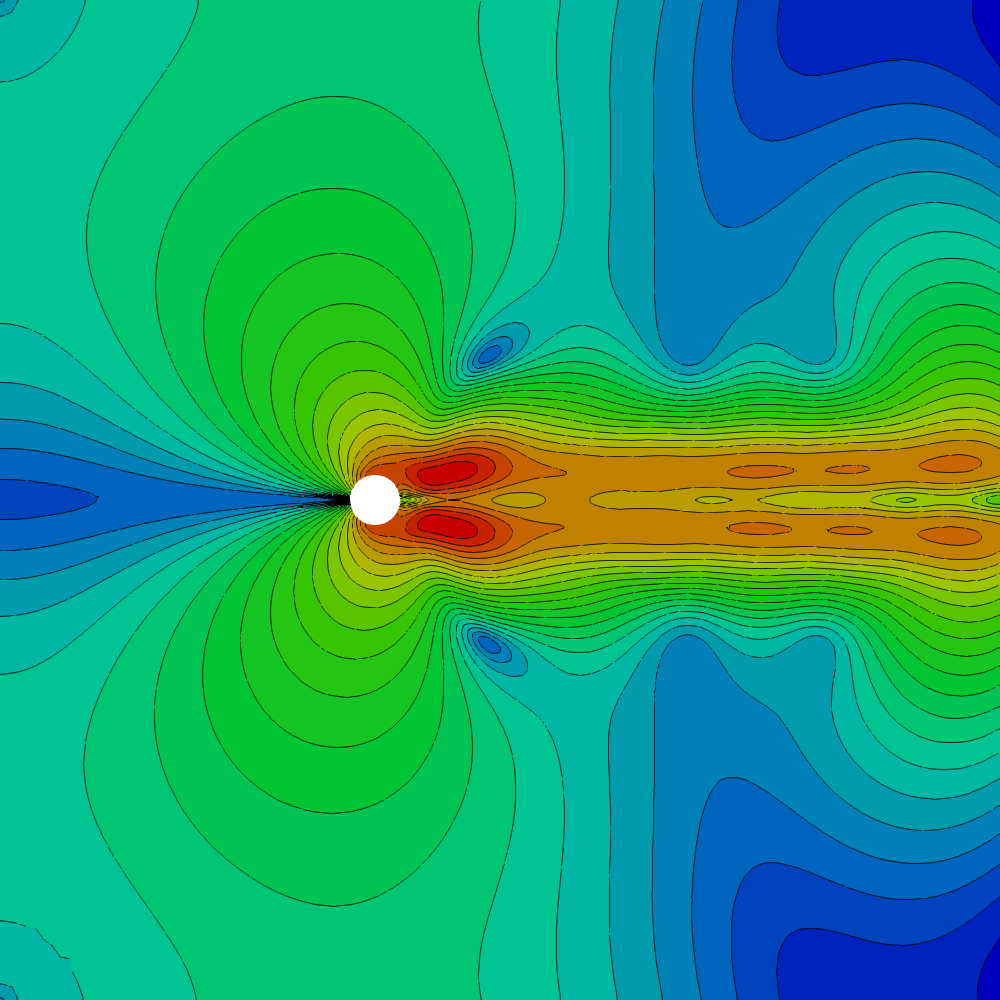};
        \addplot[domain=0:360, samples=100, dashed, thick, forget plot] ({15*cos(x) + 3}, {15*sin(x)})  node[] at (14.5, 14.5){\footnotesize $r_1$};
        \addplot[domain=0:360, samples=100, dashed, thick, forget plot] ({8*cos(x) + 1.6}, {8*sin(x)}) node[] at (8.5, 8.5){\footnotesize $r_0$};

        \coordinate (bottom) at (rel axis cs:1,0);

    \end{groupplot}

    \path (top)--(bottom) coordinate[midway] (group center);
    \node[right=1em,inner sep=0pt] at(group center -| current bounding box.east) {\pgfplotscolorbarfromname{cylinder_prms_colorbar}};

\end{tikzpicture}
    \caption[]{Laminar flow around a cylinder: contour plot of the scaled root-mean-square of the fluctuating pressure $\Delta \tilde{p}_{\text{rms}}$ in logarithmic scale and plot of the shifted circles with radii $r_0$ and $r_1$.}
    \label{fig::cylinder_rms}
\end{figure}

To properly quantify the amount of pressure disturbances in the domain, the scaled root-mean-square values of the fluctuating pressure
\begin{align*}
    \Delta \tilde{p}_{\text{rms}} = \frac{1}{p_\infty} \sqrt{\frac{1}{\overline{t}} \int_{\overline{t}} \Delta \tilde{p}^2 \, dt},
\end{align*}
is computed and compared against the reference solution. In Figure~\ref{fig::polar_p_rms} we show a polar plot of the scaled
root-mean-square of the fluctuating pressure $\Delta \tilde{p}_{\text{rms}}$ at two specific shifted radii
\begin{align*}
    r_0 & = 4.0 (1 + \Ma_{\infty} \cos(\varphi)), & r_1 & = 7.5 (1 + \Ma_{\infty} \cos(\varphi)),
\end{align*}
respectively, where the scaling $(1 + \Ma_{\infty} \cos(\varphi))$ is known as the Doppler factor \cite{inoueSoundGenerationTwodimensional2002}. In Figure~\ref{fig::cylinder_rms},
the contour plot of the scaled root-mean-square of the fluctuating pressure $\Delta \tilde{p}_{\text{rms}}$ is shown in logarithmic scale.

The far-field boundary condition (\ref*{leg::cyl::ff}~FF) shows a vast amount of disturbances propagating upstream,
due to vortices crossing the outflow boundary. Similar to the Zero-Circulation Vortex benchmark, the generalized characteristic boundary condition (\ref*{leg::cyl::grcbc_ff}~GRCBC$_{0.01 \infty}^{\TRA \VIS}$)
demonstrates a significant improvement compared to (\ref*{leg::cyl::ff}~FF), although, in agreement with the results of \cite{pirozzoliGeneralizedCharacteristicRelaxation2013},
performing a relaxation towards the far-field target state $\vec{U}_{\infty}$ at the outflow boundary remains
inappropriate. The subsonic outflow boundary condition (\ref*{leg::cyl::so}~SO) shows a minor deviation from the reference solution upstream of the cylinder $(\varphi=180^\circ)$. This deviation is due to the complete reflection at the subsonic outflow,
a phenomenon absent in settings (\ref*{leg::cyl::grcbc_so}~GRCBC$_{0.01 p_\infty}^{\TRA \VIS}$)  and (\ref*{leg::cyl::nscbc_so}~NSCBC$_{0.28 p_\infty}^{\TRA \VIS}$). 
Specifically as the distance from the cylinder increases towards the boundary, we observe a significant
deviation from the reference solution for boundary conditions (\ref*{leg::cyl::ff}~FF), (\ref*{leg::cyl::so}~SO), and (\ref*{leg::cyl::grcbc_ff}~GRCBC$_{0.01 \infty}^{\TRA \VIS}$).
In contrast, the boundary conditions (\ref*{leg::cyl::grcbc_so}~GRCBC$_{0.01 p_\infty}^{\TRA \VIS}$) and (\ref*{leg::cyl::nscbc_so}~NSCBC$_{0.28 p_\infty}^{\TRA \VIS}$)
demonstrate better agreement due to the inclusion of the transverse $\TRA$ and viscous $\VIS$ terms, particularly at the shifted radius $r_1$.

Lastly, we aim to compute aerodynamic coefficients per unit span, such as the drag coefficient $c_d$ and lift coefficient $c_l$, 
given their relevance and importance in aerodynamic applications. These quantities are computed via the following expressions
\begin{align*}
    \vec{F}_{cyl}  := \int_{\Gamma_{cyl}} \left( \TAU - p \I \right) \vec{n} \, d\bm{s},              \qquad
    c_d            := \frac{1}{\frac{1}{2}\rho_\infty u_\infty^2 D} \vec{F}_{cyl} \cdot \begin{pmatrix}
                                                                                 1 \\ 0
                                                                             \end{pmatrix}, \qquad
    c_l            := \frac{1}{\frac{1}{2} \rho_\infty u_\infty^2 D} \vec{F}_{cyl} \cdot \begin{pmatrix}
                                                                                 0 \\ 1
                                                                             \end{pmatrix}.
\end{align*}
Specifically, we are interested in their time-averaged values $\overline{c}_d, \overline{c}_l$ and first harmonic amplitudes $\tilde{c}_d, \tilde{c}_l$, denoted by
\begin{align*}
    \overline{c}_d  := \frac{1}{\overline{t}} \int_0^{\overline{t}} c_d \, dt, \qquad
    \overline{c}_l  := \frac{1}{\overline{t}} \int_0^{\overline{t}} c_l \, dt, \qquad
    \tilde{c}_d     := \max (c_d - \overline{c}_d ),                           \qquad
    \tilde{c}_l     := \max (c_l - \overline{c}_l ),
\end{align*}
respectively. In Table~\ref{tab::aerodynamic_coefficients}, we show the time-averaged drag coefficient $\overline{c}_d$, the first harmonic amplitude of 
the drag coefficient $\tilde{c}_d$ and lift coefficient $\tilde{c}_l$, along with the determined Strouhal number $\St$ for different applied outflow boundary conditions.
We omit the time-averaged lift coefficient $\overline{c}_l$, as it is zero.
\begin{table}[h]
    \centering
    \begin{tabularx}{\textwidth}{LCCCC}
                                                        & $\overline{c}_d$ & $\tilde{c}_d $ & $\tilde{c}_l$ & $\St$    \\
        \midrule
        \ref*{leg::cyl::ref}~Reference                                  & $1.333$          & $0.02562$      & $0.5177$      & $0.1822$ \\\midrule
        \ref*{leg::cyl::ff}~FF                                          & $1.404$          & $0.02688$      & $0.5403$      & $0.1883$ \\
        \ref*{leg::cyl::so}~SO                                          & $1.380$          & $0.02847$      & $0.5106$      & $0.1855$ \\
        \ref*{leg::cyl::grcbc_ff}~GRCBC$_{0.01 \infty}^{\TRA \VIS}$   & $1.402$          & $0.02706$      & $0.5461$      & $0.1883$ \\
        \ref*{leg::cyl::grcbc_so}~GRCBC$_{0.01 p_\infty}^{\TRA \VIS}$ & $1.392$          & $0.02898$      & $0.5446$      & $0.1872$ \\
        \ref*{leg::cyl::nscbc_so}~NSCBC$_{0.28 p_\infty}^{\TRA \VIS}$ & $1.388$          & $0.02754$      & $0.5265$      & $0.1864$ \\ 
    \end{tabularx}
    \caption{Laminar flow around a cylinder: aerodynamic coefficients.}
    \label{tab::aerodynamic_coefficients}
\end{table}
Compared to the (\ref*{leg::cyl::ref} Reference) solution every simulation clearly overestimates the drag coefficients $\overline{c}_d, \tilde{c}_d$,
which can be attributed to the insufficient downstream placement of the outflow boundary relative to the cylinder. For the lift coefficient $\tilde{c}_l$ and
the Strouhal number $\St$, simulations (\ref*{leg::cyl::so} SO), (\ref*{leg::cyl::grcbc_so} GRCBC$_{0.01 p_{\infty}}^{\TRA \VIS}$) and (\ref*{leg::cyl::nscbc_so} NSCBC$_{0.28 {p_\infty}}^{\TRA \VIS}$)
with the subsonic outflow target state $\vec{U}_{p_\infty}$ show slightly better results than the far-field target state $\vec{U}_{\infty}$ in simulations (\ref*{leg::cyl::ff}~FF) and (\ref*{leg::cyl::grcbc_ff}~GRCBC$_{0.01 \infty}^{\TRA \VIS}$).

We conclude that the newly implemented CBCs significantly reduce pressure disturbances across the domain compared to 
the commonly used subsonic outflow (\ref*{leg::cyl::so} SO) and far-field (\ref*{leg::cyl::ff} FF) boundary conditions, while maintaining comparable accuracy of the computed aerodynamic coefficients.
\section{Conclusion}
\label{sec::conclusion}

We incorporated the concept of characteristic boundary conditions (CBCs) into the Hybridizable Discontinuous Galerkin (HDG) framework 
and proposed a new approach to generalized characteristic relaxation boundary conditions (GRCBCs). 

The latter include commonly used boundary conditions, such as the subsonic outflow \eqref{eq::subsonic_outflow} and the far-field \eqref{eq::farfield} boundary condition as special cases, 
and gives the user the flexibility to tune the relaxation factors to the specific problem under analysis. 
Most importantly, one can further include tangential and viscous contributions, which show to have a significant impact on the computations in the presented numerical experiments.

The implementation of the CBCs integrates seamlessly into the numerical scheme without the need of additional mirror elements
or nodal shape functions as in previous publications in the Discontinuous Galerkin (DG) framework
\cite{toulopoulosArtificialBoundaryConditions2011, shehadiPolynomialcorrectionNavierStokesCharacteristic2024a, shehadiNonReflectingBoundaryConditions2024}.

Through a series of numerical experiments in a weakly compressible setting, we made the following observations:
\begin{itemize}
    \item For a purely one-dimensional acoustic perturbation (see Subsection~\ref{sec::experiments::planar_acoustic_pulse}), 
          the far-field target state $\vec{U}_\infty$ acts entirely non-reflecting, regardless of the level of relaxation in the CBCs, making it the most effective target state in this setting. 
          Conversely, fully reflecting target states, such as $\vec{U}_{p_\infty}, \vec{U}_{\rho u_\infty}, \vec{U}_{\theta_\infty}$, can be tuned from partial to zero reflectivity by appropriately adjusting the relaxation parameters. 
    \item For a two-dimensional non-linear pressure pulse injected in a fluid at rest (see Subsection~\ref{sec::experiments::oblique_pressure_pulse}), the best result was obtained
          by combining the NSCBC with the subsonic outflow target state $\vec{U}_{p_\infty}$ and the optimal determined \cite{rudyNonreflectingOutflowBoundary1980} relaxation parameter $\sigma = 0.28$ .
          However, while the CBCs using the far-field target state $\vec{U}_\infty$ tend to slightly overdamp the perturbation, likely due to the absence of transverse terms,
          the outcome is comparable to the NSCBC one. Moreover, they are less sensitive to both the relaxation parameter and the alternating inflow and outflow caused by the pulse crossing the boundary.
    \item For a common scenario in aeroacoustic applications, where vortices cross an artificial outflow boundary (see Subsection~\ref{sec::experiments::zero_circulation_vortex} and Subsection~\ref{sec::experiments::cylinder_flow}),
          the relaxation towards the far-field target state $\vec{U}_\infty$ is inappropriate, in good agreement with \cite{pirozzoliGeneralizedCharacteristicRelaxation2013}.
          Instead, the CBCs, combined with the subsonic outflow target state $\vec{U}_{p_\infty}$ and the inclusion of transverse terms $\TRA$ and viscous terms $\VIS$, show superior performance compared to conventional boundary conditions, both in inviscid and viscous flow scenarios for weakly compressible flow, by minimizing the reflection of vortices at artificial boundaries.
\end{itemize}
In summary, while the implemented CBCs are not a universal solution and may need to be complemented by other techniques,
such as sponge layers or perfectly matched layers, in cases where reflectivity poses a challenge, we conclude that they serve as an
effective tool for minimizing reflectivity at artificial boundaries in weakly compressible flow simulations. 
Particularly, at artificial subsonic outflow boundaries in aeroacoustic simulations, the novel GRCBCs show 
improved performance compared to common HDG boundary conditions in reducing the reflection of acoustic waves and vortices.
As such, they should be preferred over the latter in such scenarios.

\section{Acknowledgments}
\label{sec::acknowledgments}
We would like to thank the Austrian Science Fund (FWF) and the Vienna Scientific Cluster (VSC) for their support.
This research was funded in whole or in part by the Austrian Science Fund (FWF) [10.55776/P35931]. 
Calculations were performed using supercomputer resources provided by the Vienna Scientific Cluster (VSC). 
M.G. and A.H. acknowledge the Spanish Ministry of Science, Innovation and Universities and the 
Spanish State Research Agency MICIU/AEI/10.13039/501100011033 (Grant agreement No. PID2023-149979OB-I00) and the
Generalitat de Catalunya (Grant agreement No. 2021-SGR-01049). M.G. is Fellow of the Serra Húnter Programme of the Generalitat de Catalunya.

\bibliography{literature}
\bibliographystyle{unsrtnat}

\end{document}